\documentclass[11pt]{article}
\usepackage{makeidx}
\usepackage{amssymb,amsfonts,amsmath,graphicx,hyperref}

\newcommand{\mcm}[3]{\newcommand{#1}[#2]{{\ensuremath{#3}}}}

\usepackage{latexsym}
\usepackage{amssymb}

\mcm{\blank}{0}{(\emptybk)} \mcm{\dashbk}{0}{\mbox{---}}
\mcm{\emptybk}{0}{\:\:} \mcm{\hyph}{0}{\mbox{-}}

\mcm{\diagspace}{0}{\mbox{\hspace{2em}}}

\mcm{\cat}{1}{\mc{#1}} \mcm{\fcat}{1}{\mb{#1}}
\mcm{\mc}{1}{\mathcal{#1}} \mcm{\mr}{1}{\mathrm{#1}}
\mcm{\mi}{1}{\mathit{#1}} \mcm{\mb}{1}{\mathbf{#1}}
\mcm{\scat}{1}{\Bbb{#1}} \mcm{\twid}{1}{\widetilde{#1}}

\mcm{\elt}{0}{\in} \mcm{\sub}{0}{\,\subseteq\,}
\mcm{\such}{0}{\:|\:} \mcm{\without}{0}{\setminus}

\mcm{\atsr}{0}{\Box} \mcm{\eqv}{0}{\,\simeq\,}
\mcm{\iso}{0}{\,\cong\,}
\mcm{\of}{0}{\raisebox{0.2mm}{\ensuremath{\scriptstyle\circ}}}

\mcm{\bdry}{0}{\partial}

\mcm{\Bee}{0}{\cat{B}} \mcm{\Beep}{0}{\cat{B'}}
\mcm{\Eee}{0}{\cat{E}} \mcm{\Eeep}{0}{\cat{E'}}

\mcm{\Ess}{0}{\cat{S}} \mcm{\Tee}{0}{\cat{T}}
\mcm{\Teep}{0}{\cat{T'}} \mcm{\Stee}{0}{\scat{T}}
\mcm{\Steep}{0}{\scat{T'}}

\mcm{\blbk}{0}{\blank^{\blob}}
\mcm{\blob}{0}{\scriptscriptstyle{\bullet}}
\mcm{\stbk}{0}{\blank^{*}} \mcm{\ubl}{0}{{}^{\blob}}
\mcm{\ust}{0}{{}^{*}}

\mcm{\Cartpr}{0}{\pr{\Eee}{T}} \mcm{\Cartprp}{0}{\pr{\Eeep}{T'}}
\mcm{\Mnd}{0}{\triple{T}{\eta}{\mu}}
\mcm{\Zeropr}{0}{\pr{\Set}{\id}}

\mcm{\dopset}{0}{\ftrcat{\Delta^{\op}}{\Set}}
\mcm{\tropset}{0}{\ftrcat{\fcat{TR}^{\op}}{\Set}}

\mcm{\cod}{0}{\mr{cod}} \mcm{\dom}{0}{\mr{dom}}
\mcm{\End}{0}{\mr{End}} \mcm{\Hom}{0}{\mr{Hom}}
\mcm{\ob}{0}{\mr{ob}\,} \mcm{\op}{0}{\mr{op}}

\mcm{\comp}{0}{\mi{comp}} \mcm{\id}{0}{\mi{id}}
\mcm{\ids}{0}{\mi{ids}} \mcm{\mult}{0}{\mi{mult}}
\mcm{\unit}{0}{\mi{unit}}

\mcm{\Ab}{0}{\fcat{Ab}} \mcm{\Alg}{0}{\fcat{Alg}}
\mcm{\Bim}{1}{\fcat{Bim}(#1)} \mcm{\Cat}{0}{\fcat{Cat}}
\mcm{\Cay}{0}{\fcat{Cay}} \mcm{\Cpn}{1}{\pr{\Set/S_{#1}}{T_{#1}}}
\mcm{\fc}{0}{\fcat{fc}} \mcm{\fm}{0}{\fcat{fm}}
\mcm{\Graph}{0}{\fcat{Graph}} \mcm{\Gy}{0}{\fcat{Gy}}
\mcm{\Hpn}{1}{\pr{\Eee_{#1}}{P_{#1}}} \mcm{\Mon}{0}{\mb{Mon}}
\mcm{\Multicat}{0}{\fcat{Multicat}} \mcm{\One}{0}{\fcat{1}}
\mcm{\PD}{1}{\fcat{PD}_{#1}} \mcm{\Prof}{0}{\fcat{Prof}}
\mcm{\Set}{0}{\fcat{Set}} \mcm{\Span}{0}{\fcat{Span}}
\mcm{\Ssq}{0}{\fcat{Ssq}} \mcm{\Struc}{0}{\fcat{Struc}}
\mcm{\Sym}{0}{\fcat{Sym}} \mcm{\TR}{1}{\fcat{TR}(#1)}
\mcm{\Tr}{0}{\fcat{Tr}} \mcm{\Twocat}{0}{\fcat{2\hyph\Cat}}

\mcm{\integers}{0}{\mathbb{Z}}

\mcm{\range}{2}{#1,\,\ldots\,,#2}
\mcm{\bftuple}{2}{\tuplebts{\range{#1}{#2}}}
\mcm{\tuple}{3}{\tuplebts{\range{#1,#2}{#3}}}
\mcm{\rttuple}{1}{\tuplebts{\,\ldots\,,#1}}
\mcm{\abftuple}{2}{\atuplebts{\range{#1}{#2}}}
\mcm{\atuple}{3}{\atuplebts{\range{#1,#2}{#3}}}
\mcm{\arttuple}{1}{\atuplebts{\,\ldots\,,#1}}
\mcm{\sqbftuple}{2}{\obt\range{#1}{#2}\cbt}
\mcm{\pr}{2}{\tuplebts{#1,#2}}
\mcm{\triple}{3}{\tuplebts{#1,#2,#3}}

\mcm{\eend}{2}{#1[#2]} \mcm{\ehom}{3}{#1[#2,#3]}
\mcm{\ftrcat}{2}{[#1,#2]} \mcm{\homset}{3}{#1(#2,#3)}
\mcm{\multihom}{3}{#1(#2;#3)}
\mcm{\relhom}{5}{#1_{#2}(\range{#3}{#4};#5)}

\mcm{\go}{0}{\rTo} \mcm{\goby}{1}{\rTo^{#1}}
\mcm{\goesto}{0}{\,\longmapsto\,} \mcm{\goiso}{0}{\goby{\diso}}
\mcm{\monic}{0}{\rMonic} \mcm{\og}{0}{\lTo}
\mcm{\ogby}{1}{\lTo^{#1}}

\mcm{\gph}{2}{\spn{#1}{T #2}{#2}} \mcm{\graph}{4}{\spaan{#1}{T
#2}{#2}{#3}{#4}} \mcm{\oppair}{2}{\stackrel{\rTo^{#1}}{\lTo_{#2}}}
\mcm{\parpair}{2}{\stackrel{\rTo^{#1}}{\rTo_{#2}}}
\mcm{\spn}{3}{#2 \og #1 \go #3} \mcm{\spaan}{5}{#2 \ogby{#4} #1
\goby{#5} #3}

\mcm{\bktdvslob}{3}
    {\left(
    \begin{diagram}[height=1.5em]
    #1      \\
    \dTo>{\,#2} \\
    #3      \\
    \end{diagram}
    \right)}
\mcm{\slob}{3}{(#1 \goby{#2} #3)} \mcm{\vslob}{3}
    {\left.
    \begin{diagram}[height=1.5em]
    #1      \\
    \dTo>{\,#2} \\
    #3      \\
    \end{diagram}
    \right.}

\newenvironment{tree}
    {\begin{diagram}[height=1em,width=.75em,abut,noPS,tight]}
    {\end{diagram}}

\mcm{\enode}{0}{\circ}

\mcm{\nl}{1}{\stackrel{\textstyle #1}{\node}}
\mcm{\node}{0}{\bullet}

\mcm{\utree}{0}{\node}

\mcm{\diso}{0}{\sim}

\mcm{\vdiso}{0}{\wr}

\mcm{\nat}{0}{\mathbb{N}}

\mcm{\Onepr}{0}{\pr{\Graph}{\fc}}
\newlength{\nllwidth}
\newlength{\nllheight}
\newcommand{\stackbelow}[2]{%
\settowidth{\nllwidth}{\ensuremath{#1}\ensuremath{#2}}%
\settoheight{\nllheight}{\ensuremath{#2}}%
\addtolength{\nllheight}{.3ex}%
\mbox{%
\ensuremath{#1}%
\hspace{-.5\nllwidth}%
\raisebox{-1\nllheight}{\ensuremath{#2}}}}
\mcm{\nlal}{2}{\stackbelow{\nl{#1}}{#2}}
\mcm{\nll}{1}{\stackbelow{\node}{#1}} \mcm{\wun}{0}{\fcat{1}}
\mcm{\atuplebts}{1}{\langle #1 \rangle} \mcm{\tuplebts}{1}{(#1)}
\mcm{\bo}{0}{(} \mcm{\bc}{0}{)}
\mcm{\UBilax}{0}{\fcat{UBicat}_\mr{lax}}
\mcm{\UBiwk}{0}{\fcat{UBicat}_\mr{wk}}
\mcm{\UBistr}{0}{\fcat{UBicat}_\mr{str}}
\mcm{\Bilax}{0}{\fcat{Bicat}_\mr{lax}}
\mcm{\Biwk}{0}{\fcat{Bicat}_\mr{wk}}
\mcm{\Bistr}{0}{\fcat{Bicat}_\mr{str}} \mcm{\rotsub}{0}{\cup
\raisebox{0.1em}{$\scriptstyle{|}$}} \mcm{\pd}{0}{\fcat{pd}}
\mcm{\rep}{1}{\widehat{#1}} \mcm{\ovln}{1}{\overline{#1}}
\mcm{\Gph}{0}{\fcat{Gph}} \mcm{\tr}{0}{\fcat{tr}}

\mcm{\ladj}{0}{\,\dashv\,} \mcm{\zeropd}{0}{\node}
    {\end{diagram}}
\mcm{\END}{0}{\fcat{End}} \mcm{\HOM}{0}{\fcat{Hom}}

\newlength{\gwidth} 
\newlength{\gvert}  
\newlength{\gdrop}  
\newlength{\gbaredrop}  
\newlength{\goffset}    
\newlength{\gtemp}  

\newcommand{\present}[1]{%
\makebox[1\gwidth]{%
\rule[-1\gdrop]{0ex}{1\gvert}%
\raisebox{-1\gbaredrop}{#1}}}

\newcommand{\presentl}[1]{%
\makebox[1\gwidth][l]{%
\rule[-1\gdrop]{0ex}{1\gvert}%
\raisebox{-1\gbaredrop}{#1}}}

\newcommand{\presentr}[1]{%
\makebox[1\gwidth][r]{%
\rule[-1\gdrop]{0ex}{1\gvert}%
\raisebox{-1\gbaredrop}{#1}}}

\newcommand{\ginitdims}[2]{
\setlength{\unitlength}{1em}
\setlength{\goffset}{.25\unitlength}
\setlength{\gwidth}{#1\unitlength}
\setlength{\gvert}{#2\unitlength}
\setlength{\gdrop}{.5\gvert}
\addtolength{\gdrop}{-1\goffset}
\setlength{\gbaredrop}{1\gdrop}
\addtolength{\gvert}{.6\unitlength}
\addtolength{\gdrop}{.3\unitlength}}    

\newcommand{\cinitdims}[2]{
\setlength{\unitlength}{1em}
\setlength{\goffset}{.35\unitlength}
\setlength{\gwidth}{#1\unitlength}
\setlength{\gvert}{#2\unitlength}
\setlength{\gdrop}{.5\gvert}
\addtolength{\gdrop}{-1\goffset}
\setlength{\gbaredrop}{1\gdrop}
\addtolength{\gvert}{.6\unitlength}
\addtolength{\gdrop}{.3\unitlength}}    

\newcommand{\gsinitdims}[2]{
\setlength{\unitlength}{0.5em}
\setlength{\goffset}{.25\unitlength}
\setlength{\gwidth}{#1\unitlength}
\setlength{\gvert}{#2\unitlength}
\setlength{\gdrop}{.5\gvert}
\addtolength{\gdrop}{-1\goffset}
\setlength{\gbaredrop}{1\gdrop}
\addtolength{\gvert}{.6\unitlength}
\addtolength{\gdrop}{.3\unitlength}}    

\newcommand{\sidespic}[1]{%
\settowidth{\gtemp}{\ensuremath{#1}}%
\addtolength{\gwidth}{1\gtemp}}

\newcommand{\abovepic}[1]{%
\settoheight{\gtemp}{\ensuremath{#1}}%
\addtolength{\gvert}{1\gtemp}%
\settodepth{\gtemp}{\ensuremath{#1}}%
\addtolength{\gvert}{1\gtemp}}

\newcommand{\belowpic}[1]{%
\settoheight{\gtemp}{\ensuremath{#1}}%
\addtolength{\gvert}{1\gtemp}%
\addtolength{\gdrop}{1\gtemp}%
\settodepth{\gtemp}{\ensuremath{#1}}%
\addtolength{\gvert}{1\gtemp}%
\addtolength{\gdrop}{1\gtemp}}

\newcommand{\cell}[4]{\put(#1,#2){\makebox(0,0)[#3]{\ensuremath{#4}}}}
\mcm{\zmark}{0}{\scriptstyle{\bullet}}

\newcommand{\pregfst}[1]{%
\begin{picture}(0.5,0.2)(-0.5,-0.2)%
\cell{-0.1}{-0.2}{tr}{#1}%
\cell{0}{0}{c}{\zmark}%
\end{picture}}

\mcm{\gfst}{1}{%
\ginitdims{0.5}{0.4}%
\sidespic{#1}%
\belowpic{#1}%
\presentr{\pregfst{#1}}}

\newcommand{\preglst}[1]{%
\begin{picture}(0.5,0.2)(0,-0.2)%
\cell{0.1}{-0.2}{tl}{#1}%
\cell{0.05}{0}{c}{\zmark}%
\end{picture}}

\mcm{\glst}{1}{%
\ginitdims{.5}{.4}%
\sidespic{#1}%
\belowpic{#1}%
\presentl{\preglst{#1}}}

\newcommand{\preglft}[1]{%
\begin{picture}(0,0.2)(0,-0.2)%
\cell{-0.1}{-0.2}{tr}{#1}%
\cell{0.05}{0}{c}{\zmark}%
\end{picture}}

\mcm{\glft}{1}{%
\ginitdims{0}{.4}%
\belowpic{#1}%
\present{\preglft{#1}}}

\newcommand{\pregrgt}[1]{%
\begin{picture}(0,0.2)(0,-0.2)%
\cell{0.1}{-0.2}{tl}{#1}%
\cell{0.05}{0}{c}{\zmark}%
\end{picture}}

\mcm{\grgt}{1}{%
\ginitdims{0}{.4}%
\belowpic{#1}%
\present{\pregrgt{#1}}}

\newcommand{\pregblw}[1]{%
\begin{picture}(0,0.3)(0,-0.3)
\cell{0}{-0.3}{t}{#1}%
\cell{0.05}{0}{c}{\zmark}%
\end{picture}}

\mcm{\gblw}{1}{%
\ginitdims{0}{.6}%
\belowpic{#1}%
\present{\pregblw{#1}}}

\newcommand{\pregfbw}[1]{%
\begin{picture}(0,0.65)(0,-0.65)
\cell{0}{-0.65}{t}{#1}%
\cell{0.05}{0}{c}{\zmark}%
\end{picture}}

\mcm{\gfbw}{1}{%
\ginitdims{0}{1.3}%
\belowpic{#1}%
\present{\pregfbw{#1}}}

\newcommand{\pregzero}[1]{%
\begin{picture}(0.8,0.4)(-0.4,-0.4)
\cell{0}{-0.4}{t}{#1}%
\cell{0}{0}{c}{\zmark}%
\end{picture}}

\mcm{\gzero}{1}{%
\ginitdims{0.8}{.6}%
\belowpic{#1}%
\sidespic{#1}%
\present{\pregzero{#1}}}

\newcommand{\pregone}[1]{%
\begin{picture}(5,0.4)(0,-0.2)%
\cell{2.5}{0.2}{b}{#1}%
\put(0,0){\vector(1,0){5}}%
\end{picture}}

\mcm{\gone}{1}{%
\ginitdims{5}{0.4}%
\abovepic{#1}%
\present{\pregone{#1}}}

\newcommand{\pregtwo}[3]{%
\begin{picture}(5,3.4)(0,-0.2)%
\cell{2.5}{3.2}{b}{#1}%
\cell{2.5}{-.2}{t}{#2}%
\cell{2.7}{1.5}{l}{#3}%
\qbezier(0,1.5)(2.5,4.5)(5,1.5)%
\qbezier(0,1.5)(2.5,-1.5)(5,1.5)%
\put(5,1.5){\vector(1,-1){0}}%
\put(5,1.5){\vector(1,1){0}}%
\put(2.5,2.5){\vector(0,-1){2}}%
\end{picture}}

\mcm{\gtwo}{3}{%
\ginitdims{5}{3.4}%
\abovepic{#1}%
\belowpic{#2}%
\present{\pregtwo{#1}{#2}{#3}}}

\newcommand{\pregthree}[5]{%
\begin{picture}(5,5.4)(0,-1.2)%
\cell{2.5}{4.2}{b}{#1}%
\cell{1.5}{1.7}{b}{#2}%
\cell{2.5}{-1.2}{t}{#3}%
\cell{2.7}{2.75}{l}{#4}%
\cell{2.7}{0.25}{l}{#5}%
\qbezier(0,1.5)(2.5,6.5)(5,1.5)%
\qbezier(0,1.5)(2.5,-3.5)(5,1.5)%
\put(0,1.5){\vector(1,0){5}}%
\put(2.5,3.5){\vector(0,-1){1.5}}%
\put(2.5,1){\vector(0,-1){1.5}}%
\put(5,1.5){\vector(1,-3){0}}%
\put(5,1.5){\vector(1,3){0}}%
\end{picture}}

\mcm{\gthree}{5}{%
\ginitdims{5}{5.4}%
\abovepic{#1}%
\belowpic{#3}%
\present{\pregthree{#1}{#2}{#3}{#4}{#5}}}

\newcommand{\pregfour}[7]{%
\begin{picture}(5,8.4)(0,-2.7)%
\cell{2.5}{5.7}{b}{#1}%
\cell{1.5}{2.8}{b}{#2}%
\cell{1.5}{0.2}{t}{#3}%
\cell{2.5}{-2.7}{t}{#4}%
\cell{2.7}{4.25}{l}{#5}%
\cell{2.7}{1.5}{l}{#6}%
\cell{2.7}{-1.25}{l}{#7}%
\qbezier(0,1.5)(2.5,9.5)(5,1.5)%
\qbezier(0,1.5)(2.5,4)(5,1.5)%
\qbezier(0,1.5)(2.5,-1)(5,1.5)%
\qbezier(0,1.5)(2.5,-6.5)(5,1.5)%
\put(2.5,5.25){\vector(0,-1){2}}%
\put(2.5,2.5){\vector(0,-1){2}}%
\put(2.5,-0.25){\vector(0,-1){2}}%
\put(5,1.5){\vector(1,-4){0}}%
\put(5,1.5){\vector(4,-3){0}}%
\put(5,1.5){\vector(4,3){0}}%
\put(5,1.5){\vector(1,4){0}}%
\end{picture}}

\mcm{\gfour}{7}{%
\ginitdims{5}{8.4}%
\abovepic{#1}%
\belowpic{#4}%
\present{\pregfour{#1}{#2}{#3}{#4}{#5}{#6}{#7}}}

\newcommand{\pregthreecell}[5]{%
\begin{picture}(8,5)(-4,-2.5)%
\cell{0}{2.5}{b}{#1}%
\cell{0}{-2.5}{t}{#2}%
\cell{-1.7}{0}{r}{#3}%
\cell{1.7}{0}{l}{#4}%
\cell{0}{0.2}{b}{#5}%
\qbezier(-4,0)(0,4.2)(4,0)%
\qbezier(-4,0)(0,-4.2)(4,0)%
\qbezier(-0.5,1.8)(-2.5,0)(-0.5,-1.8)%
\qbezier(0.5,1.8)(2.5,0)(0.5,-1.8)%
\put(-1,0){\vector(1,0){2}}%
\put(4,0){\vector(1,-1){0}}%
\put(4,0){\vector(1,1){0}}%
\put(-0.5,-1.8){\vector(1,-1){0}}%
\put(0.5,-1.8){\vector(-1,-1){0}}%
\end{picture}}

\mcm{\gthreecell}{5}{%
\ginitdims{8}{5}%
\abovepic{#1}%
\belowpic{#2}%
\present{\pregthreecell{#1}{#2}{#3}{#4}{#5}}}

\newcommand{\pregthreecellu}{%
\begin{picture}(5,3.4)(-0.5,-0.2)%
\qbezier(-.5,1.5)(2,4.5)(4.5,1.5)%
\qbezier(-.5,1.5)(2,-1.5)(4.5,1.5)%
\qbezier(1.5,2.7)(0.5,1.5)(1.5,0.3)%
\qbezier(2.5,2.7)(3.5,1.5)(2.5,0.3)%
\put(1.3,1.5){\vector(1,0){1.4}}%
\put(4.5,1.5){\vector(1,-1){0}}%
\put(4.5,1.5){\vector(1,1){0}}%
\put(1.5,0.3){\vector(2,-3){0}}%
\put(2.5,0.3){\vector(-2,-3){0}}%
\end{picture}}

\mcm{\gthreecellu}{0}{%
\ginitdims{5}{3.4}%
\present{\pregthreecellu}}

\newcommand{\pregtwocentre}[3]{%
\begin{picture}(5,3.4)(0,-0.2)%
\cell{2.5}{3.2}{b}{#1}%
\cell{2.5}{-.2}{t}{#2}%
\cell{2.5}{1.5}{c}{#3}%
\qbezier(0,1.5)(2.5,4.5)(5,1.5)%
\qbezier(0,1.5)(2.5,-1.5)(5,1.5)%
\put(5,1.5){\vector(1,-1){0}}%
\put(5,1.5){\vector(1,1){0}}%
\put(2.5,2.5){\vector(0,-1){2}}%
\end{picture}}

\mcm{\gtwocentre}{3}{%
\ginitdims{5}{3.4}%
\abovepic{#1}%
\belowpic{#2}%
\present{\pregtwocentre{#1}{#2}{#3}}}

\newcommand{\pregspecialone}[9]{%
\begin{picture}(8,8)(-4,-4)%
\cell{0}{3.9}{b}{#1}%
\cell{-2}{-0.2}{t}{#2}%
\cell{0}{-3.9}{t}{#3}%
\cell{-1.5}{1.1}{r}{#4}%
\cell{0.2}{1.5}{l}{#5}%
\cell{1.5}{1.1}{l}{#6}%
\cell{0.2}{-2}{l}{#7}%
\cell{-0.9}{2.3}{b}{#8}%
\cell{0.9}{2.3}{b}{#9}%
\qbezier(-4,0)(0,8)(4,0)%
\qbezier(-4,0)(0,-8)(4,0)%
\qbezier(-0.5,3.4)(-3.5,2)(-0.5,0.6)%
\qbezier(0.5,3.4)(3.5,2)(0.5,0.6)%
\put(-4,0){\vector(1,0){8}}%
\put(0,3.4){\vector(0,-1){2.8}}%
\put(0,-0.8){\vector(0,-1){2.4}}%
\put(-1.5,2.2){\vector(1,0){1.2}}%
\put(0.3,2.2){\vector(1,0){1.2}}%
\put(4,0){\vector(1,-2){0}}%
\put(4,0){\vector(1,2){0}}%
\put(-0.5,0.6){\vector(2,-1){0}}%
\put(0.5,0.6){\vector(-2,-1){0}}%
\end{picture}}

\mcm{\gspecialone}{9}{%
\ginitdims{8}{8}%
\abovepic{#1}%
\belowpic{#3}%
\present{\pregspecialone{#1}{#2}{#3}{#4}{#5}{#6}{#7}{#8}{#9}}}

\newcommand{\pregspecialtwo}{%
\begin{picture}(5,3.4)(0,-0.2)%
\qbezier(0,1.5)(2.5,4.5)(5,1.5)%
\qbezier(0,1.5)(2.5,-1.5)(5,1.5)%
\qbezier(1.7,2.5)(0,1.5)(1.7,0.5)%
\qbezier(3.3,2.5)(5,1.5)(3.3,0.5)%
\put(5,1.5){\vector(1,-1){0}}%
\put(5,1.5){\vector(1,1){0}}%
\put(1.7,0.5){\vector(3,-2){0}}%
\put(3.3,0.5){\vector(-3,-2){0}}%
\put(2.5,2.5){\vector(0,-1){2}}%
\put(1.2,1.5){\vector(1,0){1}}%
\put(2.8,1.5){\vector(1,0){1}}%
\end{picture}}

\mcm{\gspecialtwo}{0}{%
\ginitdims{5}{3.4}%
\present{\pregspecialtwo}}

\newcommand{\pregspecialthree}{%
\begin{picture}(5,5.4)(0,-1.2)%
\qbezier(0,1.5)(2.5,6.5)(5,1.5)%
\qbezier(0,1.5)(2.5,-3.5)(5,1.5)%
\qbezier(2,3.5)(1,2.75)(2,2)%
\qbezier(3,3.5)(4,2.75)(3,2)%
\qbezier(2,1)(1,0.25)(2,-0.5)%
\qbezier(3,1)(4,0.25)(3,-0.5)%
\put(0,1.5){\vector(1,0){5}}%
\put(1.5,2.75){\vector(1,0){2}}%
\put(1.5,0.25){\vector(1,0){2}}%
\put(5,1.5){\vector(1,-3){0}}%
\put(5,1.5){\vector(1,3){0}}%
\put(2,2){\vector(1,-1){0}}%
\put(3,2){\vector(-1,-1){0}}%
\put(2,-0.5){\vector(1,-1){0}}%
\put(3,-0.5){\vector(-1,-1){0}}%
\end{picture}}

\mcm{\gspecialthree}{0}{%
\ginitdims{5}{5.4}%
\present{\pregspecialthree}}

\newcommand{\pregonew}[1]{%
\begin{picture}(8,0.4)(0,-0.2)%
\cell{4}{0.2}{b}{#1}%
\put(0,0){\vector(1,0){8}}%
\end{picture}}

\mcm{\gonew}{1}{%
\ginitdims{8}{0.4}%
\abovepic{#1}%
\present{\pregonew{#1}}}

\mcm{\gzersu}{0}{%
\gsinitdims{0}{.6}%
\present{\pregblw{}}}

\mcm{\gonesu}{0}{%
\gsinitdims{5}{0.4}%
\present{\pregone{}}}

\mcm{\gtwosu}{0}{%
\gsinitdims{5}{3.4}%
\present{\pregtwo{}{}{}}}

\mcm{\gthreesu}{0}{%
\gsinitdims{5}{5.4}%
\present{\pregthree{}{}{}{}{}}}

\mcm{\gfoursu}{0}{%
\gsinitdims{5}{8.4}%
\present{\pregfour{}{}{}{}{}{}{}}}

\newcommand{\precone}[1]{%
\begin{picture}(4.2,0.4)(-0.3,-0.2)%
\cell{1.8}{0.2}{b}{#1}%
\put(0,0){\vector(1,0){3.6}}%
\end{picture}}

\mcm{\cone}{1}{%
\cinitdims{4.2}{0.4}%
\abovepic{#1}%
\present{\precone{#1}}}

\mcm{\gfstsu}{0}{%
\gsinitdims{0.5}{0.4}%
\presentr{\pregfst{}}}

\mcm{\glstsu}{0}{%
\gsinitdims{0.5}{0.4}%
\presentl{\preglst{}}}

\newcommand{\prectwodbl}[3]%
{\begin{picture}(4.2,3.4)(-0.1,-0.2)%
\cell{2}{3.2}{b}{#1}%
\cell{2}{-0.2}{t}{#2}%
\cell{2.3}{1.5}{l}{#3}%
\qbezier(0,2)(2,4)(4,2)%
\qbezier(0,1)(2,-1)(4,1)%
\put(4,2){\vector(1,-1){0}}%
\put(4,1){\vector(1,1){0}}%
\put(1.9,2.5){\line(0,-1){1.8}}%
\put(2.1,2.5){\line(0,-1){1.8}}%
\cell{2.01}{0.4}{b}{\vee}%
\end{picture}}

\mcm{\ctwodbl}{3}{%
\cinitdims{4.2}{3.4}%
\abovepic{#1}%
\belowpic{#2}%
\present{\prectwodbl{#1}{#2}{#3}}}

\newcommand{\precthreedbl}[5]{%
\begin{picture}(4.2,5.4)(-0.1,-0.2)%
\cell{2}{5.2}{b}{#1}%
\cell{1}{2.7}{b}{#2}%
\cell{2}{-.2}{t}{#3}%
\cell{2.3}{3.75}{l}{#4}%
\cell{2.3}{1.25}{l}{#5}%
\qbezier(0,3)(2,7)(4,3)%
\qbezier(0,2)(2,-2)(4,2)%
\put(0,2.5){\vector(1,0){4}}%
\put(1.9,4.5){\line(0,-1){1.3}}%
\put(2.1,4.5){\line(0,-1){1.3}}%
\cell{2.01}{2.9}{b}{\vee}%
\put(1.9,2){\line(0,-1){1.3}}%
\put(2.1,2){\line(0,-1){1.3}}%
\cell{2.01}{0.4}{b}{\vee}%
\put(4,3){\vector(1,-3){0}}%
\put(4,2){\vector(1,3){0}}%
\end{picture}}

\mcm{\cthreedbl}{5}{%
\cinitdims{4.2}{5.4}%
\abovepic{#1}%
\belowpic{#3}%
\present{\precthreedbl{#1}{#2}{#3}{#4}{#5}}}

\newcommand{\precthreecelltrp}[5]{%
\begin{picture}(8.2,5)(-4.1,-2.5)%
\cell{0}{2.5}{b}{#1}%
\cell{0}{-2.5}{t}{#2}%
\cell{-1.8}{0}{r}{#3}%
\cell{1.8}{0}{l}{#4}%
\cell{0}{0.3}{b}{#5}%
\qbezier(-4,0.5)(0,4)(4,0.5)%
\qbezier(-4,-0.5)(0,-4)(4,-0.5)%
\qbezier(-0.6,2)(-2.6,0)(-0.6,-2)%
\qbezier(-0.4,2)(-2.4,0)(-0.5,-1.9)%
\cell{-0.6}{-2}{b}{\lrcorner}%
\qbezier(0.4,2)(2.4,0)(0.5,-1.9)%
\qbezier(0.6,2)(2.6,0)(0.6,-2)%
\cell{0.65}{-2}{b}{\llcorner}%
\put(-1,0.15){\line(1,0){1.7}}%
\put(-1,0){\line(1,0){2}}%
\put(-1,-0.15){\line(1,0){1.7}}%
\cell{1.15}{0}{r}{>}%
\put(4,0.5){\vector(1,-1){0}}%
\put(4,-0.5){\vector(1,1){0}}%
\end{picture}}

\mcm{\cthreecelltrp}{5}{%
\cinitdims{8.2}{5}%
\abovepic{#1}%
\belowpic{#2}%
\present{\precthreecelltrp{#1}{#2}{#3}{#4}{#5}}}

\newcommand{\prectwo}[3]%
{\begin{picture}(4.2,3.4)(-0.1,-0.2)%
\cell{2}{3.2}{b}{#1}%
\cell{2}{-0.2}{t}{#2}%
\cell{2.2}{1.5}{l}{#3}%
\qbezier(0,2)(2,4)(4,2)%
\qbezier(0,1)(2,-1)(4,1)%
\put(4,2){\vector(1,-1){0}}%
\put(4,1){\vector(1,1){0}}%
\put(2,2.5){\vector(0,-1){2}}%
\end{picture}}

\mcm{\ctwo}{3}{%
\cinitdims{4.2}{3.4}%
\abovepic{#1}%
\belowpic{#2}%
\present{\prectwo{#1}{#2}{#3}}}

\newcommand{\precthree}[5]{%
\begin{picture}(4.2,5.4)(-0.1,-0.2)%
\cell{2}{5.2}{b}{#1}%
\cell{1}{2.7}{b}{#2}%
\cell{2}{-.2}{t}{#3}%
\cell{2.2}{3.75}{l}{#4}%
\cell{2.2}{1.25}{l}{#5}%
\qbezier(0,3)(2,7)(4,3)%
\qbezier(0,2)(2,-2)(4,2)%
\put(0,2.5){\vector(1,0){4}}%
\put(2,4.5){\vector(0,-1){1.5}}%
\put(2,2){\vector(0,-1){1.5}}%
\put(4,3){\vector(1,-3){0}}%
\put(4,2){\vector(1,3){0}}%
\end{picture}}

\mcm{\cthree}{5}{%
\cinitdims{4.2}{5.4}%
\abovepic{#1}%
\belowpic{#3}%
\present{\precthree{#1}{#2}{#3}{#4}{#5}}}

\newcommand{\prectwoop}[3]%
{\begin{picture}(4.2,3.4)(-0.1,-0.2)%
\cell{2}{3.2}{b}{#1}%
\cell{2}{-0.2}{t}{#2}%
\cell{2.2}{1.5}{l}{#3}%
\qbezier(0,2)(2,4)(4,2)%
\qbezier(0,1)(2,-1)(4,1)%
\put(0,2){\vector(-1,-1){0}}%
\put(0,1){\vector(-1,1){0}}%
\put(2,2.5){\vector(0,-1){2}}%
\end{picture}}

\mcm{\ctwoop}{3}{%
\cinitdims{4.2}{3.4}%
\abovepic{#1}%
\belowpic{#2}%
\present{\prectwoop{#1}{#2}{#3}}}

\newcommand{\prectwopar}[4]{%
\begin{picture}(4.2,3.4)(-0.1,-0.2)%
\cell{2}{3.2}{b}{#1}%
\cell{2}{-0.2}{t}{#2}%
\cell{1.6}{1.5}{r}{#3}%
\cell{2.4}{1.5}{l}{#4}%
\qbezier(0,2)(2,4)(4,2)%
\qbezier(0,1)(2,-1)(4,1)%
\put(4,2){\vector(1,-1){0}}%
\put(4,1){\vector(1,1){0}}%
\put(1.8,2.5){\vector(0,-1){2}}%
\put(2.2,2.5){\vector(0,-1){2}}%
\end{picture}}

\mcm{\ctwopar}{4}{%
\cinitdims{4.2}{3.4}%
\abovepic{#1}%
\belowpic{#2}%
\present{\prectwopar{#1}{#2}{#3}{#4}}}

\newcommand{\precthreein}[5]{%
\begin{picture}(4.2,5.4)(-0.1,-0.2)%
\cell{2}{5.2}{b}{#1}%
\cell{1}{2.7}{b}{#2}%
\cell{2}{-.2}{t}{#3}%
\cell{2.2}{3.75}{l}{#4}%
\cell{2.2}{1.25}{l}{#5}%
\qbezier(0,3)(2,7)(4,3)%
\qbezier(0,2)(2,-2)(4,2)%
\put(0,2.5){\vector(1,0){4}}%
\put(2,4.5){\vector(0,-1){1.5}}%
\put(2,0.5){\vector(0,1){1.5}}%
\put(4,3){\vector(1,-3){0}}%
\put(4,2){\vector(1,3){0}}%
\end{picture}}

\mcm{\cthreein}{5}{%
\cinitdims{4.2}{5.4}%
\abovepic{#1}%
\belowpic{#3}%
\present{\precthreein{#1}{#2}{#3}{#4}{#5}}}

\newcommand{\precthreecell}[5]{%
\begin{picture}(8.2,5)(-4.1,-2.5)%
\cell{0}{2.5}{b}{#1}%
\cell{0}{-2.5}{t}{#2}%
\cell{-1.7}{0}{r}{#3}%
\cell{1.7}{0}{l}{#4}%
\cell{0}{0.2}{b}{#5}%
\qbezier(-4,0.5)(0,4)(4,0.5)%
\qbezier(-4,-0.5)(0,-4)(4,-0.5)%
\qbezier(-0.5,2)(-2.5,0)(-0.5,-2)%
\qbezier(0.5,2)(2.5,0)(0.5,-2)%
\put(-1,0){\vector(1,0){2}}%
\put(4,0.5){\vector(1,-1){0}}%
\put(4,-0.5){\vector(1,1){0}}%
\put(-0.5,-2){\vector(1,-1){0}}%
\put(0.5,-2){\vector(-1,-1){0}}%
\end{picture}}

\mcm{\cthreecell}{5}{%
\cinitdims{8.2}{5}%
\abovepic{#1}%
\belowpic{#2}%
\present{\precthreecell{#1}{#2}{#3}{#4}{#5}}}

\newcommand{\precthreecellpar}[6]{%
\begin{picture}(8.2,5)(-4.1,-2.5)%
\cell{0}{2.5}{b}{#1}%
\cell{0}{-2.5}{t}{#2}%
\cell{-1.7}{0}{r}{#3}%
\cell{1.7}{0}{l}{#4}%
\cell{0}{0.4}{b}{#5}%
\cell{0}{-0.4}{t}{#6}%
\qbezier(-4,0.5)(0,4)(4,0.5)%
\qbezier(-4,-0.5)(0,-4)(4,-0.5)%
\qbezier(-0.5,2)(-2.5,0)(-0.5,-2)%
\qbezier(0.5,2)(2.5,0)(0.5,-2)%
\put(-1,0.2){\vector(1,0){2}}%
\put(-1,-0.2){\vector(1,0){2}}%
\put(4,0.5){\vector(1,-1){0}}%
\put(4,-0.5){\vector(1,1){0}}%
\put(-0.5,-2){\vector(1,-1){0}}%
\put(0.5,-2){\vector(-1,-1){0}}%
\end{picture}}

\mcm{\cthreecellpar}{6}{%
\cinitdims{8.2}{5}%
\abovepic{#1}%
\belowpic{#2}%
\present{\precthreecellpar{#1}{#2}{#3}{#4}{#5}{#6}}}

\newcommand{\prectwov}[5]{%
\begin{picture}(3.4,4.2)(0.8,0.9)%
\cell{2.5}{5.1}{b}{#1}%
\cell{2.5}{0.9}{t}{#2}%
\cell{0.8}{3}{r}{#3}%
\cell{4.2}{3}{l}{#4}%
\cell{2.5}{3.2}{b}{#5}%
\qbezier(2,5)(0,3)(2,1)%
\qbezier(3,5)(5,3)(3,1)%
\put(2,1){\vector(1,-1){0}}%
\put(3,1){\vector(-1,-1){0}}%
\put(1.5,3){\vector(1,0){2}}%
\end{picture}}

\mcm{\ctwov}{5}{%
\cinitdims{3.4}{4.2}%
\abovepic{#1}%
\belowpic{#2}%
\sidespic{#3}%
\sidespic{#4}%
\present{\prectwov{#1}{#2}{#3}{#4}{#5}}}

\newcommand{\precthreecellv}[7]{%
\begin{picture}(5,8.2)(0.5,-1.6)%
\cell{3}{6.6}{b}{#1}%
\cell{3}{-1.6}{t}{#2}%
\cell{0.5}{2.5}{r}{#3}%
\cell{5.5}{2.5}{l}{#4}%
\cell{3}{4.2}{b}{#5}%
\cell{3}{0.8}{t}{#6}%
\cell{3.2}{2.5}{l}{#7}%
\qbezier(3.5,6.5)(7,2.5)(3.5,-1.5)%
\qbezier(2.5,6.5)(-1,2.5)(2.5,-1.5)%
\put(2.5,-1.5){\vector(1,-1){0}}%
\put(3.5,-1.5){\vector(-1,-1){0}}%
\qbezier(1,3)(3,5)(5,3)%
\qbezier(1,2)(3,0)(5,2)%
\put(5,3){\vector(1,-1){0}}%
\put(5,2){\vector(1,1){0}}%
\put(3,3.5){\vector(0,-1){2}}%
\end{picture}}

\mcm{\cthreecellv}{7}{%
\cinitdims{5}{8.2}%
\abovepic{#1}%
\belowpic{#2}%
\sidespic{#3}%
\sidespic{#4}%
\present{\precthreecellv{#1}{#2}{#3}{#4}{#5}{#6}{#7}}}

\newcommand{\pretopez}[2]{%
\begin{picture}(2.6,2.3)(-1.3,-2.2)%
\cell{0}{-2.2}{t}{#1}%
\cell{0}{-1.2}{c}{#2}%
\qbezier(0,0)(-2,-2)(0,-2)%
\qbezier(0,0)(2,-2)(0,-2)%
\put(0,0){\vector(-1,1){0}}%
\end{picture}}

\mcm{\topez}{2}{%
\ginitdims{2.6}{2.3}%
\belowpic{#1}%
\present{\pretopez{#1}{#2}}}

\newcommand{\pretopea}[3]{%
\begin{picture}(4,1.9)(-2,-0,2)%
\cell{0}{1.7}{b}{#1}%
\cell{0}{-0.2}{t}{#2}%
\cell{0}{0.7}{c}{#3}%
\qbezier(-2,0)(0,3)(2,0)%
\put(-2,0){\vector(1,0){4}}%
\put(2,0){\vector(2,-3){0}}%
\end{picture}}

\mcm{\topea}{3}{%
\ginitdims{4}{1.9}%
\abovepic{#1}%
\belowpic{#2}%
\present{\pretopea{#1}{#2}{#3}}}

\newcommand{\pretopeb}[4]{%
\begin{picture}(4,2.2)(-2,-0.2)%
\cell{-1.1}{1}{br}{#1}%
\cell{1.1}{1}{bl}{#2}%
\cell{0}{-0.2}{t}{#3}%
\cell{0}{0.8}{c}{#4}%
\put(-2,0){\vector(1,1){2}}%
\put(0,2){\vector(1,-1){2}}%
\put(-2,0){\vector(1,0){4}}%
\end{picture}}

\mcm{\topeb}{4}{%
\ginitdims{4}{2.2}%
\belowpic{#3}%
\present{\pretopeb{#1}{#2}{#3}{#4}}}

\newcommand{\pretopec}[5]{%
\begin{picture}(4,2.2)(-2,-0.2)%
\cell{-1.8}{1}{br}{#1}%
\cell{0}{2.2}{b}{#2}%
\cell{1.8}{1}{bl}{#3}%
\cell{0}{-0.2}{t}{#4}%
\cell{0}{0.8}{c}{#5}%
\put(-2,0){\vector(1,2){1}}%
\put(-1,2){\vector(1,0){2}}%
\put(1,2){\vector(1,-2){1}}%
\put(-2,0){\vector(1,0){4}}%
\end{picture}}

\mcm{\topec}{5}{%
\ginitdims{4}{2.2}%
\sidespic{#1}%
\abovepic{#2}%
\sidespic{#3}%
\belowpic{#4}%
\present{\pretopec{#1}{#2}{#3}{#4}{#5}}}

\newcommand{\pretoped}[6]{%
\begin{picture}(4,2.5)(-2,-0.2)%
\cell{-2}{0.6}{br}{#1}%
\cell{-0.7}{2.2}{br}{#2}%
\cell{0.7}{2.2}{bl}{#3}%
\cell{2}{0.6}{bl}{#4}%
\cell{0}{-0.2}{t}{#5}%
\cell{0}{0.8}{c}{#6}%
\put(-2,0){\vector(1,3){0.5}}%
\put(-1.5,1.5){\vector(3,2){1.5}}%
\put(0,2.5){\vector(3,-2){1.5}}%
\put(1.5,1.5){\vector(1,-3){0.5}}%
\put(-2,0){\vector(1,0){4}}%
\end{picture}}

\mcm{\toped}{6}{%
\ginitdims{4}{2.5}%
\sidespic{#1}%
\abovepic{#2}%
\abovepic{#3}%
\sidespic{#4}%
\belowpic{#5}%
\present{\pretoped{#1}{#2}{#3}{#4}{#5}{#6}}}

\newcommand{\pretopeq}[5]{%
\begin{picture}(4,2.5)(-2,-0.2)%
\cell{-2}{0.6}{br}{#1}%
\cell{-1}{2.2}{br}{#2}%
\cell{2}{0.6}{bl}{#3}%
\cell{0}{-0.2}{t}{#4}%
\cell{0}{0.8}{c}{#5}%
\put(-2,0){\vector(1,3){0.5}}%
\put(-1.5,1.5){\vector(1,1){1}}%
\cell{0.9}{2.3}{c}{\ddots}
\put(1.5,1.5){\vector(1,-3){0.5}}%
\put(-2,0){\vector(1,0){4}}%
\end{picture}}

\mcm{\topeq}{5}{%
\ginitdims{4}{2.5}%
\sidespic{#1}%
\abovepic{#2}%
\sidespic{#3}%
\belowpic{#4}%
\present{\pretopeq{#1}{#2}{#3}{#4}{#5}}}

\newcommand{\pretopebase}[1]{%
\begin{picture}(4,0.4)(0,-0.2)%
\cell{2}{0.2}{b}{#1}%
\put(0,0){\vector(1,0){4}}%
\end{picture}}

\mcm{\topebase}{1}{%
\ginitdims{4}{0.4}%
\abovepic{#1}%
\present{\pretopebase{#1}}}

\newcommand{\pretopezs}[2]{%
\begin{picture}(2.6,2.3)(-1.3,-2.2)%
\cell{0}{-2.2}{t}{#1}%
\cell{0}{-1.2}{c}{#2}%
\qbezier(0,0)(-2,-2)(0,-2)%
\qbezier(0,0)(2,-2)(0,-2)%
\end{picture}}

\mcm{\topezs}{2}{%
\ginitdims{2.6}{2.3}%
\belowpic{#1}%
\present{\pretopezs{#1}{#2}}}

\newcommand{\pretopeas}[3]{%
\begin{picture}(4,1.9)(-2,-0,2)%
\cell{0}{1.7}{b}{#1}%
\cell{0}{-0.2}{t}{#2}%
\cell{0}{0.7}{c}{#3}%
\qbezier(-2,0)(0,3)(2,0)%
\put(-2,0){\line(1,0){4}}%
\end{picture}}

\mcm{\topeas}{3}{%
\ginitdims{4}{1.9}%
\abovepic{#1}%
\belowpic{#2}%
\present{\pretopeas{#1}{#2}{#3}}}

\newcommand{\pretopebs}[4]{%
\begin{picture}(4,2.2)(-2,-0.2)%
\cell{-1.1}{1}{br}{#1}%
\cell{1.1}{1}{bl}{#2}%
\cell{0}{-0.2}{t}{#3}%
\cell{0}{0.8}{c}{#4}%
\put(-2,0){\line(1,1){2}}%
\put(0,2){\line(1,-1){2}}%
\put(-2,0){\line(1,0){4}}%
\end{picture}}

\mcm{\topebs}{4}{%
\ginitdims{4}{2.2}%
\belowpic{#3}%
\present{\pretopebs{#1}{#2}{#3}{#4}}}

\newcommand{\pretopecs}[5]{%
\begin{picture}(4,2.2)(-2,-0.2)%
\cell{-1.8}{1}{br}{#1}%
\cell{0}{2.2}{b}{#2}%
\cell{1.8}{1}{bl}{#3}%
\cell{0}{-0.2}{t}{#4}%
\cell{0}{0.8}{c}{#5}%
\put(-2,0){\line(1,2){1}}%
\put(-1,2){\line(1,0){2}}%
\put(1,2){\line(1,-2){1}}%
\put(-2,0){\line(1,0){4}}%
\end{picture}}

\mcm{\topecs}{5}{%
\ginitdims{4}{2.2}%
\sidespic{#1}%
\abovepic{#2}%
\sidespic{#3}%
\belowpic{#4}%
\present{\pretopecs{#1}{#2}{#3}{#4}{#5}}}

\newcommand{\pretopeds}[6]{%
\begin{picture}(4,2.5)(-2,-0.2)%
\cell{-2}{0.6}{br}{#1}%
\cell{-0.7}{2.2}{br}{#2}%
\cell{0.7}{2.2}{bl}{#3}%
\cell{2}{0.6}{bl}{#4}%
\cell{0}{-0.2}{t}{#5}%
\cell{0}{0.8}{c}{#6}%
\put(-2,0){\line(1,3){0.5}}%
\put(-1.5,1.5){\line(3,2){1.5}}%
\put(0,2.5){\line(3,-2){1.5}}%
\put(1.5,1.5){\line(1,-3){0.5}}%
\put(-2,0){\line(1,0){4}}%
\end{picture}}

\mcm{\topeds}{6}{%
\ginitdims{4}{2.5}%
\sidespic{#1}%
\abovepic{#2}%
\abovepic{#3}%
\sidespic{#4}%
\belowpic{#5}%
\present{\pretopeds{#1}{#2}{#3}{#4}{#5}{#6}}}

\newcommand{\pretopeqs}[5]{%
\begin{picture}(4,2.5)(-2,-0.2)%
\cell{-2}{0.6}{br}{#1}%
\cell{-1}{2.2}{br}{#2}%
\cell{2}{0.6}{bl}{#3}%
\cell{0}{-0.2}{t}{#4}%
\cell{0}{0.8}{c}{#5}%
\put(-2,0){\line(1,3){0.5}}%
\put(-1.5,1.5){\line(1,1){1}}%
\cell{0.9}{2.3}{c}{\ddots}
\put(1.5,1.5){\line(1,-3){0.5}}%
\put(-2,0){\line(1,0){4}}%
\end{picture}}

\mcm{\topeqs}{5}{%
\ginitdims{4}{2.5}%
\sidespic{#1}%
\abovepic{#2}%
\sidespic{#3}%
\belowpic{#4}%
\present{\pretopeqs{#1}{#2}{#3}{#4}{#5}}}

\newcommand{\pretopebases}[1]{%
\begin{picture}(4,0.4)(0,-0.2)%
\cell{2}{0.2}{b}{#1}%
\put(0,0){\line(1,0){4}}%
\end{picture}}

\mcm{\topebases}{1}{%
\ginitdims{4}{0.4}%
\abovepic{#1}%
\present{\pretopebases{#1}}}

\newcommand{\pregdots}[6]{%
\begin{picture}(5,8.4)(0,-2.7)%
\cell{2.5}{5.7}{b}{#1}%
\cell{1.5}{2.8}{b}{#2}%
\cell{1.5}{0.2}{t}{#3}%
\cell{2.5}{-2.7}{t}{#4}%
\cell{2.7}{4.25}{l}{#5}%
\cell{2.7}{-1.25}{l}{#6}%
\qbezier(0,1.5)(2.5,9.5)(5,1.5)%
\qbezier(0,1.5)(2.5,4)(5,1.5)%
\qbezier(0,1.5)(2.5,-1)(5,1.5)%
\qbezier(0,1.5)(2.5,-6.5)(5,1.5)%
\put(2.5,5.25){\vector(0,-1){2}}%
\put(2.5,-0.25){\vector(0,-1){2}}%
\cell{2.5}{1.7}{c}{\vdots}%
\put(5,1.5){\vector(1,-4){0}}%
\put(5,1.5){\vector(4,-3){0}}%
\put(5,1.5){\vector(4,3){0}}%
\put(5,1.5){\vector(1,4){0}}%
\end{picture}}

\mcm{\gdots}{6}{%
\ginitdims{5}{8.4}%
\abovepic{#1}%
\belowpic{#4}%
\present{\pregdots{#1}{#2}{#3}{#4}{#5}{#6}}}

\newlength{\volt}
\makeatletter

\def\diagram{\m@th\leftwidth=\z@ \rightwidth=\z@ \topheight=\z@
\botheight=\z@ \setbox\@picbox\hbox\bgroup}

\def\enddiagram{\egroup\wd\@picbox\rightwidth\unitlength
\ht\@picbox\topheight\unitlength \dp\@picbox\botheight\unitlength
\hskip\leftwidth\unitlength\box\@picbox}

\def\bfig{\begin{diagram}}
\def\efig{\end{diagram}}
\newcount\wideness \newcount\leftwidth \newcount\rightwidth
\newcount\highness \newcount\topheight \newcount\botheight

\def\ratchet#1#2{\ifnum#1<#2 \global #1=#2 \fi}

\def\putbox(#1,#2)#3{%
\horsize{\wideness}{#3} \divide\wideness by 2 {\advance\wideness
by #1 \ratchet{\rightwidth}{\wideness}} {\advance\wideness by -#1
\ratchet{\leftwidth}{\wideness}} \vertsize{\highness}{#3}
\divide\highness by 2 {\advance\highness by #2
\ratchet{\topheight}{\highness}} {\advance\highness by -#2
\ratchet{\botheight}{\highness}} \put(#1,#2){\makebox(0,0){$#3$}}}

\def\putlbox(#1,#2)#3{%
\horsize{\wideness}{#3} {\advance\wideness by #1
\ratchet{\rightwidth}{\wideness}} {\ratchet{\leftwidth}{-#1}}
\vertsize{\highness}{#3} \divide\highness by 2 {\advance\highness
by #2 \ratchet{\topheight}{\highness}} {\advance\highness by -#2
\ratchet{\botheight}{\highness}}
\put(#1,#2){\makebox(0,0)[l]{$#3$}}}

\def\putrbox(#1,#2)#3{%
\horsize{\wideness}{#3} {\ratchet{\rightwidth}{#1}}
{\advance\wideness by -#1 \ratchet{\leftwidth}{\wideness}}
\vertsize{\highness}{#3} \divide\highness by 2 {\advance\highness
by #2 \ratchet{\topheight}{\highness}} {\advance\highness by -#2
\ratchet{\botheight}{\highness}}
\put(#1,#2){\makebox(0,0)[r]{$#3$}}}

\def\adjust[#1]{} 

\newcount \coefa
\newcount \coefb
\newcount \coefc
\newcount\tempcounta
\newcount\tempcountb
\newcount\tempcountc
\newcount\tempcountd
\newcount\xext
\newcount\yext
\newcount\xoff
\newcount\yoff
\newcount\gap%
\newcount\arrowtypea
\newcount\arrowtypeb
\newcount\arrowtypec
\newcount\arrowtyped
\newcount\arrowtypee
\newcount\height
\newcount\width
\newcount\xpos
\newcount\ypos
\newcount\run
\newcount\rise
\newcount\arrowlength
\newcount\halflength
\newcount\arrowtype
\newdimen\tempdimen
\newdimen\xlen
\newdimen\ylen
\newsavebox{\tempboxa}%
\newsavebox{\tempboxb}%
\newsavebox{\tempboxc}%

\newdimen\w@dth

\def\setw@dth#1#2{\setbox\z@\hbox{\m@th$#1$}\w@dth=\wd\z@
\setbox\@ne\hbox{\m@th$#2$}\ifnum\w@dth<\wd\@ne \w@dth=\wd\@ne \fi
\advance\w@dth by 1.2em}

\def\t@^#1_#2{\allowbreak\def\n@one{#1}\def\n@two{#2}\mathrel
{\setw@dth{#1}{#2} \mathop{\hbox to
\w@dth{\rightarrowfill}}\limits \ifx\n@one\empty\else
^{\box\z@}\fi \ifx\n@two\empty\else _{\box\@ne}\fi}}
\def\t@@^#1{\@ifnextchar_{\t@^{#1}}{\t@^{#1}_{}}}
\def\to{\@ifnextchar^{\t@@}{\t@@^{}}}

\def\t@left^#1_#2{\def\n@one{#1}\def\n@two{#2}\mathrel{\setw@dth{#1}{#2}
\mathop{\hbox to \w@dth{\leftarrowfill}}\limits
\ifx\n@one\empty\else ^{\box\z@}\fi \ifx\n@two\empty\else
_{\box\@ne}\fi}}
\def\t@@left^#1{\@ifnextchar_{\t@left^{#1}}{\t@left^{#1}_{}}}
\def\toleft{\@ifnextchar^{\t@@left}{\t@@left^{}}}

\def\two@^#1_#2{\allowbreak
\def\n@one{#1}\def\n@two{#2}\mathrel{\setw@dth{#1}{#2}
\mathop{\vcenter{\lineskip\z@\baselineskip\z@
                 \hbox to \w@dth{\rightarrowfill}%
                 \hbox to \w@dth{\rightarrowfill}}%
       }\limits
\ifx\n@one\empty\else ^{\box\z@}\fi \ifx\n@two\empty\else
_{\box\@ne}\fi}}
\def\tw@@^#1{\@ifnextchar _{\two@^{#1}}{\two@^{#1}_{}}}
\def\two{\@ifnextchar ^{\tw@@}{\tw@@^{}}}

\def\tofr@^#1_#2{\def\n@one{#1}\def\n@two{#2}\mathrel{\setw@dth{#1}{#2}
\mathop{\vcenter{\hbox to \w@dth{\rightarrowfill}\kern-1.7ex
                 \hbox to \w@dth{\leftarrowfill}}%
       }\limits
\ifx\n@one\empty\else ^{\box\z@}\fi \ifx\n@two\empty\else
_{\box\@ne}\fi}}
\def\t@fr@^#1{\@ifnextchar_ {\tofr@^{#1}}{\tofr@^{#1}_{}}}
\def\tofro{\@ifnextchar^ {\t@fr@}{\t@fr@^{}}}

\def\mon{\mathop{\m@th\hbox to
      14.6\P@{\lasyb\char'51\hskip-2.1\P@$\arrext$\hss
$\mathord\rightarrow$}}\limits} 
\def\leftmono{\mathrel{\m@th\hbox to
14.6\P@{$\mathord\leftarrow$\hss$\arrext$\hskip-2.1\P@\lasyb\char'50%
}}\limits} 
\mathchardef\arrext="0200       

\def\settypes(#1,#2,#3){\arrowtypea#1 \arrowtypeb#2 \arrowtypec#3}
\def\settoheight#1#2{\setbox\@tempboxa\hbox{#2}#1\ht\@tempboxa\relax}%
\def\settodepth#1#2{\setbox\@tempboxa\hbox{#2}#1\dp\@tempboxa\relax}%
\def\settokens`#1`#2`#3`#4`{%
     \def\tokena{#1}\def\tokenb{#2}\def\tokenc{#3}\def\tokend{#4}}
\def\setsqparms[#1`#2`#3`#4;#5`#6]{%
\arrowtypea #1 \arrowtypeb #2 \arrowtypec #3 \arrowtyped #4
\width #5 \height #6 }
\def\setpos(#1,#2){\xpos=#1 \ypos#2}

\def\settriparms[#1`#2`#3;#4]{\settripairparms[#1`#2`#3`1`1;#4]}%

\def\settripairparms[#1`#2`#3`#4`#5;#6]{%
\arrowtypea #1 \arrowtypeb #2 \arrowtypec #3 \arrowtyped #4
\arrowtypee #5 \width #6 \height #6 }

\def\resetparms{\settripairparms[1`1`1`1`1;500]\width 500}

\resetparms

\def\mvector(#1,#2)#3{
\put(0,0){\vector(#1,#2){#3}}%
\put(0,0){\vector(#1,#2){26}}%
}
\def\evector(#1,#2)#3{{
\arrowlength #3
\put(0,0){\vector(#1,#2){\arrowlength}}%
\advance \arrowlength by-30
\put(0,0){\vector(#1,#2){\arrowlength}}%
}}

\def\horsize#1#2{%
\settowidth{\tempdimen}{$#2$}%
#1=\tempdimen \divide #1 by\unitlength }

\def\vertsize#1#2{%
\settoheight{\tempdimen}{$#2$}%
#1=\tempdimen
\settodepth{\tempdimen}{$#2$}%
\advance #1 by\tempdimen \divide #1 by\unitlength }

\def\putvector(#1,#2)(#3,#4)#5#6{{%
\ifnum3<\arrowtype \putdashvector(#1,#2)(#3,#4)#5\arrowtype \else
\ifnum\arrowtype<-3 \putdashvector(#1,#2)(#3,#4)#5\arrowtype \else
\xpos=#1 \ypos=#2 \run=#3 \rise=#4 \arrowlength=#5 \ifnum
\arrowtype<0
    \ifnum \run=0
        \advance \ypos by-\arrowlength
    \else
        \tempcounta \arrowlength
        \multiply \tempcounta by\rise
        \divide \tempcounta by\run
        \ifnum\run>0
            \advance \xpos by\arrowlength
            \advance \ypos by\tempcounta
        \else
            \advance \xpos by-\arrowlength
            \advance \ypos by-\tempcounta
        \fi
    \fi
    \multiply \arrowtype by-1
    \multiply \rise by-1
    \multiply \run by-1
\fi \ifcase \arrowtype
\or \put(\xpos,\ypos){\vector(\run,\rise){\arrowlength}}%
\or \put(\xpos,\ypos){\mvector(\run,\rise)\arrowlength}%
\or \put(\xpos,\ypos){\evector(\run,\rise){\arrowlength}}%
\fi\fi\fi }}

\def\putsplitvector(#1,#2)#3#4{
\xpos #1 \ypos #2 \arrowtype #4 \halflength #3 \arrowlength #3
\gap 140 \advance \halflength by-\gap \divide \halflength by2
\ifnum\arrowtype>0
   \ifcase \arrowtype
   \or \put(\xpos,\ypos){\line(0,-1){\halflength}}%
       \advance\ypos by-\halflength
       \advance\ypos by-\gap
       \put(\xpos,\ypos){\vector(0,-1){\halflength}}%
   \or \put(\xpos,\ypos){\line(0,-1)\halflength}%
       \put(\xpos,\ypos){\vector(0,-1)3}%
       \advance\ypos by-\halflength
       \advance\ypos by-\gap
       \put(\xpos,\ypos){\vector(0,-1){\halflength}}%
   \or \put(\xpos,\ypos){\line(0,-1)\halflength}%
       \advance\ypos by-\halflength
       \advance\ypos by-\gap
       \put(\xpos,\ypos){\evector(0,-1){\halflength}}%
   \fi
\else \arrowtype=-\arrowtype
   \ifcase\arrowtype
   \or \advance \ypos by-\arrowlength
       \put(\xpos,\ypos){\line(0,1){\halflength}}%
       \advance\ypos by\halflength
       \advance\ypos by\gap
       \put(\xpos,\ypos){\vector(0,1){\halflength}}%
   \or \advance \ypos by-\arrowlength
       \put(\xpos,\ypos){\line(0,1)\halflength}%
       \put(\xpos,\ypos){\vector(0,1)3}%
       \advance\ypos by\halflength
       \advance\ypos by\gap
       \put(\xpos,\ypos){\vector(0,1){\halflength}}%
   \or \advance \ypos by-\arrowlength
       \put(\xpos,\ypos){\line(0,1)\halflength}%
       \advance\ypos by\halflength
       \advance\ypos by\gap
       \put(\xpos,\ypos){\evector(0,1){\halflength}}%
   \fi
\fi }

\def\putmorphism(#1)(#2,#3)[#4`#5`#6]#7#8#9{{%
\run #2 \rise #3 \ifnum\rise=0
  \puthmorphism(#1)[#4`#5`#6]{#7}{#8}#9%
\else\ifnum\run=0
  \putvmorphism(#1)[#4`#5`#6]{#7}{#8}#9%
\else
\setpos(#1)%
\arrowlength #7 \arrowtype #8 \ifnum\run=0 \else\ifnum\rise=0
\else \ifnum\run>0
    \coefa=1
\else
   \coefa=-1
\fi \ifnum\arrowtype>0
   \coefb=0
   \coefc=-1
\else
   \coefb=\coefa
   \coefc=1
   \arrowtype=-\arrowtype
\fi \width=2 \multiply \width by\run \divide \width by\rise
\ifnum \width<0  \width=-\width\fi \advance\width by60 \if l#9
\width=-\width\fi
\putbox(\xpos,\ypos){#4}
{\multiply \coefa by\arrowlength
\advance\xpos by\coefa \multiply \coefa by\rise \divide \coefa
by\run \advance \ypos by\coefa
\putbox(\xpos,\ypos){#5} }%
{\multiply \coefa by\arrowlength
\divide \coefa by2 \advance \xpos by\coefa \advance \xpos by\width
\multiply \coefa by\rise \divide \coefa by\run \advance \ypos
by\coefa
\if l#9%
   \putrbox(\xpos,\ypos){#6}%
\else\if r#9%
   \putlbox(\xpos,\ypos){#6}%
\fi\fi }%
{\multiply \rise by-\coefc
\multiply \run by-\coefc \multiply \coefb by\arrowlength \advance
\xpos by\coefb \multiply \coefb by\rise \divide \coefb by\run
\advance \ypos by\coefb \multiply \coefc by70 \advance \ypos
by\coefc \multiply \coefc by\run \divide \coefc by\rise \advance
\xpos by\coefc \multiply \coefa by140 \multiply \coefa by\run
\divide \coefa by\rise \advance \arrowlength by\coefa
\ifcase\arrowtype
\or \put(\xpos,\ypos){\vector(\run,\rise){\arrowlength}}%
\or \put(\xpos,\ypos){\mvector(\run,\rise){\arrowlength}}%
\or \put(\xpos,\ypos){\evector(\run,\rise){\arrowlength}}%
\fi}\fi\fi\fi\fi}}

\newcount\numbdashes \newcount\lengthdash \newcount\increment

\def\howmanydashes{
\numbdashes=\arrowlength \lengthdash=40 \divide\numbdashes by
\lengthdash \lengthdash=\arrowlength \divide\lengthdash by
\numbdashes
\increment=\lengthdash \multiply\lengthdash by 3
\divide\lengthdash by 5 }

\def\putdashvector(#1)(#2,#3)#4#5{%
\ifnum#3=0 \putdashhvector(#1){#4}#5 \else \ifnum#2=0
\putdashvvector(#1){#4}#5\fi\fi}

\def\putdashhvector(#1,#2)#3#4{{%
\arrowlength=#3 \howmanydashes
\multiput(#1,#2)(\increment,0){\numbdashes}%
{\vrule height .4pt width \lengthdash\unitlength} \arrowtype=#4
\xpos=#1 \ifnum\arrowtype<0 \advance\arrowtype by 7 \fi
\ifcase\arrowtype \or \advance\xpos by 10
    \put(\xpos,#2){\vector(-1,0){\lengthdash}}
    \advance\xpos by 40
    \put(\xpos,#2){\vector(-1,0){\lengthdash}}
\or \advance \xpos by 10
    \put(\xpos,#2){\vector(-1,0){\lengthdash}}
    \advance\xpos by  \arrowlength
    \advance\xpos by  -50
    \put(\xpos,#2){\vector(-1,0){\lengthdash}}
\or \advance\xpos by 10
    \put(\xpos,#2){\vector(-1,0){\lengthdash}}
\or \advance\xpos by \arrowlength
    \advance\xpos by -\lengthdash
    \put(\xpos,#2){\vector(1,0){\lengthdash}}
\or {\advance\xpos by 10
    \put(\xpos,#2){\vector(1,0){\lengthdash}}}
    \advance\xpos by \arrowlength
    \advance\xpos by -\lengthdash
    \put(\xpos,#2){\vector(1,0){\lengthdash}}
\or \advance\xpos by \arrowlength
    \advance\xpos by -\lengthdash
    \put(\xpos,#2){\vector(1,0){\lengthdash}}
    \advance\xpos by -40
    \put(\xpos,#2){\vector(1,0){\lengthdash}}
   \fi
}}

\def\putdashvvector(#1,#2)#3#4{{%
\arrowlength=#3 \howmanydashes \ypos=#2 \advance\ypos by
-\arrowlength
\multiput(#1,#2)(0,\increment){\numbdashes}%
    {\vrule width .4pt height \lengthdash\unitlength}
\arrowtype=#4 \ypos=#2 \ifnum\arrowtype<0 \advance\arrowtype by 7
\fi \ifcase\arrowtype \or \advance\ypos by \arrowlength
\advance\ypos by -40
    \put(#1,\ypos){\vector(0,1){\lengthdash}}
    \advance\ypos by -40
    \put(#1,\ypos){\vector(0,1){\lengthdash}}
\or \advance\ypos by 10
    \put(#1,\ypos){\vector(0,1){\lengthdash}}
    \advance\ypos by \arrowlength \advance\ypos by -40
    \put(#1,\ypos){\vector(0,1){\lengthdash}}
\or \advance\ypos by \arrowlength \advance\ypos by -40
    \put(#1,\ypos){\vector(0,1){\lengthdash}}
\or \advance\ypos by 10
    \put(#1,\ypos){\vector(0,-1){\lengthdash}}
\or \advance\ypos by 10
    \put(#1,\ypos){\vector(0,-1){\lengthdash}}
    \advance\ypos by \arrowlength \advance\ypos by -40
    \put(#1,\ypos){\vector(0,-1){\lengthdash}}
\or \advance\ypos by 10
    \put(#1,\ypos){\vector(0,-1){\lengthdash}}
    \advance\ypos by 40
    \put(#1,\ypos){\vector(0,-1){\lengthdash}}
\fi }}

\def\puthmorphism(#1,#2)[#3`#4`#5]#6#7#8{{%
\xpos #1 \ypos #2 \width #6 \arrowlength #6 \arrowtype=#7
\putbox(\xpos,\ypos){#3\vphantom{#4}}%
{\advance \xpos by\arrowlength
\putbox(\xpos,\ypos){\vphantom{#3}#4}}%
\horsize{\tempcounta}{#3}%
\horsize{\tempcountb}{#4}%
\divide \tempcounta by2 \divide \tempcountb by2 \advance
\tempcounta by30 \advance \tempcountb by30 \advance \xpos
by\tempcounta \advance \arrowlength by-\tempcounta \advance
\arrowlength by-\tempcountb
\putvector(\xpos,\ypos)(1,0)\arrowlength\arrowtype \divide
\arrowlength by2 \advance \xpos by\arrowlength
\vertsize{\tempcounta}{#5}%
\divide\tempcounta by2 \advance \tempcounta by20
\if a#8 %
   \advance \ypos by\tempcounta
   \putbox(\xpos,\ypos){#5}%
\else
   \advance \ypos by-\tempcounta
   \putbox(\xpos,\ypos){#5}%
\fi}}

\def\putvmorphism(#1,#2)[#3`#4`#5]#6#7#8{{%
\xpos #1 \ypos #2 \arrowlength #6 \arrowtype #7
\settowidth{\xlen}{$#5$}%
\putbox(\xpos,\ypos){#3}%
{\advance \ypos by-\arrowlength
\putbox(\xpos,\ypos){#4}}%
{\advance\arrowlength by-140 \advance \ypos by-70 \ifdim\xlen>0pt
   \if m#8%
      \putsplitvector(\xpos,\ypos)\arrowlength\arrowtype
   \else
   \putvector(\xpos,\ypos)(0,-1)\arrowlength\arrowtype
   \fi
\else
   \putvector(\xpos,\ypos)(0,-1)\arrowlength\arrowtype
\fi}%
\ifdim\xlen>0pt
   \divide \arrowlength by2
   \advance\ypos by-\arrowlength
   \if l#8%
      \advance \xpos by-40
      \putrbox(\xpos,\ypos){#5}%
   \else\if r#8%
      \advance \xpos by40
      \putlbox(\xpos,\ypos){#5}%
   \else
      \putbox(\xpos,\ypos){#5}%
   \fi\fi
\fi }}

\def\putsquarep<#1>(#2)[#3;#4`#5`#6`#7]{{%
\setsqparms[#1]%
\setpos(#2)%
\settokens`#3`%
\puthmorphism(\xpos,\ypos)[\tokenc`\tokend`{#7}]{\width}{\arrowtyped}b%
\advance\ypos by \height
\puthmorphism(\xpos,\ypos)[\tokena`\tokenb`{#4}]{\width}{\arrowtypea}a%
\putvmorphism(\xpos,\ypos)[``{#5}]{\height}{\arrowtypeb}l%
\advance\xpos by \width
\putvmorphism(\xpos,\ypos)[``{#6}]{\height}{\arrowtypec}r%
}}

\def\putsquare{\@ifnextchar <{\putsquarep}{\putsquarep%
   <\arrowtypea`\arrowtypeb`\arrowtypec`\arrowtyped;\width`\height>}}
\def\square{\@ifnextchar< {\squarep}{\squarep
   <\arrowtypea`\arrowtypeb`\arrowtypec`\arrowtyped;\width`\height>}}
\def\squarep<#1>[#2`#3`#4`#5;#6`#7`#8`#9]{{
\setsqparms[#1]
\diagram
\putsquarep<\arrowtypea`\arrowtypeb`\arrowtypec`
\arrowtyped;\width`\height>
(0,0)[#2`#3`#4`{#5};#6`#7`#8`{#9}]
\enddiagram
}}                                                 
\def\putptrianglep<#1>(#2,#3)[#4`#5`#6;#7`#8`#9]{{%
\settriparms[#1]%
\xpos=#2 \ypos=#3 \advance\ypos by \height
\puthmorphism(\xpos,\ypos)[#4`#5`{#7}]{\height}{\arrowtypea}a%
\putvmorphism(\xpos,\ypos)[`#6`{#8}]{\height}{\arrowtypeb}l%
\advance\xpos by\height
\putmorphism(\xpos,\ypos)(-1,-1)[``{#9}]{\height}{\arrowtypec}r%
}}

\def\putptriangle{\@ifnextchar <{\putptrianglep}{\putptrianglep
   <\arrowtypea`\arrowtypeb`\arrowtypec;\height>}}
\def\ptriangle{\@ifnextchar <{\ptrianglep}{\ptrianglep
   <\arrowtypea`\arrowtypeb`\arrowtypec;\height>}}
\def\ptrianglep<#1>[#2`#3`#4;#5`#6`#7]{{
\settriparms[#1]
\diagram
\putptrianglep<\arrowtypea`\arrowtypeb`
\arrowtypec;\height>
(0,0)[#2`#3`#4;#5`#6`{#7}]
\enddiagram
}}                                            

\def\putqtrianglep<#1>(#2,#3)[#4`#5`#6;#7`#8`#9]{{%
\settriparms[#1]%
\xpos=#2 \ypos=#3 \advance\ypos by\height
\puthmorphism(\xpos,\ypos)[#4`#5`{#7}]{\height}{\arrowtypea}a%
\putmorphism(\xpos,\ypos)(1,-1)[``{#8}]{\height}{\arrowtypeb}l%
\advance\xpos by\height
\putvmorphism(\xpos,\ypos)[`#6`{#9}]{\height}{\arrowtypec}r%
}}

\def\putqtriangle{\@ifnextchar <{\putqtrianglep}{\putqtrianglep
   <\arrowtypea`\arrowtypeb`\arrowtypec;\height>}}
\def\qtriangle{\@ifnextchar <{\qtrianglep}{\qtrianglep
   <\arrowtypea`\arrowtypeb`\arrowtypec;\height>}}
\def\qtrianglep<#1>[#2`#3`#4;#5`#6`#7]{{
\settriparms[#1]
\width=\height                                
\diagram
\putqtrianglep<\arrowtypea`\arrowtypeb`
\arrowtypec;\height>
(0,0)[#2`#3`#4;#5`#6`{#7}]
\enddiagram
}}

\def\putdtrianglep<#1>(#2,#3)[#4`#5`#6;#7`#8`#9]{{%
\settriparms[#1]%
\xpos=#2 \ypos=#3
\puthmorphism(\xpos,\ypos)[#5`#6`{#9}]{\height}{\arrowtypec}b%
\advance\xpos by \height \advance\ypos by\height
\putmorphism(\xpos,\ypos)(-1,-1)[``{#7}]{\height}{\arrowtypea}l%
\putvmorphism(\xpos,\ypos)[#4``{#8}]{\height}{\arrowtypeb}r%
}}

\def\putdtriangle{\@ifnextchar <{\putdtrianglep}{\putdtrianglep
   <\arrowtypea`\arrowtypeb`\arrowtypec;\height>}}
\def\dtriangle{\@ifnextchar <{\dtrianglep}{\dtrianglep
   <\arrowtypea`\arrowtypeb`\arrowtypec;\height>}}
\def\dtrianglep<#1>[#2`#3`#4;#5`#6`#7]{{
\settriparms[#1]
\width=\height                                
\diagram
\putdtrianglep<\arrowtypea`\arrowtypeb`
\arrowtypec;\height>
(0,0)[#2`#3`#4;#5`#6`{#7}]
\enddiagram
}}

\def\putbtrianglep<#1>(#2,#3)[#4`#5`#6;#7`#8`#9]{{%
\settriparms[#1]%
\xpos=#2 \ypos=#3
\puthmorphism(\xpos,\ypos)[#5`#6`{#9}]{\height}{\arrowtypec}b%
\advance\ypos by\height
\putmorphism(\xpos,\ypos)(1,-1)[``{#8}]{\height}{\arrowtypeb}r%
\putvmorphism(\xpos,\ypos)[#4``{#7}]{\height}{\arrowtypea}l%
}}

\def\putbtriangle{\@ifnextchar <{\putbtrianglep}{\putbtrianglep
   <\arrowtypea`\arrowtypeb`\arrowtypec;\height>}}
\def\btriangle{\@ifnextchar <{\btrianglep}{\btrianglep
   <\arrowtypea`\arrowtypeb`\arrowtypec;\height>}}
\def\btrianglep<#1>[#2`#3`#4;#5`#6`#7]{{
\settriparms[#1]
\width=\height                               
\diagram
\putbtrianglep<\arrowtypea`\arrowtypeb`
\arrowtypec;\height>
(0,0)[#2`#3`#4;#5`#6`{#7}]
\enddiagram
}}

\def\putAtrianglep<#1>(#2,#3)[#4`#5`#6;#7`#8`#9]{{%
\settriparms[#1]%
\xpos=#2 \ypos=#3 {\multiply \height by2
\puthmorphism(\xpos,\ypos)[#5`#6`{#9}]{\height}{\arrowtypec}b}%
\advance\xpos by\height \advance\ypos by\height
\putmorphism(\xpos,\ypos)(-1,-1)[#4``{#7}]{\height}{\arrowtypea}l%
\putmorphism(\xpos,\ypos)(1,-1)[``{#8}]{\height}{\arrowtypeb}r%
}}

\def\putAtriangle{\@ifnextchar <{\putAtrianglep}{\putAtrianglep
   <\arrowtypea`\arrowtypeb`\arrowtypec;\height>}}
\def\Atriangle{\@ifnextchar <{\Atrianglep}{\Atrianglep
   <\arrowtypea`\arrowtypeb`\arrowtypec;\height>}}
\def\Atrianglep<#1>[#2`#3`#4;#5`#6`#7]{{
\settriparms[#1]
\width=\height                                     
\diagram
\putAtrianglep<\arrowtypea`\arrowtypeb`
\arrowtypec;\height>
(0,0)[#2`#3`#4;#5`#6`{#7}]
\enddiagram
}}

\def\putAtrianglepairp<#1>(#2)[#3;#4`#5`#6`#7`#8]{{%
\settripairparms[#1]%
\setpos(#2)%
\settokens`#3`%
\puthmorphism(\xpos,\ypos)[\tokenb`\tokenc`{#7}]{\height}{\arrowtyped}b%
\advance\xpos by\height
\puthmorphism(\xpos,\ypos)[\phantom{\tokenc}`\tokend`{#8}]%
{\height}{\arrowtypee}b%
\advance\ypos by\height
\putmorphism(\xpos,\ypos)(-1,-1)[\tokena``{#4}]{\height}{\arrowtypea}l%
\putvmorphism(\xpos,\ypos)[``{#5}]{\height}{\arrowtypeb}m%
\putmorphism(\xpos,\ypos)(1,-1)[``{#6}]{\height}{\arrowtypec}r%
}}

\def\putAtrianglepair{\@ifnextchar <{\putAtrianglepairp}{\putAtrianglepairp%
   <\arrowtypea`\arrowtypeb`\arrowtypec`\arrowtyped`\arrowtypee;\height>}}
\def\Atrianglepair{\@ifnextchar <{\Atrianglepairp}{\Atrianglepairp%
   <\arrowtypea`\arrowtypeb`\arrowtypec`\arrowtyped`\arrowtypee;\height>}}

\def\Atrianglepairp<#1>[#2;#3`#4`#5`#6`#7]{{
\settripairparms[#1]
\settokens`#2`
\width=\height                                
\diagram
\putAtrianglepairp                            
<\arrowtypea`\arrowtypeb`\arrowtypec`
\arrowtyped`\arrowtypee;\height>
(0,0)[{#2};#3`#4`#5`#6`{#7}]
\enddiagram
}}

\def\putVtrianglep<#1>(#2,#3)[#4`#5`#6;#7`#8`#9]{{%
\settriparms[#1]%
\xpos=#2 \ypos=#3 \advance\ypos by\height {\multiply\height by2
\puthmorphism(\xpos,\ypos)[#4`#5`{#7}]{\height}{\arrowtypea}a}%
\putmorphism(\xpos,\ypos)(1,-1)[`#6`{#8}]{\height}{\arrowtypeb}l%
\advance\xpos by\height \advance\xpos by\height
\putmorphism(\xpos,\ypos)(-1,-1)[``{#9}]{\height}{\arrowtypec}r%
}}

\def\putVtriangle{\@ifnextchar <{\putVtrianglep}{\putVtrianglep
   <\arrowtypea`\arrowtypeb`\arrowtypec;\height>}}
\def\Vtriangle{\@ifnextchar <{\Vtrianglep}{\Vtrianglep
   <\arrowtypea`\arrowtypeb`\arrowtypec;\height>}}
\def\Vtrianglep<#1>[#2`#3`#4;#5`#6`#7]{{
\settriparms[#1]
\width=\height                                 
\diagram
\putVtrianglep<\arrowtypea`\arrowtypeb`
\arrowtypec;\height>
(0,0)[#2`#3`#4;#5`#6`{#7}]
\enddiagram
}}

\def\putVtrianglepairp<#1>(#2)[#3;#4`#5`#6`#7`#8]{{
\settripairparms[#1]%
\setpos(#2)%
\settokens`#3`%
\advance\ypos by\height
\putmorphism(\xpos,\ypos)(1,-1)[`\tokend`{#6}]{\height}{\arrowtypec}l%
\puthmorphism(\xpos,\ypos)[\tokena`\tokenb`{#4}]{\height}{\arrowtypea}a%
\advance\xpos by\height
\puthmorphism(\xpos,\ypos)[\phantom{\tokenb}`\tokenc`{#5}]%
{\height}{\arrowtypeb}a%
\putvmorphism(\xpos,\ypos)[``{#7}]{\height}{\arrowtyped}m%
\advance\xpos by\height
\putmorphism(\xpos,\ypos)(-1,-1)[``{#8}]{\height}{\arrowtypee}r%
}}

\def\putVtrianglepair{\@ifnextchar <{\putVtrianglepairp}{\putVtrianglepairp%
    <\arrowtypea`\arrowtypeb`\arrowtypec`\arrowtyped`\arrowtypee;\height>}}
\def\Vtrianglepair{\@ifnextchar <{\Vtrianglepairp}{\Vtrianglepairp%
    <\arrowtypea`\arrowtypeb`\arrowtypec`\arrowtyped`\arrowtypee;\height>}}
\def\Vtrianglepairp<#1>[#2;#3`#4`#5`#6`#7]{{
\settripairparms[#1]
\settokens`#2`
\diagram
\putVtrianglepairp                             
<\arrowtypea`\arrowtypeb`\arrowtypec`
\arrowtyped`\arrowtypee;\height>
(0,0)[{#2};#3`#4`#5`#6`{#7}]
\enddiagram
}}

\def\putCtrianglep<#1>(#2,#3)[#4`#5`#6;#7`#8`#9]{{%
\settriparms[#1]%
\xpos=#2 \ypos=#3 \advance\ypos by\height
\putmorphism(\xpos,\ypos)(1,-1)[``{#9}]{\height}{\arrowtypec}l%
\advance\xpos by\height \advance\ypos by\height
\putmorphism(\xpos,\ypos)(-1,-1)[#4`#5`{#7}]{\height}{\arrowtypea}l%
{\multiply\height by 2
\putvmorphism(\xpos,\ypos)[`#6`{#8}]{\height}{\arrowtypeb}r}%
}}

\def\putCtriangle{\@ifnextchar <{\putCtrianglep}{\putCtrianglep
    <\arrowtypea`\arrowtypeb`\arrowtypec;\height>}}
\def\Ctriangle{\@ifnextchar <{\Ctrianglep}{\Ctrianglep
    <\arrowtypea`\arrowtypeb`\arrowtypec;\height>}}
\def\Ctrianglep<#1>[#2`#3`#4;#5`#6`#7]{{
\settriparms[#1]
\width=\height                               
\diagram
\putCtrianglep<\arrowtypea`\arrowtypeb`
\arrowtypec;\height>
(0,0)[#2`#3`#4;#5`#6`{#7}]
\enddiagram
}}                                           
\def\putDtrianglep<#1>(#2,#3)[#4`#5`#6;#7`#8`#9]{{%
\settriparms[#1]%
\xpos=#2 \ypos=#3 \advance\xpos by\height \advance\ypos by\height
\putmorphism(\xpos,\ypos)(-1,-1)[``{#9}]{\height}{\arrowtypec}r%
\advance\xpos by-\height \advance\ypos by\height
\putmorphism(\xpos,\ypos)(1,-1)[`#5`{#8}]{\height}{\arrowtypeb}r%
{\multiply\height by 2
\putvmorphism(\xpos,\ypos)[#4`#6`{#7}]{\height}{\arrowtypea}l}%
}}

\def\putDtriangle{\@ifnextchar <{\putDtrianglep}{\putDtrianglep
    <\arrowtypea`\arrowtypeb`\arrowtypec;\height>}}
\def\Dtriangle{\@ifnextchar <{\Dtrianglep}{\Dtrianglep
   <\arrowtypea`\arrowtypeb`\arrowtypec;\height>}}
\def\Dtrianglep<#1>[#2`#3`#4;#5`#6`#7]{{
\settriparms[#1]
\width=\height                              
\diagram
\putDtrianglep<\arrowtypea`\arrowtypeb`
\arrowtypec;\height>
(0,0)[#2`#3`#4;#5`#6`{#7}]
\enddiagram
}}                                          
\def\setrecparms[#1`#2]{\width=#1 \height=#2}%

\def\recursep<#1`#2>[#3;#4`#5`#6`#7`#8]{{\m@th
\width=#1 \height=#2 \settokens`#3`
\settowidth{\tempdimen}{$\tokena$} \ifdim\tempdimen=0pt
  \savebox{\tempboxa}{\hbox{$\tokenb$}}%
  \savebox{\tempboxb}{\hbox{$\tokend$}}%
  \savebox{\tempboxc}{\hbox{$#6$}}%
\else
  \savebox{\tempboxa}{\hbox{$\hbox{$\tokena$}\times\hbox{$\tokenb$}$}}%
  \savebox{\tempboxb}{\hbox{$\hbox{$\tokena$}\times\hbox{$\tokend$}$}}%
  \savebox{\tempboxc}{\hbox{$\hbox{$\tokena$}\times\hbox{$#6$}$}}%
\fi \ypos=\height \divide\ypos by 2 \xpos=\ypos \advance\xpos by
\width \bfig
\putCtrianglep<-1`1`1;\ypos>(0,0)[`\tokenc`;#5`#6`{#7}]%
\puthmorphism(\ypos,0)[\tokend`\usebox{\tempboxb}`{#8}]{\width}{-1}b%
\puthmorphism(\ypos,\height)[\tokenb`\usebox{\tempboxa}`{#4}]{\width}{-1}a%
\advance\ypos by \width
\putvmorphism(\ypos,\height)[``\usebox{\tempboxc}]{\height}1r%
\efig }}

\def\recurse{\@ifnextchar <{\recursep}{\recursep<\width`\height>}}

\def\puttwohmorphisms(#1,#2)[#3`#4;#5`#6]#7#8#9{{%
\puthmorphism(#1,#2)[#3`#4`]{#7}0a \ypos=#2 \advance\ypos by 20
\puthmorphism(#1,\ypos)[\phantom{#3}`\phantom{#4}`#5]{#7}{#8}a
\advance\ypos by -40
\puthmorphism(#1,\ypos)[\phantom{#3}`\phantom{#4}`#6]{#7}{#9}b }}

\def\puttwovmorphisms(#1,#2)[#3`#4;#5`#6]#7#8#9{{%
\putvmorphism(#1,#2)[#3`#4`]{#7}0a \xpos=#1 \advance\xpos by -20
\putvmorphism(\xpos,#2)[\phantom{#3}`\phantom{#4}`#5]{#7}{#8}l
\advance\xpos by 40
\putvmorphism(\xpos,#2)[\phantom{#3}`\phantom{#4}`#6]{#7}{#9}r }}

\def\puthcoequalizer(#1)[#2`#3`#4;#5`#6`#7]#8#9{{%
\setpos(#1)%
\puttwohmorphisms(\xpos,\ypos)[#2`#3;#5`#6]{#8}11%
\advance\xpos by #8
\puthmorphism(\xpos,\ypos)[\phantom{#3}`#4`#7]{#8}1{#9} }}

\def\putvcoequalizer(#1)[#2`#3`#4;#5`#6`#7]#8#9{{%
\setpos(#1)%
\puttwovmorphisms(\xpos,\ypos)[#2`#3;#5`#6]{#8}11%
\advance\ypos by -#8
\putvmorphism(\xpos,\ypos)[\phantom{#3}`#4`#7]{#8}1{#9} }}

\def\putthreehmorphisms(#1)[#2`#3;#4`#5`#6]#7(#8)#9{{%
\setpos(#1) \settypes(#8)
\if a#9 %
     \vertsize{\tempcounta}{#5}%
     \vertsize{\tempcountb}{#6}%
     \ifnum \tempcounta<\tempcountb \tempcounta=\tempcountb \fi
\else
     \vertsize{\tempcounta}{#4}%
     \vertsize{\tempcountb}{#5}%
     \ifnum \tempcounta<\tempcountb \tempcounta=\tempcountb \fi
\fi \advance \tempcounta by 60
\puthmorphism(\xpos,\ypos)[#2`#3`#5]{#7}{\arrowtypeb}{#9}
\advance\ypos by \tempcounta
\puthmorphism(\xpos,\ypos)[\phantom{#2}`\phantom{#3}`#4]{#7}{\arrowtypea}{#9}
\advance\ypos by -\tempcounta \advance\ypos by -\tempcounta
\puthmorphism(\xpos,\ypos)[\phantom{#2}`\phantom{#3}`#6]{#7}{\arrowtypec}{#9}
}}

\def\setarrowtoks[#1`#2`#3`#4`#5`#6]{%
\def\toka{#1}
\def\tokb{#2}
\def\tokc{#3}
\def\tokd{#4}
\def\toke{#5}
\def\tokf{#6}
}
\def\hex{\@ifnextchar <{\hexp}{\hexp<1000`400>}}
\def\hexp<#1`#2>[#3`#4`#5`#6`#7`#8;#9]{%
\setarrowtoks[#9] \yext=#2 \advance \yext by #2 \xext=#1
\advance\xext by \yext \bfig
\putCtriangle<-1`0`1;#2>(0,0)[`#5`;\tokb``\tokd] \xext=#1
\yext=#2 \advance \yext by #2
\putsquare<1`0`0`1;\xext`\yext>(#2,0)[#3`#4`#7`#8;\toka```\tokf]
\advance \xext by #2
\putDtriangle<0`1`-1;#2>(\xext,0)[`#6`;`\tokc`\toke] \efig }

\makeatother

\input{tcilatex}
\begin{document}

\title{Lecture Notes in Lie Groups}
\author{Vladimir G. Ivancevic\thanks{%
Land Operations Division, Defence Science \& Technology Organisation, P.O.
Box 1500, Edinburgh SA 5111, Australia \quad\quad (e-mail:
~Vladimir.Ivancevic@dsto.defence.gov.au)} \and Tijana T. Ivancevic\thanks{%
Tesla Science Evolution Institute \& QLIWW IP Pty Ltd., Adelaide, Australia
\quad\quad (e-mail: ~tijana.ivancevic@alumni.adelaide.edu.au)}}
\date{}
\maketitle

\begin{abstract}
These lecture notes in Lie Groups are designed for a 1--semester third year or graduate
course in mathematics, physics, engineering, chemistry or
biology. This landmark theory of the 20th Century mathematics and physics
gives a rigorous foundation to modern dynamics, as well as field and gauge
theories in physics, engineering and biomechanics. We give both physical and medical examples of Lie groups. 
The only necessary background for comprehensive reading of these notes are advanced calculus and linear
algebra.
\end{abstract}

\tableofcontents

\section{Preliminaries: Sets, Maps and Diagrams}

\subsection{Sets}

Given a map (or, a function) $f:A\rightarrow B$, the set $A$ is called the
\emph{domain} of $f$, and denoted $\limfunc{Dom}f$. The set $B$ is called
the \emph{codomain} of $f$, and denoted $\limfunc{Cod}f.$ The codomain is
not to be confused with the \emph{range} of $f(A)$, which is in general only
a subset of $B$ (see \cite{GaneshSprBig,GaneshADG}).

A map $f:X\rightarrow Y$ is called \emph{injective}, or 1--1, or an \emph{%
injection}, iff for every $y$ in the codomain $Y$ there is \emph{at most}
one $x$ in the domain $X$ with $f(x)=y$. Put another way, given $x$ and $%
x^{\prime }$ in $X$, if $f(x)=f(x^{\prime })$, then it follows that $%
x=x^{\prime }$. A map $f:X\rightarrow Y$ is called \emph{surjective}, or
\emph{onto}, or a \emph{surjection}, iff for every $y $ in the codomain $%
\limfunc{Cod}f$ there is \emph{at least} one $x$ in the \emph{domain} $X$
with $f(x)=y$. Put another way, the \emph{range} $f(X)$ is equal to the
codomain $Y$. A map is \emph{bijective} iff it is both injective and
surjective. Injective functions are called \emph{monomorphisms}, and
surjective functions are called \emph{epimorphisms} in the \emph{category of
sets} (see below). Bijective functions are called \emph{isomorphisms}.

A \emph{relation} is any subset of a \emph{Cartesian product} (see below).
By definition, an \emph{equivalence relation} $\alpha$ on a set $X$ is a
relation which is \emph{reflexive, symmetrical} and \emph{transitive}, i.e.,
relation that satisfies the following three conditions:

\begin{enumerate}
\item \emph{Reflexivity}: each element $x\in X$ is equivalent to itself, i.e.%
$,$ $x\alpha x$;

\item \emph{Symmetry}: for any two elements $a,b\in X$, $a\alpha b$ implies $%
b\alpha a$; ~and

\item \emph{Transitivity}: $a\alpha b$\ and $b\alpha c$ implies $a\alpha c$.
\end{enumerate}

Similarly, a relation $\leq $ defines a \emph{partial order} on a set $S$ if
it has the following properties:

\begin{enumerate}
\item \emph{Reflexivity}: $a\leq a$\ for all $a\in S$;

\item \emph{Antisymmetry}: $a\leq b$\ and $b\leq a$ implies $a=b$; ~and

\item \emph{Transitivity}: $a\leq b$\ and $b\leq c$ implies $a\leq c$.
\end{enumerate}

A \emph{partially ordered set} (or \emph{poset}) is a set taken together
with a partial order on it. Formally, a partially ordered set is defined as
an ordered pair $P=(X,\leq )$, where $X$ is called the \emph{ground set} of $%
P$ and $\leq $\ is the partial order of $P$.

\subsection{Maps}

Let $f$ and $g$ be maps with domains $A$ and $B$. Then the maps $f+g$, $f-g$%
, $fg$, and $f/g$ are defined as follows (see \cite{GaneshSprBig,GaneshADG}%
):
\begin{eqnarray*}
(f+g)(x) &=&f(x)+g(x)\text{ \ \ \ \ \ \ \ \ \ \ \ domain }=A\cap B, \\
(f-g)(x) &=&f(x)-g(x)\text{ \ \ \ \ \ \ \ \ \ \ \ domain }=A\cap B, \\
(fg)(x) &=&f(x)\,g(x)\text{ \ \ \ \ \ \ \ \ \ \ \ domain }=A\cap B, \\
\left( \frac{f}{g}\right) (x) &=&\frac{f(x)}{g(x)}\text{ \ \ \ \ \ \ \ \ \ \
\ domain }=\{x\in A\cap B:g(x)\neq 0\}.
\end{eqnarray*}

Given two maps $f$ and $g$, the composite map $f\circ g$, called the \emph{%
composition} of $f$ and $g$, is defined by
\begin{equation*}
(f\circ g)(x)=f(g(x)).
\end{equation*}
The $(f\circ g)-$machine is composed of the $g-$machine (first) and then the
$f-$machine,
\begin{equation*}
x\rightarrow [[g]]\rightarrow g(x)\rightarrow [[f]]\rightarrow f(g(x)).
\end{equation*}
For example, suppose that $y=f(u)=\sqrt{u}$ and $u=g(x)=x^{2}+1$. Since $y$
is a function of $u$ and $u$ is a function of $x$, it follows that $y$ is
ultimately a function of $x$. We calculate this by substitution
\begin{equation*}
y=f(u)=f\circ g=f(g(x))=f(x^{2}+1)=\sqrt{x^{2}+1}.
\end{equation*}

If $f$ and $g$ are both differentiable (or smooth, i.e., $C^{\infty}$) maps
and $h=f\circ g$ is the composite map defined by $h(x)=f(g(x))$, then $h$ is
differentiable and $h^{\prime }$ is given by the product:
\begin{equation*}
h^{\prime }(x)=f^{\prime }(g(x))\,g^{\prime }(x).
\end{equation*}
In Leibniz notation, if $y=f(u)$ and $u=g(x)$ are both differentiable maps,
then
\begin{equation*}
\frac{dy}{dx}=\frac{dy}{du}\frac{du}{dx}.
\end{equation*}
The reason for the name \emph{chain rule} becomes clear if we add another
link to the chain. Suppose that we have one more differentiable map $x=h(t)$%
. Then, to calculate the derivative of $y$ with respect to $t$, we use the
chain rule twice,
\begin{equation*}
\frac{dy}{dt}=\frac{dy}{du}\frac{du}{dx}\frac{dx}{dt}.
\end{equation*}

Given a 1--1 continuous (i.e., $C^{0}$) map $F$ with a nonzero \emph{Jacobian%
} $\left| \frac{\partial (x,...)}{\partial (u,...)}\right|$\ that maps a
region $S$ onto a region $R$, we have the following substitution formulas:%
\newline
\smallskip

\noindent 1. For a single integral,
\begin{equation*}
\int_{R}f(x)\,dx=\int_{S}f(x(u))\frac{\partial x}{\partial u}du;
\end{equation*}
2. For a double integral,
\begin{equation*}
\iint_{R}f(x,y)\,dA=\iint_{S}f(x(u,v),y(u,v))\left| \frac{\partial (x,y)}{%
\partial (u,v)}\right| dudv;
\end{equation*}
3. For a triple integral,
\begin{equation*}
\iiint_{R}f(x,y,z)\,dV=\iiint_{S}f(x(u,v,w),y(u,v,w),z(u,v,w))\left| \frac{%
\partial (x,y,z)}{\partial (u,v,w)}\right| dudvdw;
\end{equation*}
4. Generalization to $n-$tuple integrals is obvious.

\subsection{Commutative Diagrams}

Many properties of mathematical systems can be unified and simplified by a
presentation with \emph{commutative diagrams of arrows}. Each arrow $%
f:X\rightarrow Y$ represents a function (i.e., a map, transformation,
operator); that is, a source (domain) set $X$, a target (codomain) set $Y$,
and a rule $x\mapsto f(x)$ which assigns to each element $x\in X$ an element
$f(x)\in Y$. A typical diagram of sets and functions is (see \cite%
{GaneshSprBig,GaneshADG}):
\begin{equation*}
\qtriangle[X`Y`Z;f`h`g\qquad\qquad \text{or} \qquad\qquad] %
\qtriangle[X`f(X)`g(f(X));f`h`g]
\end{equation*}
This diagram is \emph{commutative} iff $h=g\circ f$, where $g\circ f$ is the
usual composite function $g\circ f:X\rightarrow Z$, defined by $x\mapsto
g(f(x))$.

Less formally, composing maps is like following directed paths from one
object to another (e.g., from set to set). In general, a diagram is
commutative iff any two paths along arrows that start at the same point and
finish at the same point yield the same `homomorphism' via compositions
along successive arrows. Commutativity of the whole diagram follows from
commutativity of its triangular components. Study of commutative diagrams is
popularly called `diagram chasing', and provides a powerful tool for
mathematical thought.

Many properties of mathematical constructions may be represented by \emph{%
universal properties} of diagrams. Consider the \emph{Cartesian product} $%
X\times Y$ of two sets, consisting as usual of all ordered pairs $\langle
x,y\rangle $ of elements $x\in X$ and $y\in Y$. The projections $\langle
x,y\rangle \mapsto x,\,\,\langle x,y\rangle \mapsto y$ of the product on its
`axes' $X$ and $Y$ are functions $p:X\times Y\rightarrow X,\,\,q:X\times
Y\rightarrow Y$. Any function $h:W\rightarrow X\times Y$ from a third set $W$
is uniquely determined by its composites $p\circ h$ and $q\circ h$.
Conversely, given $W$ and two functions $f$ and $g$ as in the diagram below,
there is a unique function $h$ which makes the following diagram commute:
\begin{equation*}
\Atrianglepair<1`1`1`-1`1;>[W`X`X\times Y`Y;f`h`g`p`q]
\end{equation*}%
This property describes the Cartesian product $X\times Y$ uniquely; the same
diagram, read in the category of topological spaces or of groups, describes
uniquely the Cartesian product of spaces or of the direct product of groups.

\section{Groups}

A \emph{group} is a pointed set $(G,e)$ with a \emph{multiplication} $\mu
:G\times G\rightarrow G$ and an \emph{inverse} $\nu :G\rightarrow G$ such
that the following diagrams commute (see \cite%
{Switzer,GaneshSprBig,GaneshADG}):

\begin{enumerate}
\item
\begin{equation*}
\Vtrianglepair<1`1`1`1`-1;>[G`G\times G`G`G;(e,1)`(1,e)`1`\mu `1]
\end{equation*}
($e$ is a two--sided identity)

\item
\begin{equation*}
\square <1`1`1`1;900`500>[G\times G\times G`G\times G`G\times G`G;\mu \times
1`1\times \mu `\mu `\mu ]
\end{equation*}
(associativity)

\item
\begin{equation*}
\Vtrianglepair<1`1`1`1`-1;>[G`G\times G`G`G;(\nu ,1)`(1,\nu )`e`\mu `e]
\end{equation*}
(inverse).
\end{enumerate}

Here $e:G\rightarrow G$ is the constant map $e(g)=e$ for all $g\in G$. $%
(e,1) $ means the map such that $(e,1)(g)=(e,g)$, etc. A group $G$ is called
\emph{commutative} or \emph{Abelian group} if in addition the following
diagram commutes
\begin{equation*}
\Vtriangle<1`1`1;>[G\times G`G\times G`G;T`\mu `\mu ]
\end{equation*}
where $T:G\times G\rightarrow G\times G$ is the switch map $%
T(g_{1},g_{2})=(g_{2},g_{1}),$ for all $(g_{1},g_{2})\in G\times G.$

A group $G$ \emph{acts} (on the left) on a set $A$ if there is a function $%
\alpha :G\times A\rightarrow A$ such that the following diagrams commute:

\begin{enumerate}
\item
\begin{equation*}
\qtriangle<1`1`1;>[A`G\times A`A;(e,1)`1`\alpha ]
\end{equation*}

\item
\begin{equation*}
\square <1`1`1`1;900`500>[G\times G\times A`G\times A`G\times A`A;1\times
\alpha `\mu \times 1`\alpha `\alpha ]
\end{equation*}%
where $(e,1)(x)=(e,x)$ for all $x\in A$. The \emph{orbit}\emph{s} of the
action are the sets $Gx=\{gx:g\in G\}$\ for all $x\in A$.
\end{enumerate}

Given two groups $(G,\ast )$ and $(H,\cdot )$, a \emph{group homomorphism}
from $(G,\ast )$ to $(H,\cdot )$ is a function $h:G\rightarrow H$ such that
for all $x$ and $y$ in $G$ it holds that\
\begin{equation*}
h(x\ast y)=h(x)\cdot h(y).
\end{equation*}%
From this property, one can deduce that $h$ maps the identity element $e_{G}$
of $G$ to the identity element $e_{H}$ of $H$, and it also maps inverses to
inverses in the sense that $h(x^{-1})=h(x)^{-1}$. Hence one can say that $h$
is \emph{compatible} with the \emph{group structure}.

The \emph{kernel} $\limfunc{Ker}h$\ of a group homomorphism $h:G\rightarrow
H $ consists of all those elements of $G$ which are sent by $h$ to the
identity element $e_{H}$ of $H$, i.e.,
\begin{equation*}
\limfunc{Ker}h=\{x\in G:h(x)=e_{H}\}.
\end{equation*}

The \emph{image} $\limfunc{Im}h$\ of a group homomorphism $h:G\rightarrow H$
consists of all elements of $G$ which are sent by $h$ to $H$, i.e.,
\begin{equation*}
\limfunc{Im}h=\{h(x):x\in G\}.
\end{equation*}

The kernel is a \emph{normal subgroup} of $G$ and the image is a \emph{%
subgroup} of $H$. The homomorphism $h$ is \emph{injective} (and called a
\emph{group monomorphism}) iff $\limfunc{Ker}h=e_{G}$, i.e., iff the kernel
of $h$ consists of the identity element of $G$ only.

\section{Manifolds}

A \textit{manifold} is an abstract mathematical space, which locally (i.e.,
in a close--up view) resembles the spaces described by \textit{Euclidean
geometry}, but which globally (i.e., when viewed as a whole) may have a more
complicated structure (see \cite{de Rham}). For example, the \textit{surface
of Earth} is a manifold; locally it seems to be flat, but viewed as a whole
from the outer space (globally) it is actually round. A manifold can be
constructed by `gluing' separate \textit{Euclidean spaces} together; for
example, a world map can be made by gluing many maps of local regions
together, and accounting for the resulting distortions.

As main pure--mathematical references for manifolds we recommend popular
graduate textbooks by two ex--\emph{Bourbaki} members, \emph{Serge Lang}
\cite{LangMan,LangDg} and Jean Dieudonne \cite{Dieudonne,Dieudonne2}.
Besides, the reader might wish to consult some other `classics', including
\cite{de Rham,Spivak,SpivakDG,Bruhat,Choquet,Bott,Abraham}. Finally, as
first--order applications, we recommend three popular textbooks in
mechanics, \cite{AbrahamMeh,Arnold,Marsden}, as well as our own geometrical
monographs \cite{GaneshSprBig,GaneshADG}.

Another example of a manifold is a \textit{circle} $S^{1}$. A small piece of
a circle appears to be like a slightly--bent part of a straight line
segment, but overall the circle and the segment are different 1D manifolds.
A circle can be formed by bending a straight line segment and gluing the
ends together.\footnote{%
Locally, the circle looks like a line. It is 1D, that is, only one
coordinate is needed to say where a point is on the circle locally.
Consider, for instance, the top part of the circle, where the $y-$coordinate
is positive. Any point in this part can be described by the $x-$coordinate.
So, there is a continuous \textit{bijection} $\chi _{top}$ (a mapping which
is 1--1 both ways), which maps the top part of the circle to the open
interval $(-1,1)$, by simply projecting onto the first coordinate: ~$\chi
_{top}(x,y)=x$.~ Such a function is called a \textit{chart}. Similarly,
there are charts for the bottom, left , and right parts of the circle.
Together, these parts \textit{cover} the whole circle and the four charts
form an \textit{atlas} (see the next subsection) for the circle. The top and
right charts overlap: their intersection lies in the quarter of the circle
where both the $x-$ and the $y-$coordinates are positive. The two charts $%
\chi _{top}$ and $\chi _{right}$ map this part bijectively to the interval $%
(0,1)$. Thus a function $T$ from $(0,1)$ to itself can be constructed, which
first inverts the top chart to reach the circle and then follows the right
chart back to the interval:
\begin{equation*}
T(a)=\chi _{\mathrm{right}}\left( \chi _{\mathrm{top}}^{-1}(a)\right) =\chi
_{\mathrm{right}}\left( a,\sqrt{1-a^{2}}\right) =\sqrt{1-a^{2}}.
\end{equation*}%
Such a function is called a \textit{transition map}. The top, bottom, left,
and right charts show that the circle is a manifold, but they do not form
the only possible atlas. Charts need not be geometric projections, and the
number of charts is a matter of choice. $T$ and the other transition
functions are differentiable on the interval $(0,1)$. Therefore, with this
atlas the circle is a \emph{differentiable,} or \emph{smooth manifold.}}

The surfaces of a \textit{sphere}\footnote{%
The surface of the sphere $S^{2}$ can be treated in almost the same way as
the circle $S^{1}$. It can be viewed as a subset of $\mathbb{R}^{3}$,
defined by: ~$S=\{(x,y,z)\in \mathbb{R}^{3}|x^{2}+y^{2}+z^{2}=1\}.$~ The
sphere is 2D, so each chart will map part of the sphere to an open subset of
$\mathbb{R}^{2}$. Consider the northern hemisphere, which is the part with
positive $z$ coordinate. The function $\chi $ defined by $\chi (x,y,z)=(x,y)$%
, maps the northern hemisphere to the open unit disc by projecting it on the
$(x,y)-$plane. A similar chart exists for the southern hemisphere. Together
with two charts projecting on the $(x,z)-$plane and two charts projecting on
the $(y,z)-$plane, an atlas of six charts is obtained which covers the
entire sphere. This can be easily generalized to an $n$D sphere $%
S^{n}=\{(x_{1},x_{2},...,x_{n})\in \mathbb{R}%
^{n}|x_{1}^{2}+x_{2}^{2}+...+x_{n}^{2}=1\}$.
\par
An $n-$sphere $S^{n}$ can be also constructed by gluing together two copies
of $\mathbb{R}^{n}$. The transition map between them is defined as $\mathbb{R%
}^{n}\setminus \{0\}\rightarrow \mathbb{R}^{n}\setminus \{0\}:x\mapsto
x/\Vert x\Vert ^{2}.$ This function is its own inverse, so it can be used in
both directions. As the transition map is a $(C^{\infty })-$\emph{smooth
function}, this atlas defines a \emph{smooth manifold}.} and a \textit{torus}%
\footnote{%
A torus (pl. tori), denoted by $T^{2}$, is a doughnut--shaped surface of
revolution generated by revolving a circle about an axis coplanar with the
circle. The sphere $S^{2}$ is a special case of the torus obtained when the
axis of rotation is a diameter of the circle. If the axis of rotation does
not intersect the circle, the torus has a hole in the middle and resembles a
ring doughnut, a hula hoop or an inflated tire. The other case, when the
axis of rotation is a chord of the circle, produces a sort of squashed
sphere resembling a round cushion.
\par
A torus can be defined parametrically by:
\begin{equation*}
x(u,v)=(R+r\cos {v})\cos {u},\qquad y(u,v)=(R+r\cos {v})\sin {u},\qquad
z(u,v)=r\sin {v},
\end{equation*}
where $u,v\in \lbrack 0,2\pi ],$ $R$ is the distance from the center of the
tube to the center of the torus, and $r$ is the radius of the tube.
According to a broader definition, the generator of a torus need not be a
circle but could also be an ellipse or any other conic section.
\par
Topologically, a torus is a closed surface defined as product of two
circles: $T^{2}=S^{1}\times S^{1}$. The surface described above, given the
relative topology from $\mathbb{R}^{3}$, is \emph{homeomorphic} to a
topological torus as long as it does not intersect its own axis.
\par
One can easily generalize the torus to arbitrary dimensions. An $n-$torus $%
T^{n}$ is defined as a product of $n$ circles: ~$T^{n}=S^{1}\times
S^{1}\times \cdots \times S^{1}$.~ Equivalently, the $n-$torus is obtained
from the $n-$cube (the $\mathbb{R}^{n}-$generalization of the ordinary cube
in $\mathbb{R}^{3}$) by gluing the opposite faces together.
\par
An $n-$torus $T^{n}$ is an example of an $n$D \textit{compact manifold}. It
is also an important example of a \textit{Lie group} (see below).} are
examples of 2D manifolds. Manifolds are important objects in mathematics,
physics and control theory, because they allow more complicated structures
to be expressed and understood in terms of the well--understood properties
of simpler Euclidean spaces (see \cite{GaneshADG}).

The \textit{Cartesian product} of manifolds is also a manifold (note that
not every manifold can be written as a product). The dimension of the
product manifold is the sum of the dimensions of its factors. Its topology
is the product topology, and a Cartesian product of charts is a chart for
the product manifold. Thus, an atlas for the product manifold can be
constructed using atlases for its factors. If these atlases define a
differential structure on the factors, the corresponding atlas defines a
differential structure on the product manifold. The same is true for any
other structure defined on the factors. If one of the factors has a
boundary, the product manifold also has a boundary. Cartesian products may
be used to construct tori and cylinders, for example, as $S^1 \times S^1$
and $S^1 \times [0,1]$, respectively.

Manifolds need not be \emph{connected} (all in `one piece'): a pair of
separate circles is also a \textit{topological manifold}(see below).
Manifolds need not be \emph{closed}: a line segment without its ends is a
manifold. Manifolds need not be \emph{finite}: a parabola is a topological
manifold.

Manifolds\footnote{%
{}
\par
Additional structures are often defined on manifolds. Examples of manifolds
with additional structure include:
\par
\begin{itemize}
\item \emph{differentiable} (or, \textit{smooth manifolds}, on which one can
do calculus;
\par
\item \textit{Riemannian manifolds}, on which \emph{distances} and \emph{%
angles} can be defined; they serve as the \textit{configuration spaces} in
mechanics;
\par
\item \textit{symplectic manifolds}, which serve as the \textit{phase spaces}
in mechanics and physics;
\par
\item 4D pseudo--Riemannian manifolds which model \textit{space--time} in
general relativity.
\end{itemize}
} can be viewed using either extrinsic or intrinsic view. In the \textit{%
extrinsic view}, usually used in geometry and topology of surfaces, an $n$D
manifold $M$ is seen as embedded in an $(n+1)$D Euclidean space $\mathbb{R}%
^{n+1}$. Such a manifold is called a `codimension 1 space'. With this view
it is easy to use intuition from Euclidean spaces to define additional
structure. For example, in a Euclidean space it is always clear whether a
vector at some point is tangential or normal to some surface through that
point. On the other hand, the \textit{intrinsic view} of an $n$D manifold $M$
is an abstract way of considering $M$ as a topological space by itself,
without any need for surrounding $(n+1)$D Euclidean space. This view is more
flexible and thus it is usually used in high--dimensional mechanics and
physics (where manifolds used represent configuration and phase spaces of
dynamical systems), can make it harder to imagine what a tangent vector
might be.

\subsection{Definition of a Manifold}

\begin{figure}[h]
\centerline{\includegraphics[width=9cm]{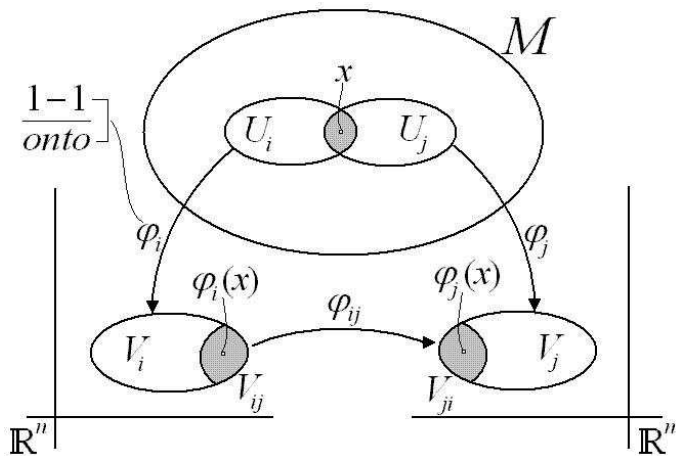}}
\caption{Geometric picture of the manifold concept.}
\label{Manifold1}
\end{figure}

Consider a set $M$ (see Figure \ref{Manifold1}) which is a \emph{candidate}
for a manifold. Any point $x\in M$ has its \textit{Euclidean chart}, given
by a 1--1 and \emph{onto} map $\varphi _{i}:M\rightarrow \mathbb{R}^{n}$,
with its \textit{Euclidean image} $V_{i}=\varphi _{i}(U_{i})$. More
precisely, a chart $\varphi _{i}$ is defined by (see \cite%
{GaneshSprBig,GaneshADG})
\begin{equation*}
\varphi _{i}:M\supset U_{i}\ni x\mapsto \varphi _{i}(x)\in V_{i}\subset
\mathbb{R}^{n},
\end{equation*}%
where $U_{i}\subset M$ and $V_{i}\subset \mathbb{R}^{n}$ are open sets.

Clearly, any point $x\in M$ can have several different charts (see Figure %
\ref{Manifold1}). Consider a case of two charts, $\varphi _{i},\varphi
_{j}:M\rightarrow \mathbb{R}^{n}$, having in their images two open sets, $%
V_{ij}=\varphi _{i}(U_{i}\cap U_{j})$ and $V_{ji}=\varphi _{j}(U_{i}\cap
U_{j})$. Then we have \textit{transition functions} $\varphi _{ij}$ between
them,
\begin{equation*}
\varphi _{ij}=\varphi _{j}\circ \varphi _{i}^{-1}:V_{ij}\rightarrow
V_{ji},\qquad \text{locally given by\qquad }\varphi _{ij}(x)=\varphi
_{j}(\varphi _{i}^{-1}(x)).
\end{equation*}
If transition functions $\varphi _{ij}$ exist, then we say that two charts, $%
\varphi _{i}$ and $\varphi _{j}$ are \emph{compatible}. Transition functions
represent a general (nonlinear) \emph{transformations of coordinates}, which
are the core of classical \emph{tensor calculus}.

A set of compatible charts $\varphi _{i}:M\rightarrow \mathbb{R}^{n},$ such
that each point $x\in M$ has its Euclidean image in at least one chart, is
called an \textit{atlas}. Two atlases are \emph{equivalent} iff all their
charts are compatible (i.e., transition functions exist between them), so
their union is also an atlas. A \textit{manifold structure} is a class of
equivalent atlases.

Finally, as charts $\varphi _{i}:M\rightarrow \mathbb{R}^{n}$ were supposed
to be 1-1 and onto maps, they can be either \emph{homeomorphism}\emph{s}, in
which case we have a \emph{topological} ($C^{0}$) manifold, or \emph{%
diffeomorphism}\emph{s}, in which case we have a \emph{smooth} ($C^{k}$)
manifold.

\subsection{Formal Definition of a Smooth Manifold}

Given a chart $( U,\varphi ) $, we call the set $U$ a \textit{coordinate
domain}, or a coordinate neighborhood of each of its points. If in addition $%
\varphi (U)$ is an open ball in $\mathbb{R}^{n}$, then $U$ is called a
\textit{coordinate ball}. The map $\varphi $ is called a (\emph{local})
\textit{coordinate map}, and the component functions $( x^{1},...,x^{n}) $
of $\varphi $, defined by $\varphi (m)=( x^{1}(m),...,x^{n}(m)) $, are
called \emph{local coordinates} on $U$ \cite{GaneshSprBig,GaneshADG}.

Two charts $( U_{1},\varphi _{1}) $ and $( U_{2},\varphi _{2}) $ such that $%
U_{1}\cap U_{2}\neq \varnothing $ are called \emph{compatible} if $\varphi
_{1}(U_{1}\cap U_{2})$ and $\varphi _{2}(U_{2}\cap U_{1})$ are open subsets
of $\mathbb{R}^{n}$. A family $( U_{\alpha },\varphi _{\alpha }) _{\alpha
\in A}$ of compatible charts on $M$ such that the $U_{\alpha }$ form a \emph{%
covering} of $M$ is called an \emph{atlas}. The maps $\varphi _{\alpha \beta
}=\varphi _{\beta }\circ \varphi _{\alpha }^{-1}:\varphi _{\alpha }(
U_{\alpha \beta }) \rightarrow \varphi _{\beta }( U_{\alpha \beta }) $ are
called the \emph{transition maps}, for the atlas\ $( U_{\alpha },\varphi
_{\alpha }) _{\alpha \in A},$ where $U_{\alpha \beta }=U_{\alpha }\cap
U_{\beta }$, so that we have a commutative triangle:
\begin{equation*}
\Atriangle<1`1`1;500>[U_{\alpha \beta }\subseteq M`\varphi _{\alpha }(
U_{\alpha \beta }) `\varphi _{\beta }( U_{\alpha \beta }) ;\varphi _{\alpha
}`\varphi _{\beta }`\varphi _{\alpha \beta }]
\end{equation*}

An atlas $( U_{\alpha },\varphi _{\alpha }) _{\alpha \in A}$ for a manifold $%
M$ is said to be a $C^{k}-$\emph{atlas}, if all transition maps $\varphi
_{\alpha \beta }:\varphi _{\alpha }( U_{\alpha \beta }) \rightarrow \varphi
_{\beta }( U_{\alpha \beta }) $ are of class $C^{k}$. Two $C^{k}$ atlases
are called $C^{k}-$\emph{equivalent}, if their union is again a $C^{k}-$%
atlas for $M$. An equivalence class of $C^{k}-$atlases is called a $C^{k}-$%
\emph{structure} on $M$. In other words, a smooth structure on $M$ is a
\emph{maximal} smooth atlas on $M$, i.e., such an atlas that is not
contained in any strictly larger smooth atlas. By a $C^{k}-$\emph{manifold} $%
M$, we mean a topological manifold together with a $C^{k}-$structure and a
chart on $M$ will be a chart belonging to some atlas of the $C^{k}-$%
structure. Smooth manifold means $C^{\infty }-$manifold, and the word `\emph{%
smooth}' is used synonymously with $C^{\infty }$.

Sometimes the terms `local coordinate system' or `parametrization' are used
instead of charts. That $M$ is not defined with any particular atlas, but
with an equivalence class of atlases, is a mathematical formulation of the
\emph{general covariance} principle. Every suitable coordinate system is
equally good. A Euclidean chart may well suffice for an open subset of $%
\mathbb{R}^{n}$, but this coordinate system is not to be preferred to the
others, which may require many charts (as with polar coordinates), but are
more convenient in other respects.

For example, the atlas of an $n-$sphere $S^{n}$ has two charts. If $%
N=(1,0,...,0)$ and $S=(-1,...,0,0)$ are the north and south poles of $S^{n}$
respectively, then the two charts are given by the stereographic projections
from $N$ and $S$:
\begin{eqnarray*}
\varphi _{1} &:&S^{n}\backslash \{N\}\rightarrow \mathbb{R}^{n},\varphi
_{1}(x^{1},...,x^{n+1})=(x^{2}/(1-x^{1}),\ldots ,x^{n+1}/(1-x^{1})),\;\,
\text{and} \\
\varphi _{2} &:&S^{n}\backslash \{S\}\rightarrow \mathbb{R}^{n},\varphi
_{2}(x^{1},...,x^{n+1})=(x^{2}/(1+x^{1}),\ldots ,x^{n+1}/(1+x^{1})),
\end{eqnarray*}
while the overlap map $\varphi _{2}\circ \varphi _{1}^{-1}:\mathbb{R}%
^{n}\backslash \{0\}\rightarrow \mathbb{R}^{n}\backslash \{0\}$ is given by
the diffeomorphism $(\varphi _{2}\circ \varphi _{1}^{-1})(z)=z/||z||^{2}$,
for $z$ in $\mathbb{R}^{n}\backslash \{0\}$, from $\mathbb{R}^{n}\backslash
\{0\}$ to itself.

Various \emph{additional structures} can be imposed on $\mathbb{R}^{n}$, and
the corresponding manifold $M$ will inherit them through its covering by
charts. For example, if a covering by charts takes their values in a \textit{%
Banach space} $E$, then $E$ is called the \textit{model space} and $M$ is
referred to as a $C^{k}-$\textit{Banach manifold} modelled on $E$.
Similarly, if a covering by charts takes their values in a \textit{Hilbert
space} $\mathcal{H}$, then $\mathcal{H}$ is called the \emph{model space}
and $M$ is referred to as a $C^{k}-$\textit{Hilbert manifold} modelled on $%
\mathcal{H}$. If not otherwise specified, we will consider $M$ to be an
Euclidean manifold, with its covering by charts taking their values in $%
\mathbb{R}^{n}$.

For a Hausdorff $C^{k}-$manifold the following properties are equivalent:
(i) it is paracompact; (ii) it is metrizable; (iii) it admits a Riemannian
metric;\footnote{%
Recall the corresponding properties of a \textit{Euclidean metric} $d$. For
any three points $x,y,z\in \mathbb{R}^{n}$, the following axioms are valid:
\begin{eqnarray*}
M_{1} &:&d(x,y)>0,\text{ \ for \ }x\neq y;\qquad\text{and}\qquad d(x,y)=0,%
\text{ \ for \ }x=y; \\
M_{2} &:&d(x,y)=d(y,x); \qquad\qquad\quad M_{3}~: d(x,y)\leq d(x,z)+d(z,y).
\end{eqnarray*}%
} (iv) each connected component is separable.

\subsection{Smooth Maps Between Smooth Manifolds}

A map $\varphi :M\rightarrow N$ between two manifolds $M$ and $N$, with $%
M\ni m\mapsto \varphi (m)\in N$, is called a \emph{smooth map}, or $C^{k}-$%
map, if we have the following charting \cite{GaneshSprBig,GaneshADG}:\newline
\bigbreak

\begin{equation*}
\bfig \putsquare<1`1`1`1;2000`1000>(0,0)[\odot`\odot`\odot`\odot;
\varphi`\phi`\psi`\psi \circ \varphi \circ \phi ^{-1}]
\put(50,1050){\oval(1000,600)}\put(1900,1050){\oval(1000,600)}
\put(-100,1050){\oval(450,400)}\put(2050,1050){\oval(450,400)}
\put(-150,1100){$U$}\put(-130,930){$m$} \put(2000,1100){$V$}\put(2020,930){$%
\varphi(m)$} \put(400,1100){$M$}\put(1450,1100){$N$}
\put(-20,-30){\oval(500,530)}\put(2000,-30){\oval(500,530)} \put(-220,100){$%
\phi(U)$}\put(2050,100){$\psi(V)$} \put(-120,-150){$\phi(m)$}\put(1830,-150){%
$\psi(\varphi(m))$} \putsquare<0`-1`0`1;1000`900>(-450,-450)[``\mathbb{R}%
^m`;```] \putsquare<0`0`-1`-1;1000`900>(1400,-450)[```\mathbb{R}^n;```] \efig
\end{equation*}

\noindent This means that for each $m\in M$ and each chart $\left( V,\psi
\right) $ on $N$ with $\varphi \left( m\right) \in V$ there is a chart $%
\left( U,\phi \right) $ on $M$ with $m\in U,\varphi \left( U\right)
\subseteq V$, and $\Phi =\psi \circ \varphi \circ \phi ^{-1}$ is $C^{k}$,
that is, the following diagram commutes:
\begin{equation*}
\square <1`1`1`1;1000`500>[M\supseteq U`V\subseteq N`\phi (U)`\psi
(V);\varphi `\phi `\psi `\Phi ]
\end{equation*}

Let $M$ and $N$ be smooth manifolds and let $\varphi :M\rightarrow N$ be a
smooth map. The map $\varphi $ is called a \emph{covering}, or equivalently,
$M$ is said to \emph{cover} $N$, if $\varphi $ is surjective and each point $%
n\in N$ admits an open neighborhood $V$ such that $\varphi ^{-1}(V)$ is a
union of disjoint open sets, each diffeomorphic via $\varphi $ to $V$.

A $C^{k}-$map $\varphi :M\rightarrow N$ is called a $C^{k}-$\textit{%
diffeomorphism} if $\varphi $ is a bijection, $\varphi ^{-1}:N\rightarrow M$
exists and is also $C^{k}$. Two manifolds are called diffeomorphic if there
exists a diffeomorphism between them. All smooth manifolds and smooth maps
between them form the category $\mathcal{M}$.

\subsection{Tangent Bundle and Lagrangian Dynamics}

The tangent bundle of a smooth $n-$manifold is the place where tangent
vectors live, and is itself a smooth $2n-$manifold. Vector--fields are
cross-sections of the tangent bundle. The \textit{Lagrangian} is a natural
energy function on the tangent bundle (see \cite{GaneshSprBig,GaneshADG}).

\begin{figure}[h]
\centerline{\includegraphics[width=4.5cm]{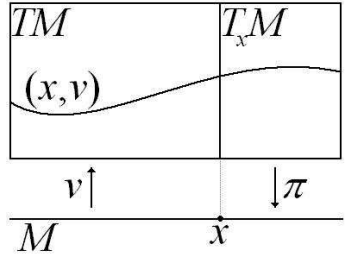}}
\caption{A sketch of a tangent bundle $TM$ of a smooth manifold $M$.}
\label{TBun1}
\end{figure}

In mechanics, to each $n$D \textit{configuration manifold} $M$ there is
associated its $2n$D \textit{velocity phase--space manifold}, denoted by $TM$
and called the tangent bundle of $M$ (see Figure \ref{TBun1}). The original
smooth manifold $M$ is called the \emph{base} of $TM$. There is an onto map $%
\pi:TM\rightarrow M$, called the \emph{projection}. Above each point $x\in M$
there is a \textit{tangent space} $T_x M=\pi^{-1}(x)$ to $M$ at $x$, which
is called a \textit{fibre}. The fibre $T_x M\subset TM$ is the subset of $TM$%
, such that the total tangent bundle, $TM=\dbigsqcup\limits_{m\in M}T_x M$,
is a \emph{disjoint union} of tangent spaces $T_x M$ to $M$ for all points $%
x\in M$. From dynamical perspective, the most important quantity in the
tangent bundle concept is the smooth map $v:M\rightarrow TM$, which is an
inverse to the projection $\pi$, i.e, $\pi\circ v=\func{Id}_M,\;\pi(v(x))=x$%
. It is called the \textit{velocity vector--field}. Its graph $(x,v(x))$
represents the \textit{cross--section} of the tangent bundle $TM$. This
explains the dynamical term \textit{velocity phase--space}, given to the
tangent bundle $TM$ of the manifold $M$.

If $[a,b]$ is a closed interval, a $C^{0}-$map $\gamma :[a,b]\rightarrow M$
is said to be \emph{differentiable} at the endpoint $a$ if there is a chart $%
( U,\phi ) $ at $\gamma (a)$ such that the following limit exists and is
finite:
\begin{equation}
{\frac{d}{dt}}(\phi \circ \gamma )(a)\equiv(\phi \circ \gamma )^{\prime
}(a)=\lim_{t\rightarrow a}\frac{(\phi \circ \gamma )(t)-(\phi \circ \gamma
)(a)}{t-a}.  \label{limit}
\end{equation}
Generalizing (\ref{limit}), we get the notion of the \emph{curve on a
manifold}. For a smooth manifold $M$ and a point $m\in M$ a curve at $m$ is
a $C^{0}-$map $\gamma :I\rightarrow M$ from an interval $I\subset \mathbb{R}$
into $M$ with $0\in I$ and $\gamma (0)=m$.

Two curves $\gamma _{1}$ and $\gamma _{2}$ passing though a point $m\in U$
are \emph{tangent at} $m$ with respect to the chart $( U,\phi ) $ if $(\phi
\circ \gamma _{1})^{\prime }(0)=(\phi \circ \gamma _{2})^{\prime }(0)$.
Thus, two curves are tangent if they have identical tangent vectors (same
direction and speed) in a local chart on a manifold.

For a smooth manifold $M$ and a point $m\in M,$ the \emph{tangent space} $%
T_{m}M$ to $M$ at $m$ is the \emph{set of equivalence classes} of curves at $%
m$:
\begin{equation*}
T_{m}M=\{ [ \gamma ] _{m}:\gamma \text{ is a curve at a point }m\in M\} .
\end{equation*}

A $C^{k}-$map $\varphi :M\ni m\mapsto \varphi (m)\in N$ between two
manifolds $M$ and $N$ induces a linear map $T_{m}\varphi :T_{m}M\rightarrow
T_{\varphi (m)}N$ for each point $m\in M$, called a \textit{tangent map}, if
we have:\newline
$\quad$\newline

\begin{equation*}
\bfig \putsquare<1`1`1`1;2000`1000>(0,0)[\,\odot`\,\,\odot`\,\odot`\,\,%
\odot; T(\varphi)`\pi_M`\pi_N`\varphi]
\put(-500,750){\framebox(1000,600){}}\put(1500,750){\framebox(1030,600){}}
\put(0,0){\oval(1000,400)}\put(2000,0){\oval(1000,400)} \put(-100,-100){$m$%
}\put(2050,-100){$\varphi(m)$} \put(-400,0){$M$}\put(2300,0){$N$}
\multiput(-500,750)(100,0){10}{\line(0,1){600}}
\multiput(1500,750)(100,0){10}{\line(0,1){600}} \put(270,1400){$TM$%
}\put(1500,1400){$T(N)$} \put(-150,1400){$T_m(M)$}\put(1850,1400){$%
T_{\varphi(m)}(N)$} \linethickness{0.3mm}\put(0,750){\line(0,1){600}}
\put(2000,750){\line(0,1){600}} \efig
\end{equation*}
\bigbreak\bigbreak \noindent i.e., the following diagram commutes:
\begin{equation*}
\square <1`1`1`1;1000`500>[T_{m}M`T_{\varphi (m)}N`M\ni m`\varphi (m)\in
N;T_{m}\varphi `\pi _{M}`\pi _{N}`\varphi ]
\end{equation*}
with the \textit{natural projection} $\pi _{M}:TM\rightarrow M,$ given by $%
\pi_{M}( T_{m}M) =m,$ that takes a tangent vector $v$ to the point $m\in M$
at which the vector $v$ is attached i.e., $v\in T_{m}M$.

For an $n$D smooth manifold $M$, its $n$D \textit{tangent bundle} $TM$ is
the disjoint union of all its tangent spaces $T_{m}M$ at all points $m\in M$%
, $TM=\dbigsqcup\limits_{m\in M}T_{m}M$.

To define the smooth structure on $TM$, we need to specify how to construct
local coordinates on $TM$. To do this, let $( x^{1}(m),...,x^{n}(m)) $ be
local coordinates of a point $m$ on $M$ and let $( v^{1}(m),...,v^{n}(m)) $
be components of a tangent vector in this coordinate system. Then the $2n$
numbers $( x^{1}(m),...,x^{n}(m),\,v^{1}(m),...,v^{n}(m)) $ give a \emph{%
local coordinate system} on $TM$.

$TM=\dbigsqcup\limits_{m\in M}T_{m}M$ defines a family of vector spaces
parameterized by $M$. The inverse image $\pi _{M}^{-1}(m)$ of a point $m\in
M $ under the natural projection $\pi _{M}$ is the tangent space $T_{m}M$.
This space is called the \emph{fibre} of the tangent bundle over the point $%
m\in M$.

A $C^{k}-$map $\varphi :M\rightarrow N$ between two manifolds $M$ and $N$
induces a linear \textit{tangent map} $T\varphi:TM\rightarrow TN$ between
their tangent bundles, i.e., the following diagram commutes:
\begin{equation*}
\square <1`1`1`1;900`500>[TM`TN`M`N;T\varphi`\pi _{M}`\pi _{N}`\varphi ]
\end{equation*}

All tangent bundles and their tangent maps form the category $\mathcal{TB}$.
The category $\mathcal{TB}$ is the natural framework for \textit{Lagrangian
dynamics}.

Now, we can formulate the \emph{global version of the} \emph{chain rule}. If
$\varphi :M\rightarrow N$ and $\psi :N\rightarrow P$ are two smooth maps,
then we have $T(\psi \circ \varphi )=T\psi \circ T\varphi $. In other words,
we have a functor $T:\mathcal{M\Rightarrow TB}$ from the category $\mathcal{M%
}$ of smooth manifolds to the category $\mathcal{TB}$ of their tangent
bundles:
\begin{equation*}
\Atriangle<1`1`1;500>[M`N`P;\varphi `(\psi \circ \varphi )\qquad \overset{T}{%
\Longrightarrow }`\psi ]\qquad \Atriangle<1`1`1;500>[TM`TN`TP;T\varphi
`T(\psi \circ \varphi )`T\psi ]
\end{equation*}

\subsection{Cotangent Bundle and Hamiltonian Dynamics}

The cotangent bundle of a smooth $n-$manifold is the place is where 1--forms
live, and is itself a smooth $2n-$manifold. Covector--fields (1--forms) are
cross-sections of the cotangent bundle. The \textit{Hamiltonian} is a
natural energy function on the cotangent bundle (see \cite%
{GaneshSprBig,GaneshADG}).

A \emph{dual} notion to the tangent space $T_{m}M$ to a smooth manifold $M$
at a point $m$ is its \textit{cotangent space} $T_{m}^{\ast }M$ at the same
point $m$. Similarly to the tangent bundle, for a smooth manifold $M$ of
dimension $n$, its \textit{cotangent bundle} $T^{\ast }M$ is the disjoint
union of all its cotangent spaces $T_{m}^{\ast }M$ at all points $m\in M$,
i.e., $T^{\ast }M=\dbigsqcup\limits_{m\in M}T_{m}^{\ast }M$. Therefore, the
cotangent bundle of an $n-$manifold $M$ is the vector bundle $T^{\ast
}M=(TM)^{\ast }$, the (real) dual of the tangent bundle $TM$.

If $M$ is an $n-$manifold, then $T^{\ast }M$ is a $2n-$manifold. To define
the smooth structure on $T^{\ast }M$, we need to specify how to construct
local coordinates on $T^{\ast }M$. To do this, let $( x^{1}(m),...,x^{n}(m))
$ be local coordinates of a point $m$ on $M$ and let $(
p_{1}(m),...,p_{n}(m)) $ be components of a covector in this coordinate
system. Then the $2n$ numbers $(
x^{1}(m),...,x^{n}(m),\,p_{1}(m),...,p_{n}(m)) $ give a local coordinate
system on $T^{\ast }M$. This is the basic idea one uses to prove that indeed
$T^{\ast }M$ is a $2n-$manifold.

$T^{\ast }M=\dbigsqcup\limits_{m\in M}T_{m}^{\ast }M$ defines a family of
vector spaces parameterized by $M$, with the \textit{conatural projection}, $%
\pi _{M}^{\ast }:T^{\ast }M\rightarrow M,$ given by $\pi _{M}^{\ast }(
T_{m}^{\ast }M) =m,$ that takes a covector $p$ to the point $m\in M$ at
which the covector $p$ is attached i.e., $p\in T_{m}^{\ast }M$. The inverse
image $\pi _{M}^{-1}(m)$ of a point $m\in M$ under the conatural projection $%
\pi _{M}^{\ast }$ is the cotangent space $T_{m}^{\ast }M$. This space is
called the \emph{fibre} of the cotangent bundle over the point $m\in M$.

In a similar way, a $C^{k}-$map $\varphi :M\rightarrow N$ between two
manifolds $M$ and $N$ induces a linear \emph{cotangent map} $T^{\ast
}\varphi :T^{\ast }M\rightarrow T^{\ast }N$ between their cotangent bundles,
i.e., the following diagram commutes:
\begin{equation*}
\square <1`1`1`1;900`500>[T^{\ast }M`T^{\ast }N`M`N;T^{\ast }\varphi `\pi
_{M}^{\ast }`\pi _{N}^{\ast }`\varphi ]
\end{equation*}

All cotangent bundles and their cotangent maps form the category $\mathcal{%
T^{\ast }B}$. The category $\mathcal{T^{\ast }B}$ is the natural stage for
\textit{Hamiltonian dynamics}.

Now, we can formulate the \emph{dual version of the global chain rule}. If $%
\varphi :M\rightarrow N$ and $\psi :N\rightarrow P$ are two smooth maps,
then we have $T^{\ast }(\psi \circ \varphi )=T^{\ast }\psi \circ T^{\ast
}\varphi $. In other words, we have a cofunctor $T^{\ast }:\mathcal{%
M\Rightarrow T^{\ast }B}$ from the category $\mathcal{M}$ of smooth
manifolds to the category $\mathcal{T^{\ast }B}$ of their cotangent bundles:
\begin{equation*}
\Atriangle<1`1`1;500>[M`N`P;\varphi `(\psi \circ \varphi )\qquad \overset{%
T^{\ast }}{\Longrightarrow }`\psi ]\qquad \Atriangle<-1`-1`-1;500>[T^{\ast
}M`T^{\ast }N`T^{\ast }P;T^{\ast }\varphi `T^{\ast }(\psi \circ \varphi
)`T^{\ast }\psi ]
\end{equation*}

\section{Lie Groups}

In this section we present the basics of \textit{classical theory of Lie
groups} and their Lie algebras, as developed mainly by Sophus Lie, Elie
Cartan, Felix Klein, Wilhelm Killing and Hermann Weyl. For more
comprehensive treatment see e.g., \cite%
{Chevalley,Helgason,Gilmore,Fulton,Bourbaki}.

In the middle of the 19th Century S. Lie made a far reaching discovery that
techniques designed to solve particular unrelated types of ODEs, such as
separable, homogeneous and exact equations, were in fact all special cases
of a general form of integration procedure based on the invariance of the
differential equation under a continuous group of symmetries. Roughly
speaking a symmetry group of a system of differential equations is a group
that transforms solutions of the system to other solutions. Once the
symmetry group has been identified a number of techniques to solve and
classify these differential equations becomes possible. In the classical
framework of Lie, these groups were local groups and arose locally as groups
of transformations on some Euclidean space. The passage from the local Lie
group to the present day definition using manifolds was accomplished by E.
Cartan at the end of the 19th Century, whose work is a striking synthesis of
Lie theory, classical geometry, differential geometry and topology.

These continuous groups, which originally appeared as symmetry groups of
differential equations, have over the years had a profound impact on diverse
areas such as algebraic topology, differential geometry, numerical analysis,
control theory, classical mechanics, quantum mechanics etc. They are now
universally known as Lie groups.

A Lie group is smooth manifold which also carries a group structure whose
product and inversion operations are smooth as maps of manifolds. These
objects arise naturally in describing physical symmetries.\footnote{%
Here are a few examples of Lie groups and their relations to other areas of
mathematics and physics:
\par
\begin{enumerate}
\item Euclidean space $\mathbb{R}^{n}$ is an Abelian Lie group (with
ordinary vector addition as the group operation).
\par
\item The group $GL_{n}(\mathbb{R})$ of invertible matrices (under matrix
multiplication) is a Lie group of dimension $n^{2}$. It has a subgroup $%
SL_{n}(\mathbb{R})$ of matrices of determinant 1 which is also a Lie group.
\par
\item The group $O_{n}(\mathbb{R})$ generated by all rotations and
reflections of an $n$D vector space is a Lie group called the \textit{%
orthogonal group}. It has a subgroup of elements of determinant 1, called
the special orthogonal group $SO(n)$, which is the \textit{group of rotations%
} in $\mathbb{R}^{n}$.
\par
\item Spin groups are double covers of the special orthogonal groups (used
e.g., for studying fermions in quantum field theory).
\par
\item The group $Sp_{2n}(\mathbb{R})$ of all matrices preserving a
symplectic form is a Lie group called the \textit{symplectic group}.
\par
\item The Lorentz group and the Poincar\'{e} group of isometries of
space--time are Lie groups of dimensions 6 and 10 that are used in special
relativity.
\par
\item The Heisenberg group is a Lie group of dimension 3, used in quantum
mechanics.
\par
\item The unitary group $U(n)$ is a compact group of dimension $n^{2}$
consisting of unitary matrices. It has a subgroup of elements of determinant
1, called the special unitary group $SU(n)$.
\par
\item The group $U(1)\times SU(2)\times SU(3)$ is a Lie group of dimension $%
1+3+8=12$ that is the \textit{gauge group} of the \textit{Standard Model} of
elementary particles, whose dimension corresponds to: 1 photon + 3 vector
bosons + 8 gluons.
\end{enumerate}
}

A Lie group is a group whose elements can be continuously parametrized by
real numbers, such as the \textit{rotation group}\ $SO(3)$, which can be
parametrized by the \textit{Euler angles}. More formally, a Lie group is an
analytic real or complex manifold that is also a group, such that the group
operations multiplication and inversion are analytic maps. Lie groups are
important in mathematical analysis, physics and geometry because they serve
to describe the symmetry of analytical structures. They were introduced by
\emph{Sophus} \textit{Lie} in 1870 in order to study symmetries of
differential equations.

While the Euclidean space $\mathbb{R}^{n}$ is a \textit{real Lie group}
(with ordinary vector addition as the group operation), more typical
examples are given by matrix Lie groups, i.e., groups of invertible matrices
(under matrix multiplication). For instance, the group $SO(3)$ of all
rotations in $\mathbb{R}^{3}$ is a matrix Lie group.

One classifies Lie groups regarding their \textit{algebraic properties}%
\footnote{%
If $G$ and $H$ are Lie groups (both real or both complex), then a \textit{%
Lie--group--homomorphism} $f:G\rightarrow H$ is a \textit{group homomorphism}
which is also an \textit{analytic map} (one can show that it is equivalent
to require only that f be continuous). The composition of two such
homomorphisms is again a homomorphism, and the class of all (real or
complex) Lie groups, together with these morphisms, forms a \textit{category}%
. The two Lie groups are called \textit{isomorphic} iff there exists a
bijective homomorphism between them whose inverse is also a homomorphism.
Isomorphic Lie groups do not need to be distinguished for all practical
purposes; they only differ in the notation of their elements.} (simple,
semisimple, solvable, nilpotent, Abelian), their \textit{connectedness}
(connected or simply connected) and their \textit{compactness}.\footnote{%
An $n-$torus ~$T^n=S^1 \times S^1\times\cdots\times S^1$~ (as defined above)
is an example of a \textit{compact Abelian Lie group}. This follows from the
fact that the unit circle $S^1$ is a compact Abelian Lie group (when
identified with the unit complex numbers with multiplication). Group
multiplication on $T^n$ is then defined by coordinate--wise multiplication.
\par
Toroidal groups play an important part in the theory of compact Lie groups.
This is due in part to the fact that in any compact Lie group one can always
find a maximal torus; that is, a closed subgroup which is a torus of the
largest possible dimension.}

To every Lie group, we can associate a \textit{Lie algebra} which completely
captures the local structure of the group (at least if the Lie group is
connected).\footnote{%
Conventionally, one can regard any field $X$ of tangent vectors on a Lie
group as a partial differential operator, denoting by $Xf$ the \textit{Lie
derivative} (the \textit{directional derivative}) of the scalar field $f$ in
the direction of $X$. Then a vector--field on a Lie group $G$ is said to be
left--invariant if it commutes with left translation, which means the
following. Define $L_{g}[f](x)=f(gx)$ for any analytic function $%
f:G\rightarrow \mathbb{R}$ and all $g,x\in G$. Then the vector--field $X$ is
left--invariant iff $XL_{g}=L_{g}X$ for all $g\in G$. Similarly, instead of $%
\mathbb{R}$, we can use $\mathbb{C}$. The set of all vector--fields on an
analytic manifold is a \textit{Lie algebra} over $\mathbb{R}$ (or $\mathbb{C}
$).
\par
On a Lie group $G$, the left--invariant vector--fields form a subalgebra,
the Lie algebra $\mathfrak{g}$\ associated with $G$. This Lie algebra is
finite--dimensional (it has the same dimension as the manifold $G$) which
makes it susceptible to classification attempts. By classifying $\mathfrak{g}
$, one can also get a handle on the group $G$. The representation theory of
simple Lie groups is the best and most important example.
\par
Every element $v$ of the tangent space $T_{e}$ at the identity element $e$
of $G$ determines a unique left--invariant vector--field whose value at the
element $g$ of $G$ is denoted by $gv$; the vector space underlying the Lie
algebra $\mathfrak{g}$\ may therefore be identified with $T_{e}$.
\par
Every vector--field $v$ in the Lie algebra $\mathfrak{g}$ determines a
function $c:\mathbb{R}\rightarrow G$ whose derivative everywhere is given by
the corresponding left--invariant vector--field: ~$c^{\prime}(t)=TL_{c(t)}v$%
~ and which has the property: ~$c(s+t)=c(s)c(t),\qquad (\text{for all $s$
and $t$})$ (the operation on the r.h.s. is the group multiplication in $G$).
The formal similarity of this formula with the one valid for the elementary
exponential function justifies the definition: ~$\mathrm{exp}(v)=c(1).$~
This is called the \textit{exponential map}, and it maps the Lie algebra $%
\mathfrak{g}$ into the Lie group $G$. It provides a \textit{diffeomorphism}
between a neighborhood of $0$ in $\mathfrak{g}$\ and a neighborhood of $e$
in $G$. This exponential map is a generalization of the exponential function
for real numbers (since $\mathbb{R}$ is the Lie algebra of the Lie group of
positive real numbers with multiplication), for complex numbers (since $%
\mathbb{C}$ is the Lie algebra of the Lie group of non--zero complex numbers
with multiplication) and for matrices (since $M(n,\mathbb{R})$ with the
regular commutator is the Lie algebra of the Lie group $GL(n,\mathbb{R})$ of
all invertible matrices). As the exponential map is surjective on some
neighborhood $N$ of $e$, it is common to call elements of the Lie algebra
\textit{infinitesimal generators} of the group $G$.
\par
The exponential map and the Lie algebra determine the local group structure
of every connected Lie group, because of the \textit{%
Baker--Campbell--Hausdorff formula}: there exists a neighborhood $U$ of the
zero element of the Lie algebra $\mathfrak{g}$, such that for $u,v\in U$ we
have
\begin{equation*}
\mathrm{exp}(u)\mathrm{exp}(v)=\mathrm{exp}%
(u+v+1/2[u,v]+1/12[[u,v],v]-1/12[[u,v],u]-...),
\end{equation*}
where the omitted terms are known and involve \textit{Lie bracket}\emph{s}
of four or more elements. In case $u$ and $v$ commute, this formula reduces
to the familiar \textit{exponential law}:
\begin{equation*}
\mathrm{exp}(u)\mathrm{exp}(v)=\mathrm{exp}(u+v).
\end{equation*}%
\par
Every homomorphism $f:G\rightarrow H$ of Lie groups induces a homomorphism
between the corresponding Lie algebras $\mathfrak{g}$ and $\mathfrak{h}$.
The association $G\Longrightarrow \mathfrak{g}$ is called the \textit{Lie
Functor}.}

\subsection{Definition of a Lie Group}

A \textit{Lie group} is a smooth (Banach) manifold $M$ that has at the same
time a group $G-$structure consistent with its manifold $M-$structure in the
sense that \textit{group multiplication} ~$\mu :G\times G\rightarrow G,~~
(g,h)\mapsto gh$ and the \textit{group inversion} ~$\nu :G\rightarrow G,~~
g\mapsto g^{-1}$ are $C^{k}-$maps. A point $e\in G$ is called the \textit{%
group identity element} (see e.g., \cite{Chevalley,Helgason,Arnold,Abraham}).

For example, any $n$D Banach vector space $V$ is an Abelian Lie group with
group operations $\mu :V\times V\rightarrow V$, $\mu (x,y)=x+y$, and $\nu
:V\rightarrow V$, $\nu (x)=-x$. The identity is just the zero vector. We
call such a Lie group a \emph{vector group}.

Let $G$ and $H$ be two Lie groups. A map $G\rightarrow H$ is said to be a
\emph{morphism} of Lie groups (or their \textit{smooth homomorphism}) if it
is their homomorphism as abstract groups and their smooth map as manifolds.

Similarly, a group $G$ which is at the same time a topological space is said
to be a \textit{topological group} if both maps ($\mu ,\nu $) are
continuous, i.e., $C^{0}-$maps for it. The homomorphism $G\rightarrow H$ of
topological groups is said to be continuous if it is a continuous map.

A topological group (as well as a smooth manifold) is not necessarily
Hausdorff. A topological group $G$ is Hausdorff iff its identity is closed.
As a corollary we have that every Lie group is a Hausdorff topological group.

For every $g$ in a Lie group $G$, the two maps, \
\begin{eqnarray*}
L_{g} &:&G\rightarrow G,\qquad h\mapsto gh, \\
R_{h} &:&G\rightarrow G,\qquad g\mapsto gh,
\end{eqnarray*}
are called \emph{left} and \textit{right translation} maps. Since $%
L_{g}\circ L_{h}=L_{gh}$, and $R_{g}\circ R_{h}=R_{gh}$, it follows that $%
\left( L_{g}\right) ^{-1}=L_{g^{-1}}$ and $\left( R_{g}\right)
^{-1}=R_{g^{-1}}$, so both $L_{g}$ and $R_{g}$ are diffeomorphisms. Moreover
$L_{g}\circ R_{h}=R_{h}\circ L_{g}$, i.e., left and right translation
commute.

A vector--field $X$ on $G$ is called \emph{left--invariant vector--field} if
for every $g\in G$, $L_{g}^{\ast }X=X$, that is, if $(T_{h}L_{g})X(h)=X(gh)$
for all $h\in G$, i.e., the following diagram commutes:
\begin{equation*}
\square <1`-1`-1`1;900`500>[TG`TG`G`G;TL_{g}`X`X`L_{g}]
\end{equation*}

A \textit{Riemannian metric} on a Lie group $G$ is called left-invariant if
it is preserved by all left translations $L_{g}$, i.e., if the derivative of
left translation carries every vector to a vector of the same length.
Similarly, a vector field $X$ on $G$ is called left--invariant if (for every
$g\in G$) $L_{g}^{\ast }X=X$.

\subsection{Lie Algebra}

An \emph{algebra} $A$ is a vector space with a product. The product must
have the property that
\begin{equation*}
a(uv)=(au)v=u(av),
\end{equation*}%
for every $a\in \mathbb{R}$ and $u,v\in A$. A map $\phi :A\rightarrow
A^{\prime }$ between algebras is called an \textit{algebra homomorphism} if $%
\phi (u\cdot v)=\phi (u)\cdot \phi (v)$. A vector subspace $\mathfrak{I}$ of
an algebra $A$ is called a \textit{left ideal} (resp. \emph{right ideal}) if
it is closed under algebra multiplication and if $u\in A$ and $i\in
\mathfrak{I}$ implies that $ui\in \mathfrak{I}$ (resp. $iu\in \mathfrak{I}$%
). A subspace $\mathfrak{I}$ is said to be a \emph{two--sided ideal} if it
is both a left and right ideal. An ideal may not be an algebra itself, but
the quotient of an algebra by a two--sided ideal inherits an algebra
structure from $A$.

A \emph{Lie algebra} is an algebra $A$ where the multiplication, i.e., the
\emph{Lie bracket} $(u,v)\mapsto \lbrack u,v]$, has the following properties:

LA 1. $[u,u]=0$ for every $u\in A$, and

LA 2. $[u,[v,w]]+[w,[u,v]]+[v,w,u]]=0$ for all $u,v,w\in A$.

The condition LA 2 is usually called \textit{Jacobi identity}. A subspace $%
E\subset A$ of a Lie algebra is called a \textit{Lie subalgebra} if $%
[u,v]\in E$ for every $u,v\in E$. A map $\phi:A\rightarrow A^{\prime }$
between Lie algebras is called a \textit{Lie algebra homomorphism} if $\phi
([u,v])=[\phi (u),\phi (v)]$ for each $u,v\in A$.

All Lie algebras (over a given field $\mathbb{K}$) and all smooth
homomorphisms between them form the category $\mathcal{LAL}$, which is
itself a complete subcategory of the category $\mathcal{AL}$ of all algebras
and their homomorphisms.

Let $\mathcal{X}_{L}(G)$ denote the set of left--invariant vector--fields on
$G $; it is a Lie subalgebra of $\mathcal{X}(G)$, the set of all
vector--fields on $G$, since $L_{g}^{\ast }[X,Y]=[L_{g}^{\ast }X,L_{g}^{\ast
}Y]=[X,Y]$, so the Lie bracket $[X,Y]\in \mathcal{X}_{L}(G)$.

Let $e$ be the identity element of $G$. Then for each $\xi $ on the tangent
space $T_{e}G$ we define a vector--field $X_{\xi }$ on $G$ by $X_{\xi
}(g)=T_{e}L_{g}(\xi )$. $\mathcal{X}_{L}(G)$ and $T_{e}G$ are isomorphic as
vector spaces. Define the Lie bracket on $T_{e}G$ by $\lbrack \xi ,\eta ]=%
\left[ X_{\xi },X_{\eta }\right] (e)$ for all $\xi ,\eta \in T_{e}G$. This
makes $T_{e}G$ into a Lie algebra. Also, by construction, we have $\left[
X_{\xi },X_{\eta }\right] =X_{[\xi ,\eta ]}$; this defines a bracket in $%
T_{e}G$ via \emph{left extension}. The vector space $T_{e}G$ with the above
algebra structure is called the Lie algebra of the Lie group $G$ and is
denoted $\mathfrak{g}$.

For example, let $V$ be a $n$D vector space. Then $T_{e}V\simeq V$ and the
left--invariant vector--field defined by $\xi \in T_{e}V$ is the constant
vector--field $X_{\xi }(\eta )=\xi $, for all $\eta \in V$. The Lie algebra
of $V$ is $V$ itself.

Since any two elements of an Abelian Lie group $G$ commute, it follows that
all adjoint operators $Ad_{g}$, $g\in G$, equal the identity. Therefore, the
Lie algebra $g$ is Abelian; that is, $[\xi ,\eta ]=0$ for all $\xi ,\eta \in
\mathfrak{g}$.

For example, $G=SO(3)$ is the group of rotations of 3D Euclidean space, i.e.
the configuration space of a rigid body fixed at a point. A motion of the
body is then described by a curve $g=g(t)$ in the group $SO(3)$. Its Lie
algebra $\mathfrak{g}=\mathfrak{so}(3)$\ is the 3D vector space of angular
velocities of all possible rotations. The commutator in this algebra is the
usual vector (cross) product (see, e.g. \cite{Arnold,Abraham,GaneshADG}).

A rotation velocity $\dot{g}$ of the rigid body (fixed at a point) is a
tangent vector to the Lie group $G=SO(3)$ at the point $g\in G$. To get the
angular velocity, we must carry this vector to the tangent space $TG_{e}$ of
the group at the identity, i.e. to its Lie algebra $\mathfrak{g}=\mathfrak{so%
}(3)$. This can be done in two ways: by left and right translation, $L_{g}$
and $R_{g}$. As a result, we obtain two different vector fields in the Lie
algebra $\mathfrak{so}(3):$

\begin{equation*}
\omega _{c}=L_{g^{-1}\ast }\dot{g}\in \mathfrak{so}(3)\qquad \text{and\qquad
}\omega _{x}=R_{g^{-1}\ast }\dot{g}\in \mathfrak{so}(3),
\end{equation*}%
which are called the `angular velocity in the body' and the `angular
velocity in space,' respectively.

The dual space $\mathfrak{g}^{\ast }$ to the Lie algebra $\mathfrak{g}=%
\mathfrak{so}(3)$ is the space of angular momenta $\mathbf{\pi }$. The
kinetic energy $T$ of a body is determined by the vector field of angular
velocity in the body and does not depend on the position of the body in
space. Therefore, kinetic energy gives a left-invariant Riemannian metric on
the rotation group $G=SO(3)$.

\subsection{One-Parameter Subgroup}

Let $X_{\xi }$ be a left--invariant vector--field on $G$ corresponding to $%
\xi $ in $\mathfrak{g}$. Then there is a unique integral curve $\gamma _{\xi
}:\mathbb{R}\rightarrow G$ of $X_{\xi }$ starting at $e$, i.e., (see, e.g.
\cite{GaneshSprBig,GaneshADG})
\begin{equation*}
\dot{\gamma}_{\xi }(t)=X_{\xi }\left( \gamma _{\xi }(t)\right) ,\qquad
\gamma _{\xi }(0)=e
\end{equation*}
$\gamma _{\xi }(t)$ is a smooth \emph{one--parameter subgroup} of $G$, i.e.,
$\gamma _{\xi }(t+s)=\gamma _{\xi }(t)\cdot \gamma _{\xi }(s)$, since, as
functions of $t$ both sides equal $\gamma _{\xi }(s)$ at $t=0$ and both
satisfy differential equation $\dot{\gamma}(t)=X_{\xi }\left( \gamma _{\xi
}(t)\right) $ by left invariance of $X_{\xi }$, so they are equal. Left
invariance can be also used to show that $\gamma _{\xi }(t)$ is defined for
all $t\in \mathbb{R}$. Moreover, if $\phi :\mathbb{R}\rightarrow G$ is a
one--parameter subgroup of $G$, i.e., a \emph{smooth homomorphism} of the
additive group $\mathbb{R}$ into $G$, then $\phi =\gamma _{\xi }$ with $\xi =%
\dot{\phi}(0)$, since taking derivative at $s=0$ in the relation
\begin{equation*}
\phi (t+s)=\phi (t)\cdot \phi (s)\qquad \text{gives\qquad }\dot{\phi}(t)=X_{%
\dot{\phi}(0)}\left( \phi (t)\right) ,
\end{equation*}%
so $\phi =\gamma _{\xi }$ since both equal $e$ at $t=0$. Therefore, all
one--parameter subgroups of $G$ are of the form $\gamma _{\xi }(t)$ for some
$\xi \in \mathfrak{g}$.

\subsection{Exponential Map}

The map $\exp :\mathfrak{g}\rightarrow G$, given by (see, e.g. \cite%
{Abraham,GaneshSprBig,GaneshADG}):
\begin{equation*}
\exp (\xi )=\gamma _{\xi }(1),\qquad \exp (0) = e
\end{equation*}
is called the \textit{exponential map} of the Lie algebra $\mathfrak{g}$ of $%
G$ into $G$. $\exp $ is a $C^{k}-$-map, similar to the projection $\pi$ of
tangent and cotangent bundles; $\exp $ is locally a diffeomorphism from a
neighborhood of zero in $\mathfrak{g}$ onto a neighborhood of $e$ in $G$; if
$f:G\rightarrow H$ is a smooth homomorphism of Lie groups, then
\begin{equation*}
f\circ \exp _{G}=\exp _{H}\circ T_{e}f\,.
\end{equation*}

Also, in this case
\begin{equation*}
\exp (s\xi )=\gamma _{\xi }(s).
\end{equation*}%
Indeed, for fixed $s\in \mathbb{R}$, the curve $t\mapsto \gamma _{\xi }(ts)$%
, which at $t=0$ passes through $e$, satisfies the differential equation
\begin{equation*}
\frac{d}{dt}\gamma _{\xi }(ts)=sX_{\xi }\left( \gamma _{\xi }(ts)\right)
=X_{s\xi }\left( \gamma _{\xi }(ts)\right) .
\end{equation*}
Since $\gamma _{s\xi }(t)$ satisfies the same differential equation and
passes through $e$ at $t=0$, it follows that $\gamma _{s\xi }(t)=\gamma
_{\xi }(st)$. Putting $t=1$ induces $\exp (s\xi )=\gamma _{\xi }(s)$.

Hence $\exp $ maps the line $s\xi $ in $\mathfrak{g}$ onto the
one--parameter subgroup $\gamma _{\xi }(s)$ of $G$, which is tangent to $\xi
$ at $e$. It follows from left invariance that the flow $F_{t}^{\xi }$ of $X$
satisfies $F_{t}^{\xi }(g)=g\exp (s\xi )$.

Globally, the exponential map $\exp$ is a natural operation, i.e., for any
morphism $\varphi:G\rightarrow H$ of Lie groups $G$ and $H$ and a Lie
functor $\mathcal{F}$, the following diagram commutes:
\begin{equation*}
\square <1`1`1`1;900`500>[\mathcal{F}(G)`\mathcal{F}(H)`G`H;\mathcal{F}%
(\varphi) `\exp`\exp`\varphi ]
\end{equation*}

Let $G_{1}$ and $G_{2}$ be Lie groups with Lie algebras $\mathfrak{g}_{1}$
and $\mathfrak{g}_{2}.$\ $\ $Then $G_{1}\times G_{2}$ is a Lie group with
Lie algebra $\mathfrak{g}_{1}\times \mathfrak{g}_{2},$ and the exponential
map is given by:
\begin{equation*}
\exp :\mathfrak{g}_{1}\times \mathfrak{g}_{2}\rightarrow G_{1}\times
G_{2},\qquad (\xi _{1},\xi _{2})\mapsto \left( \exp _{1}(\xi _{1}),\exp
_{2}(\xi _{2})\right) .
\end{equation*}

For example, in case of a $n$D vector space, or infinite--dimensional Banach
space, the exponential map is the identity.

The unit circle in the complex plane $S^{1}=\{z\in \mathbb{C}:\left\vert
z\right\vert =1\}$ is an Abelian Lie group under multiplication. The tangent
space $T_{e}S^{1}$ is the imaginary axis, and we identify $\mathbb{R}$ with $%
T_{e}S^{1}$ by $t\mapsto 2\pi it$. With this identification, the exponential
map $\exp :\mathbb{R}\rightarrow S^{1}$ is given by $\exp (t)=\mathrm{e}%
^{2\pi it}$.

The $n$D torus $T^{n}=S^{1}\times $\textperiodcentered \textperiodcentered
\textperiodcentered $\times S^{1}$ ($n$ times) is an Abelian Lie group. The
exponential map $\exp :\mathbb{R}^{n}\rightarrow T^{n}$ is given by
\begin{equation*}
\exp (t_{1},...,t_{n})=(\mathrm{e}^{2\pi it_{1}},...,\mathrm{e}^{2\pi
it_{n}}).
\end{equation*}%
Since $S^{1}=\mathbb{R}/\mathbb{Z}$, it follows that $T^{n}=\mathbb{R}^{n}/%
\mathbb{Z}^{n}$, the projection $\mathbb{R}^{n}\rightarrow T^{n}$ being
given by the $\exp $ map.

\subsection{Adjoint Representation}

For every $g\in G$, the map (see, e.g. \cite%
{Arnold,Abraham,GaneshSprBig,GaneshADG}):
\begin{equation*}
Ad_{g}=T_{e}\left( R_{g^{-1}}\circ L_{g}\right) :\mathfrak{g}\rightarrow
\mathfrak{g}
\end{equation*}
is called the \textit{adjoint map}, or \textit{adjoint operator} associated
with $g$.

For each $\xi \in \mathfrak{g}$ and $g\in G$ we have
\begin{equation*}
\exp \left( Ad_{g}\xi \right) =g\left( \exp \xi \right) g^{-1}.
\end{equation*}

The relation between the adjoint map and the Lie bracket is the following:
For all $\xi ,\eta \in \mathfrak{g}$ we have
\begin{equation*}
\left. \frac{d}{dt}\right| _{t=0}Ad_{\exp (t\xi )}\eta =[\xi ,\eta ].
\end{equation*}

Left and right translations induce operators on the cotangent space $T^{\ast
}G_{g}$\ dual to $L_{g\ast }$ and $R_{g\ast },$ denoted by (for every $h\in
G $):%
\begin{equation*}
L_{g}^{\ast }:T^{\ast }G_{gh}\rightarrow T^{\ast }G_{h},\qquad R_{g}^{\ast
}:T^{\ast }G_{hg}\rightarrow T^{\ast }G_{h}.
\end{equation*}%
The transpose operators $Ad_{g}^{\ast }:\mathfrak{g}\rightarrow \mathfrak{g}$
satisfy the relations $Ad_{gh}^{\ast }=Ad_{h}^{\ast }Ad_{g}^{\ast }$ (for
every $g,h\in G$) and constitute the \textit{co-adjoint representation} of
the Lie group $G$. The co-adjoint representation plays an important role in
all questions related to (left) invariant metrics on the Lie group.
According to A. Kirillov, the orbit of any vector field $X$ in a Lie algebra
$\mathfrak{g}$ in a co-adjoint representation $Ad_{g}^{\ast }$ is itself a
symplectic manifold and therefore a phase space for a Hamiltonian mechanical
system.

A Lie subgroup $H$ of $G$ is a subgroup $H$ of $G$ which is also a
submanifold of $G$. Then $\mathfrak{h}$ is a Lie subalgebra of $\mathfrak{g}$
and moreover $\mathfrak{h}=\{\xi \in \mathfrak{g}|\exp (t\xi )\in H$, for
all $t\in \mathbb{R}\}.$

One can characterize \textit{Lebesgue measure} up to a multiplicative
constant on $\mathbb{R}^{n}$ by its invariance under translations.
Similarly, on a locally compact group there is a unique (up to a nonzero
multiplicative constant) left--invariant measure, called \textit{Haar measure%
}. For Lie groups the existence of such measures is especially simple: Let $%
G $ be a Lie group. Then there is a volume form $U{b5}$, unique up to
nonzero multiplicative constants, that is left--invariant. If $G$ is
compact, $U{b5}$ is right invariant as well.

\subsection{Actions of Lie Groups on Smooth Manifolds}

Let $M$ be a smooth manifold. An \textit{action of a Lie group} $G$ (with
the unit element $e$) on $M$ is a smooth map $\phi :G\times M\rightarrow M,$
such that for all $x\in M$ and $g,h\in G$, (i) $\phi (e,x)=x$ and (ii) $\phi
\left( g,\phi (h,x)\right) =\phi (gh,x).$ In other words, letting $\phi
_{g}:x\in M\mapsto \phi _{g}(x)=\phi (g,x)\in M$, we have (i') $\phi
_{e}=id_{M}$ and (ii') $\phi _{g}\circ \phi _{h}=\phi _{gh}$. $\phi _{g}$ is
a diffeomorphism, since $(\phi _{g})^{-1}=\phi _{g^{-1}}$. We say that the
map $g\in G\mapsto \phi _{g}\in Diff(M)$ is a homomorphism of $G$ into the
group of diffeomorphisms of $M$. In case that $M$ is a vector space and each
$\phi _{g}$ is a linear operator, the function of $G$ on $M$ is called a
representation of $G$ on $M$ (see, e.g. \cite%
{Arnold,Abraham,GaneshSprBig,GaneshADG}).

An action $\phi$ of $G$ on $M$ is said to be \textit{transitive group action}%
, if for every $x,y\in M$, there is $g\in G$ such that $\phi (g,x)=y$;
\textit{effective group action}, if $\phi _{g}=id_{M}$ implies $g=e$, that
is $g\mapsto \phi _{g}$ is 1--1; and \textit{free group action}, if for each
$x\in M $, $g\mapsto \phi _{g}(x)$ is 1--1.

For example,

\begin{enumerate}
\item $G=\mathbb{R}$ acts on $M=\mathbb{R}$ by translations; explicitly,
\begin{equation*}
\phi :G\times M\rightarrow M,\qquad \phi (s,x)=x+s.
\end{equation*}
Then for $x\in \mathbb{R}$, $O_{x}=\mathbb{R}$. Hence $M/G$ is a single
point, and the action is transitive and free.

\item A complete flow $\phi _{t}$ of a vector--field $X$ on $M$ gives an
action of $\mathbb{R}$ on $M$, namely
\begin{equation*}
(t,x)\in \mathbb{R}\times M\mapsto \phi _{t}(x)\in M.
\end{equation*}

\item Left translation $L_{g}:G\rightarrow G$ defines an effective action of
$G$ on itself. It is also transitive.

\item The coadjoint action of $G$ on $\mathfrak{g}^{\ast }$ is given by
\begin{equation*}
Ad^{\ast }:(g,\alpha )\in G\times \mathfrak{g}^{\ast }\mapsto
Ad_{g^{-1}}^{\ast }(\alpha )=\left( T_{e}(R_{g^{-1}}\circ L_{g})\right)
^{\ast }\alpha \in \mathfrak{g}^{\ast }.
\end{equation*}
\end{enumerate}

Let $\phi $ be an action of $G$ on $M$. For $x\in M$ the \emph{orbit} of $x $
is defined by
\begin{equation*}
O_{x}=\{\phi _{g}(x)|g\in G\}\subset M
\end{equation*}
and the \textit{isotropy group} of $\phi $ at $x$ is given by
\begin{equation*}
G_{x}=\{g\in G|\phi (g,x)=x\}\subset G.
\end{equation*}

An action $\phi $ of $G$ on a manifold $M$ defines an equivalence relation
on $M$ by the relation belonging to the same orbit; explicitly, for $x,y\in
M $, we write $x\sim y$ if there exists a $g\in G$ such that $\phi (g,x)=y$,
that is, if $y\in O_{x}.$ The set of all orbits $M/G$ is called the \textit{%
group orbit space} (see, e.g. \cite{Arnold,Abraham,GaneshSprBig,GaneshADG}).

For example, let $M=\mathbb{R}^{2}\backslash \{0\}$, $G=SO(2)$, the group of
rotations in plane, and the action of $G$ on $M$ given by
\begin{equation*}
\left( \left[
\begin{array}{cc}
\cos \theta & -\sin \theta \\
\sin \theta & \cos \theta%
\end{array}
\right] ,(x,y)\right) \longmapsto (x\cos \theta -y\sin \theta ,\,x\sin
\theta +y\cos \theta ).
\end{equation*}
The action is always free and effective, and the orbits are concentric
circles, thus the orbit space is $M/G\simeq \mathbb{R}_{+}^{\ast }.$

A crucial concept in mechanics is the \textit{infinitesimal description of
an action}. Let $\phi :G\times M\rightarrow M$ be an action of a Lie group $%
G $ on a smooth manifold $M$. \ For each $\xi \in \mathfrak{g},$%
\begin{equation*}
\phi _{\xi }:\mathbb{R}\times M\rightarrow M,\qquad \phi _{\xi }(t,x)=\phi
\left( \exp (t\xi ),x\right)
\end{equation*}
is an $\mathbb{R}-$-action on $M$. Therefore, $\phi _{\exp (t\xi
)}:M\rightarrow M$ is a flow on $M$; the corresponding vector--field on $M$,
given by
\begin{equation*}
\xi _{M}(x)=\left. \frac{d}{dt}\right| _{t=0}\phi _{\exp (t\xi )}(x)
\end{equation*}
is called the infinitesimal generator of the action, corresponding to $\xi $
in $\mathfrak{g}.$

The tangent space at $x$ to an orbit $O_{x}$ is given by
\begin{equation*}
T_{x}O_{x}=\{\xi _{M}(x)|\xi \in \mathfrak{g}\}.
\end{equation*}

Let $\phi :G\times M\rightarrow M$ be a smooth $G--$action. For all $g\in G$%
, all $\xi ,\eta \in \mathfrak{g}$ and all $\alpha ,\beta \in \mathbb{R}$,
we have:

$\left( Ad_{g}\xi \right) _{M}=\phi _{g^{-1}}^{\ast }\xi _{M}$, $\left[ \xi
_{M},\eta _{M}\right] =-\left[ \xi ,\eta \right] _{M}$, and $(\alpha \xi
+\beta \eta )_{M}=\alpha \xi _{M}+\beta \eta _{M}$.

Let $M$ be a smooth manifold, $G$ a Lie group and $\phi :G\times
M\rightarrow M$ a $G-$action on $M$. We say that a smooth map $%
f:M\rightarrow M$ is with respect to this action if for all $g\in G$,
\begin{equation*}
f\circ \phi _{g}=\phi _{g}\circ f\text{.}
\end{equation*}
Let $f:M\rightarrow M$ be an equivariant smooth map. Then for any $\xi \in
\mathfrak{g}$ we have
\begin{equation*}
Tf\circ \xi _{M}=\xi _{M}\circ f.
\end{equation*}

\subsection{Basic Tables of Lie Groups and Their Lie Algebras}

One classifies Lie groups regarding their algebraic properties (simple,
semisimple, solvable, nilpotent, Abelian), their connectedness (connected or
simply connected) and their compactness (see Tables A.1--A.3). This is the
content of the \textit{Hilbert 5th problem}.\newpage

\noindent\textbf{Some real Lie groups and their Lie algebras:}\newline

\noindent{\footnotesize
\begin{tabular}{|p{36pt}|p{65pt}|p{65pt}|p{25pt}|p{65pt}|p{11pt}|}
\hline
\textbf{Lie group} & \textbf{Description} & \textbf{Remarks} & \textbf{Lie
\,\,\, algb.} & \textbf{Description} & \textbf{dim /$\mathbb{R}$} \\
\hline\hline
$\mathbb{R}^{n}$ & Euclidean space with addition & Abelian, simply
connected, not compact & $\mathbb{R}^{n}$ & the Lie bracket is zero & $n$ \\
\hline
$\mathbb{R}^{\mathrm{\times} }$ & nonzero real numbers with multiplication &
Abelian, not connected, not compact & $\mathbb{R}$ & the Lie bracket is zero
& 1 \\ \hline
$\mathbb{R}^{\mathrm{>} \mathrm{0}}$ & positive real numbers with
multiplication & Abelian, simply connected, not compact & $\mathbb{R}$ & the
Lie bracket is zero & 1 \\ \hline
$S^1 = \mathbb{R}/\mathbb{Z}$ & complex numbers of absolute value 1, with
multiplication & Abelian, connected, not simply connected, compact & $%
\mathbb{R} $ & the Lie bracket is zero & 1 \\ \hline
$\mathbb{H}^{\mathrm{\times} }$ & non--zero quaternions with multiplication
& simply connected, not compact & $\mathbb{H}$ & quaternions, with Lie
bracket the commutator & 4 \\ \hline
$S^{\mathrm{3}}$ & quaternions of absolute value 1, with multiplication; a $%
3-$sphere & simply connected, compact, simple and semi--simple, isomorphic
to $SU(2)$, $SO(3)$ and to $Spin(3)$ & $\mathbb{R}^{\mathrm{3}}$ & real $3-$%
vectors, with Lie bracket the cross product; isomorphic to $\mathfrak{su}(2)$
and to $\mathfrak{so}(3)$ & 3 \\ \hline
$GL(n,\mathbb{R})$ & general linear group: invertible $n-$by-$n$ real
matrices & not connected, not compact & M($n,\mathbb{R}$) & $n-$by-$n$
matrices, with Lie bracket the commutator & $n^{\mathrm{2}}$ \\ \hline
$GL^+(n,\mathbb{R})$ & $n-$by-$n$ real matrices with positive determinant &
simply connected, not compact & M($n,\mathbb{R}$) & $n-$by-$n$ matrices,
with Lie bracket the commutator & $n^{\mathrm{2}}$ \\ \hline\hline
\end{tabular}
}\newpage

\noindent\textbf{Classical real Lie groups and their Lie algebras:}\newline

\noindent{\footnotesize
\begin{tabular}{|p{36pt}|p{65pt}|p{65pt}|p{25pt}|p{65pt}|p{13pt}|}
\hline
\textbf{Lie group} & \textbf{Description} & \textbf{Remarks} & \textbf{Lie
\,\,\, algb.} & \textbf{Description} & \textbf{dim /$\mathbb{R}$} \\
\hline\hline
$SL(n,\mathbb{R})$ & special linear group: real matrices with determinant 1
& simply connected, not compact if $n>1$ & $\mathfrak{sl}(n,\mathbb{R})$ &
square matrices with trace 0, with Lie bracket the commutator & $n^2-1$ \\
\hline
$O(n,\mathbb{R})$ & orthogonal group: real orthogonal matrices & not
connected, compact & $\mathfrak{so}(n,\mathbb{R})$ & skew--symmetric square
real matrices, with Lie bracket the commutator; $\mathfrak{so}(3,\mathbb{R})$
is isomorphic to $\mathfrak{su}(2)$ and to $\mathbb{R }^{\mathrm{3}}$ with
the cross product & $n(n-1)/2$ \\ \hline
$SO(n,\mathbb{R})$ & special orthogonal group: real orthogonal matrices with
determinant 1 & connected, compact, for $n \ge 2$: not simply connected, for
$n=3$ and $n \ge 5$: simple and semisimple & $\mathfrak{so}(n,\mathbb{R})$ &
skew--symmetric square real matrices, with Lie bracket the commutator & $%
n(n-1)/2$ \\ \hline
$Spin(n)$ & spinor group & simply connected, compact, for $n=3$ and $n \ge 5$%
: simple and semisimple & $\mathfrak{so}(n,\mathbb{R})$ & skew--symmetric
square real matrices, with Lie bracket the commutator & $n(n-1)/2$ \\ \hline
$U(n)$ & unitary group: complex unitary $n-$by-$n$ matrices & isomorphic to $%
S^1$ for $n=1$, not simply connected, compact & $\mathfrak{u}(n)$ & square
complex matrices $A$ satisfying $A = -A^\ast$, with Lie bracket the
commutator & $n^{\mathrm{2}}$ \\ \hline
$SU(n)$ & special unitary group: complex unitary $n-$by-$n$ matrices with
determinant 1 & simply connected, compact, for $n \ge 2$: simple and
semisimple & $\mathfrak{su}(n)$ & square complex matrices $A$ with trace 0
satisfying $A = -A^\ast$, with Lie bracket the commutator & $n^2-1$ \\
\hline\hline
\end{tabular}
}\newpage

\noindent \textbf{Basic complex Lie groups and their Lie algebras:}\footnote{%
The dimensions given are dimensions over $\mathbb{C}$. Note that every
complex Lie group/algebra can also be viewed as a real Lie group/algebra of
twice the dimension.}\newline

\noindent{\footnotesize
\begin{tabular}{|p{36pt}|p{65pt}|p{65pt}|p{25pt}|p{65pt}|p{11pt}|}
\hline
\textbf{Lie group} & \textbf{Description} & \textbf{Remarks} & \textbf{Lie
\,\,\, algb.} & \textbf{Description} & \textbf{dim /$\mathbb{C}$} \\
\hline\hline
$\mathbb{C}^{n}$ & group operation is addition & Abelian, simply connected,
not compact & $\mathbb{C}^{n}$ & the Lie bracket is zero & $n$ \\ \hline
$\mathbb{C}^\times$ & nonzero complex numbers with multiplication & Abelian,
not simply connected, not compact & $\mathbb{C}$ & the Lie bracket is zero &
1 \\ \hline
$GL(n,\mathbb{C})$ & general linear group: invertible $n-$by-$n$ complex
matrices & simply connected, not compact, for $n=1$: isomorphic to $\mathbb{C%
}^{\mathrm{\times} }$ & $M(n,\mathbb{C})$ & $n-$by-$n$ matrices, with Lie
bracket the commutator & $n^{\mathrm{2}}$ \\ \hline
$SL(n,\mathbb{C})$ & special linear group: complex matrices with determinant
1 & simple, semisimple, simply connected, for $n \ge 2$: not compact & $%
\mathfrak{sl}(n, \mathbb{C})$ & square matrices with trace 0, with Lie
bracket the commutator & $n^{\mathrm{2}}-1$ \\ \hline
$O(n,\mathbb{C})$ & orthogonal group: complex orthogonal matrices & not
connected, for $n \ge 2$: not compact & $\mathfrak{so}(n,\mathbb{C})$ &
skew--symmetric square complex matrices, with Lie bracket the commutator & $%
n(n-1)/2$ \\ \hline
$SO(n,\mathbb{C})$ & special orthogonal group: complex orthogonal matrices
with determinant 1 & for $n \ge 2$: not compact, not simply connected, for $%
n=3$ and $n \ge 5$: simple and semisimple & $\mathfrak{so}(n,\mathbb{C})$ &
skew--symmetric square complex matrices, with Lie bracket the commutator & $%
n(n-1)/2$ \\ \hline\hline
\end{tabular}
}

\subsection{Representations of Lie groups}

The idea of a \textit{representation of a Lie group} plays an important role
in the study of continuous symmetry (see, e.g., \cite{Helgason}). A great
deal is known about such representations, a basic tool in their study being
the use of the corresponding 'infinitesimal' representations of Lie algebras.

Formally, a representation of a Lie group $G$ on a vector space $V$ (over a
field $K$) is a group homomorphism $G\rightarrow Aut(V)$ from $G$ to the
automorphism group of $V$. If a basis for the vector space $V$ is chosen,
the representation can be expressed as a homomorphism into $GL(n,K)$. This
is known as a \textit{matrix representation}.

On the Lie algebra level, there is a corresponding linear map from the Lie
algebra of $G$ to $End(V)$ preserving the Lie bracket $[\cdot,\cdot]$.

If the homomorphism is in fact an monomorphism, the representation is said
to be \emph{faithful}.

A unitary representation is defined in the same way, except that $G$ maps to
unitary matrices; the Lie algebra will then map to skew--Hermitian matrices.

Now, if $G$ is a semisimple group, its finite--dimensional representations
can be decomposed as direct sums of irreducible representations. The
irreducibles are indexed by highest weight; the allowable (\emph{dominant})
highest weights satisfy a suitable positivity condition. In particular,
there exists a set of \emph{fundamental weights}, indexed by the vertices of
the \textit{Dynkin diagram} of $G$ (see below), such that dominant weights
are simply non--negative integer linear combinations of the fundamental
weights.

If $G$ is a commutative compact Lie group, then its irreducible
representations are simply the continuous characters of $G$. A \textit{%
quotient representation} is a quotient module of the group ring.

\subsection{Root Systems and Dynkin Diagrams}

A \textit{root system} is a special configuration in Euclidean space that
has turned out to be fundamental in Lie theory as well as in its
applications. Also, the classification scheme for root systems, by \emph{%
Dynkin diagrams}, occurs in parts of mathematics with no overt connection to
Lie groups (such as singularity theory, see e.g., \cite{Helgason}).

\subsubsection{Definitions}

Formally, a \emph{root system} is a finite set $\Phi $ of non--zero vectors (%
\emph{roots}) spanning a finite--dimensional Euclidean space $V$ and
satisfying the following properties:

\begin{enumerate}
\item The only scalar multiples of a root $\alpha $ in $V$ which belong to $%
\Phi $ are $\alpha $ itself and -$\alpha $.

\item For every root $\alpha $ in $V$, the set $\Phi $ is symmetric under
reflection through the hyperplane of vectors perpendicular to $\alpha $.

\item If $\alpha $ and $\beta $ are vectors in $\Phi $, the projection of $%
2\beta $ onto the line through $\alpha $ is an integer multiple of $\alpha $.
\end{enumerate}

The \emph{rank} of a root system $\Phi $ is the dimension of $V$. Two root
systems may be combined by regarding the Euclidean spaces they span as
mutually orthogonal subspaces of a common Euclidean space. A root system
which does not arise from such a combination, such as the systems A$_{%
\mathrm{2}}$, B$_{\mathrm{2}}$, and G$_{\mathrm{2}}$ in Figure \ref{Root},
is said to be \emph{irreducible}.

Two irreducible root systems $(V_1,\Phi _1)$ and $(V_ 2,\Phi _2)$ are
considered to be the same if there is an invertible linear transformation $%
V_1\rightarrow V_2 $ which preserves distance up to a scale factor and which
sends $\Phi _{\mathrm{1}}$ to $\Phi _{\mathrm{2}}$.

The group of isometries of $V$ generated by reflections through hyperplanes
associated to the roots of $\Phi $ is called the Weyl group of $\Phi $ as it
acts faithfully on the finite set $\Phi $, the Weyl group is always finite.

\subsubsection{Classification}

\begin{figure}[tbh]
\centerline{\includegraphics[width=9cm]{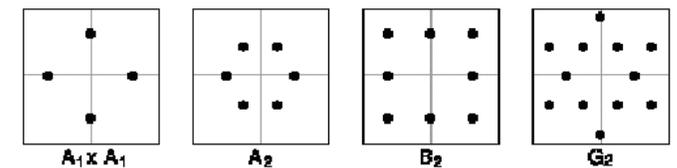}}
\caption{Classification of root systems of rank 2.}
\label{Root}
\end{figure}

It is not too difficult to classify the root systems of rank 2 (see Figure %
\ref{Root}).

Whenever $\Phi $ is a root system in $V$ and $W$ is a subspace of $V$
spanned by $\Psi =\Phi \cap W$, then $\Psi $ is a root system in $W$. Thus,
our exhaustive list of root systems of rank 2 shows the geometric
possibilities for any two roots in a root system. In particular, two such
roots meet at an angle of 0, 30, 45, 60, 90, 120, 135, 150, or 180 degrees.

In general, irreducible root systems are specified by a family (indicated by
a letter $A$ to $G$) and the rank (indicated by a subscript $n$). There are
four \emph{infinite families}:

\begin{itemize}
\item $A_{n}\,(n\geq 1),$ which corresponds to the special unitary group, $%
SU(n+1)$;

\item $B_{n}\,(n\geq 2),$ which corresponds to the special orthogonal group,
$SO(2n+1)$;

\item $C_{n}\,(n\geq 3),$ which corresponds to the symplectic group, $Sp(2n)
$;

\item $D_{n}\,(n\geq 4),$ which corresponds to the special orthogonal group,
$SO(2n)$,
\end{itemize}

as well as five \emph{exceptional cases}: $E_{\mathrm{6}},E_{\mathrm{7}},E_{%
\mathrm{8}},F_{\mathrm{4}},G_{\mathrm{2}}.$

\subsubsection{Dynkin Diagrams}

\begin{figure}[tbh]
\centerline{\includegraphics[width=10cm]{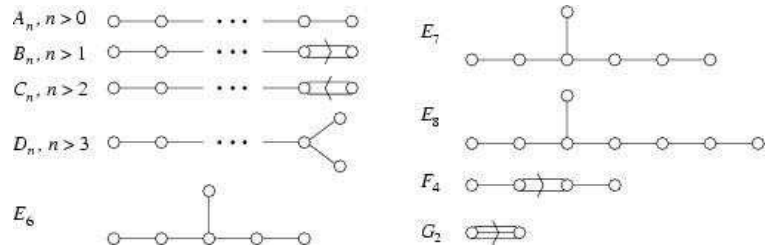}}
\caption{The problem of classifying irreducible root systems reduces to the
problem of classifying connected Dynkin diagrams.}
\label{Dynkin}
\end{figure}

A Dynkin diagram is a graph with a few different kinds of possible edges
(see Figure \ref{Dynkin}). The connected components of the graph correspond
to the irreducible subalgebras of $\mathfrak{g}$. So a simple Lie algebra's
Dynkin diagram has only one component. The rules are restrictive. In fact,
there are only certain possibilities for each component, corresponding to
the classification of semi--simple Lie algebras (see, e.g., \cite{Conway}).

The \emph{roots} of a complex Lie algebra form a lattice of rank $k$ in a
\textit{Cartan subalgebra} $\mathfrak{h\subset g}$, where $k$ is the Lie
algebra rank of $\mathfrak{g}$. Hence, the \emph{root lattice} can be
considered a lattice in $\mathbb{R}^{k}$. A vertex, or node, in the Dynkin
diagram is drawn for each \textit{Lie algebra simple root}, which
corresponds to a generator of the root lattice. Between two nodes $\alpha $\
and $\beta $, an edge is drawn if the simple roots are not perpendicular.
One line is drawn if the angle between them is $2\pi /3$, two lines if the
angle is $3\pi /4$, and three lines are drawn if the angle is $5\pi /6$.
There are no other possible angles between Lie algebra simple roots.
Alternatively, the number of lines $N$ between the simple roots $\alpha $\
and $\beta $ is given by
\begin{equation*}
N=A_{\alpha \beta }A_{\beta \alpha }=\frac{2\left\langle \alpha ,\beta
\right\rangle }{|\alpha |^{2}}\frac{2\left\langle \beta ,\alpha
\right\rangle }{|\beta |^{2}}=4\cos ^{2}\theta ,
\end{equation*}%
where $A_{\alpha \beta }=\frac{2\left\langle \alpha ,\beta \right\rangle }{%
|\alpha |^{2}}$\ is an entry in the \emph{Cartan matrix }$(A_{\alpha \beta })
$ (for details on Cartan matrix see, e.g., \cite{Helgason}). In a Dynkin
diagram, an arrow is drawn from the longer root to the shorter root (when
the angle is $3\pi /4$\ or $5\pi /6$).

Here are some properties of \textit{admissible Dynkin diagrams}:

\begin{enumerate}
\item A diagram obtained by removing a node from an admissible diagram is
admissible.

\item An admissible diagram has no loops.

\item No node has more than three lines attached to it.

\item A sequence of nodes with only two single lines can be collapsed to
give an admissible diagram.

\item The only connected diagram with a triple line has two nodes.
\end{enumerate}

A \textit{Coxeter--Dynkin diagram}, also called a \textit{Coxeter graph}, is
the same as a Dynkin diagram, but without the arrows. The Coxeter diagram is
sufficient to characterize the algebra, as can be seen by enumerating
connected diagrams.

The simplest way to recover a simple Lie algebra from its Dynkin diagram is
to first reconstruct its Cartan matrix $(A_{ij})$. The $i$th node and $j$th
node are connected by $A_{ij}A_{ji}$\ lines. Since $A_{ij}=0$\ iff $A_{ji}=0$%
, and otherwise $A_{ij}\in \{-3,-2,-1\}$, it is easy to find $A_{ij}$\ and $%
A_{ji}$, up to order, from their product. The arrow in the diagram indicates
which is larger. For example, if node 1 and node 2 have two lines between
them, from node 1 to node 2, then $A_{12}=-1$\ and $A_{21}=-2$.

However, it is worth pointing out that each simple Lie algebra can be
constructed concretely. For instance, the infinite families $A_{n}$, $B_{n}$%
, $C_{n}$, and $D_{n}$\ correspond to the special linear Lie algebra $%
\mathfrak{gl}(n+1,\mathbb{C})$, the odd orthogonal Lie algebra $\mathfrak{so}%
(2n+1,\mathbb{C})$, the symplectic Lie algebra $\mathfrak{sp}(2n,\mathbb{C})$%
, and the even orthogonal Lie algebra $\mathfrak{so}(2n,\mathbb{C})$. The
other simple Lie algebras are called \emph{exceptional Lie algebras}, and
have constructions related to the \emph{octonions}.

To prove this classification Theorem, one uses the angles between pairs of
roots to encode the root system in a much simpler combinatorial object, the
Dynkin diagram. The Dynkin diagrams can then be classified according to the
scheme given above.

To every root system is associated a corresponding Dynkin diagram.
Otherwise, the Dynkin diagram can be extracted from the root system by
choosing a \emph{base}, that is a subset $\Delta $ of $\Phi $ which is a
basis of $V$ with the special property that every vector in $\Phi $ when
written in the basis $\Delta $ has either all coefficients $\geq 0$ or else
all $\leq 0$.

The vertices of the Dynkin diagram correspond to vectors in $\Delta $. An
edge is drawn between each non--orthogonal pair of vectors; it is a double
edge if they make an angle of 135 degrees, and a triple edge if they make an
angle of 150 degrees. In addition, double and triple edges are marked with
an angle sign pointing toward the shorter vector.

Although a given root system has more than one base, the Weyl group acts
transitively on the set of bases. Therefore, the root system determines the
Dynkin diagram. Given two root systems with the same Dynkin diagram, we can
match up roots, starting with the roots in the base, and show that the
systems are in fact the same.

Thus the problem of classifying root systems reduces to the problem of
classifying possible Dynkin diagrams, and the problem of classifying
irreducible root systems reduces to the problem of classifying connected
Dynkin diagrams. Dynkin diagrams encode the inner product on $E$ in terms of
the basis $\Delta $, and the condition that this inner product must be
positive definite turns out to be all that is needed to get the desired
classification (see Figure \ref{Dynkin}).

In detail, the individual root systems can be realized case--by--case, as in
the following paragraphs:

\textbf{A$_{n}$.} Let $V$ be the subspace of $\mathbb{R}^{n\mathrm{+}\mathrm{%
1}}$ for which the coordinates sum to 0, and let $\Phi $ be the set of
vectors in $V$ of length $\sqrt{2}$ and with integer coordinates in $\mathbb{%
R}^{n\mathrm{+}\mathrm{1}}$. Such a vector must have all but two coordinates
equal to 0, one coordinate equal to $1$, and one equal to -$1$, so there are
$n^{\mathrm{2}}+n$ roots in all.

\textbf{B$_{n}$.} Let $V=\mathbb{R}^{n}$, and let $\Phi $ consist of all
integer vectors in $V$ of length 1 or $\sqrt{2}$. The total number of roots
is $2n^{\mathrm{2}}$.

\textbf{C$_{n}$:} Let $V=\mathbb{R}^{n}$, and let $\Phi $ consist of all
integer vectors in $V$ of $\sqrt{2}$ together with all vectors of the form $%
2\lambda $, where $\lambda $ is an integer vector of length 1. The total
number of roots is $2n^{\mathrm{2}}$. The total number of roots is $2n^{%
\mathrm{2}}$.

\textbf{D$_{n}$.} Let $V=\mathbb{R}^{n}$, and let $\Phi $ consist of all
integer vectors in $V$ of length $\sqrt{2}$. The total number of roots is $%
2n(n-1)$.

\textbf{E$_{n}$.} For $V_{\mathrm{8}}$, let $V=\mathbb{R}^{8}$, and let $%
E_{8}$ denote the set of vectors $\alpha $ of length $\sqrt{2}$ such that
the coordinates of $2\alpha $ are all integers and are either all even or
all odd. Then $E_{\mathrm{7}}$ can be constructed as the intersection of $E_{%
\mathrm{8}}$ with the hyperplane of vectors perpendicular to a fixed root $%
\alpha $ in $E_{\mathrm{8}}$, and $E_{\mathrm{6}}$ can be constructed as the
intersection of $E_{\mathrm{8}}$ with two such hyperplanes corresponding to
roots $\alpha $ and $\beta $ which are neither orthogonal to one another nor
scalar multiples of one another. The root systems $E_{\mathrm{6}}$, $E_{%
\mathrm{7}},$ and $E_{\mathrm{8}}$ have 72, 126, and 240 roots respectively.

\textbf{F$_{4}$.} For $F_{\mathrm{4}}$, let $V=\mathbb{R}^{4}$, and let $%
\Phi $ denote the set of vectors $\alpha $ of length 1 or $\sqrt{2}$ such
that the coordinates of $2\alpha $ are all integers and are either all even
or all odd. There are 48 roots in this system.

\textbf{G$_{2}$.} There are 12 roots in $G_{\mathrm{2}}$, which form the
vertices of a \emph{hexagram}.

\subsubsection{Irreducible Root Systems}

Irreducible root systems classify a number of related objects in Lie theory,
notably:

\begin{enumerate}
\item Simple complex Lie algebras;

\item Simple complex Lie groups;

\item Simply connected complex Lie groups which are simple modulo centers;
and

\item Simple compact Lie groups.
\end{enumerate}

In each case, the roots are non--zero weights of the adjoint representation.

A root system can also be said to describe a \emph{plant's root} and
associated systems.

\subsection{Simple and Semisimple Lie Groups and Algebras}

A \textit{simple Lie group} is a Lie group which is also a simple group.
These groups, and groups closely related to them, include many of the
so--called \emph{classical groups} of geometry, which lie behind projective
geometry and other geometries derived from it by the \textit{Erlangen
programme} of Felix Klein. They also include some \emph{exceptional groups},
that were first discovered by those pursuing the classification of simple
Lie groups. The exceptional groups account for many special examples and
configurations in other branches of mathematics. In particular the
classification of finite simple groups depended on a thorough prior
knowledge of the `exceptional' possibilities.

The complete listing of the simple Lie groups is the basis for the theory of
the semisimple Lie groups and reductive groups, and their representation
theory. This has turned out not only to be a major extension of the theory
of compact Lie groups (and their representation theory), but to be of basic
significance in mathematical physics.

Such groups are classified using the prior classification of the complex
simple Lie algebras. It has been shown that a simple Lie group has a simple
Lie algebra that will occur on the list given there, once it is complexified
(that is, made into a complex vector space rather than a real one). This
reduces the classification to two further matters.

The groups $SO(p,q,\mathbb{R})$ and $SO(p+q,\mathbb{R})$, for example, give
rise to different real Lie algebras, but having the same Dynkin diagram. In
general there may be different \emph{real forms} of the same complex Lie
algebra.

Secondly, the Lie algebra only determines uniquely the simply connected
(universal) cover $G^{\ast }$ of the component containing the identity of a
Lie group $G$. It may well happen that $G^{\ast }$ is not actually a simple
group, for example having a non--trivial center. We have therefore to worry
about the global topology, by computing the fundamental group of $G$ (an
Abelian group: a Lie group is an $H-$space). This was done by Elie Cartan.

For an example, take the special orthogonal groups in even dimension. With $%
-I$ a scalar matrix in the center, these are not actually simple groups; and
having a two--fold spin cover, they aren't simply--connected either. They
lie `between' $G^{\ast }$ and $G$, in the notation above.

Recall that a \emph{semisimple module} is a module in which each submodule
is a direct summand. In particular, a \textit{semisimple representation} is
completely reducible, i.e., is a direct sum of irreducible representations
(under a descending chain condition). Similarly, one speaks of an Abelian
category as being semisimple when every object has the corresponding
property. Also, a semisimple ring is one that is semisimple as a module over
itself.

A \emph{semisimple matrix} is diagonalizable over any algebraically closed
field containing its entries. In practice this means that it has a diagonal
matrix as its Jordan normal form.

A \emph{Lie algebra} $\mathfrak{g}$\ is called \emph{semisimple} when it is
a direct sum of \emph{simple Lie algebras}, i.e., non--trivial Lie algebras $%
\mathfrak{L}$ whose only ideals are $\{0\}$ and $\mathfrak{L}$ itself. An
equivalent condition is that the \textit{Killing form}
\begin{equation*}
\mathcal{B}(X,Y)=\limfunc{Tr}(Ad(X)\,Ad(Y))
\end{equation*}
is non--degenerate \cite{Schafer}. The following properties can be proved
equivalent for a finite--dimensional algebra $\mathfrak{L}$ over a field of
characteristic $0$:

1. $\mathfrak{L}$ is semisimple.

2. $\mathfrak{L}$ has no nonzero Abelian ideal.

3. $\mathfrak{L}$ has zero radical (the radical is the biggest solvable
ideal).

4. Every representation of $\mathfrak{L}$ is fully reducible, i.e., is a sum
of irreducible representations.

5. $\mathfrak{L}$ is a (finite) direct product of simple Lie algebras (a Lie
algebra is called simple if it is not Abelian and has no nonzero ideal ).

A \textit{connected Lie group} is called \emph{semisimple} when its Lie
algebra is semisimple; and the same holds for algebraic groups. Every finite
dimensional representation of a semisimple Lie algebra, Lie group, or
algebraic group in characteristic $0$ is semisimple, i.e., completely
reducible, but the converse is not true. Moreover, in characteristic $p>0$,
semisimple Lie groups and Lie algebras have finite dimensional
representations which are not semisimple. An element of a semisimple Lie
group or Lie algebra is itself semisimple if its image in every
finite--dimensional representation is semisimple in the sense of matrices.

Every semisimple Lie algebra $\mathfrak{g}$\ can be classified by its Dynkin
diagram\emph{\ }\cite{Helgason}.

\section{Some Classical Examples of Lie Groups}

\subsection{Galilei Group}

The \textit{Galilei group} is the group of transformations in space and time
that connect those Cartesian systems that are termed `inertial frames' in
Newtonian mechanics. The most general relationship between two such frames
is the following. The origin of the time scale in the inertial frame $%
S^{\prime }$ may be shifted compared with that in $S$; the orientation of
the Cartesian axes in $S^{\prime }$ may be different from that in $S$; the
origin $O$ of the Cartesian frame in $S^{\prime }$ may be moving relative to
the origin $O$ in $S$ at a uniform velocity. The transition from $S$ to $%
S^{\prime }$ involves ten parameters; thus the Galilei group is a ten
parameter group. The basic assumption inherent in Galilei--Newtonian
relativity is that there is an absolute time scale, so that the only way in
which the time variables used by two different `inertial observers' could
possibly differ is that the zero of time for one of them may be shifted
relative to the zero of time for the other (see, e.g. \cite%
{Arnold,GaneshSprBig,GaneshADG}).

Galilei space--time structure involves the following three elements:

\begin{enumerate}
\item \emph{World}, as a 4D affine space $A^{4}$. The points of $A^{4}$ are
called \emph{world points} or \emph{events}. The parallel transitions of the
world $A^{4}$ form a linear (i.e., Euclidean) space $\mathbb{R}^{4}$.

\item \emph{Time}, as a linear map $t:\mathbb{R}^{4}\rightarrow \mathbb{R}$
of the linear space of the world parallel transitions onto the real `time
axes'. Time interval from the event $a\in A^{4}$ to $b\in A^{4}$ is called
the number $t(b-a)$; if $t(b-a)=0$ then the events $a$ and $b$ are called
synchronous. The set of all mutually synchronous events consists a 3D affine
space $A^{3}$, being a subspace of the world $A^{4}$. The kernel of the
mapping $t$ consists of the parallel transitions of $A^{4}$ translating
arbitrary (and every) event to the synchronous one; it is a linear 3D
subspace $\mathbb{R}^{3}$ of the space $\mathbb{R}^{4}$.

\item \emph{Distance} (\emph{metric}) between the synchronous events,
\begin{equation*}
\rho (a,b)=\parallel a-b\parallel ,\qquad \text{for all}\quad a,b\in \mathop %
A\nolimits^{3},
\end{equation*}%
given by the scalar product in $\mathbb{R}^{3}$. The distance transforms
arbitrary space of synchronous events into the well known 3D Euclidean space
$E^{3}$.
\end{enumerate}

The space $A^4$, with the Galilei space--time structure on it, is called
Galilei space. Galilei group is the group of all possible transformations of
the Galilei space, preserving its structure. The elements of the Galilei
group are called Galilei transformations. Therefore, Galilei transformations
are affine transformations of the world $A^4$ preserving the time intervals
and distances between the synchronous events.

The direct product $\mathbb{R}\times \mathbb{R}^{3}$, of the time axes with
the 3D linear space R3 with a fixed Euclidean structure, has a natural
Galilei structure. It is called Galilei coordinate system.

\subsection{General Linear Group}

The group of linear isomorphisms of $\mathbb{R}^{n}$ to $\mathbb{R}^{n}$ is
a Lie group of dimension $n^{2}$, called the \textit{general linear group}
and denoted $Gl(n,\mathbb{R})$. It is a smooth manifold, since it is a
subset of the vector space $L(\mathbb{R}^{n},\mathbb{R}^{n})$ of all linear
maps of $\mathbb{R}^{n}$ to $\mathbb{R}^{n}$, as $Gl(n,\mathbb{R})$ is the
inverse image of $\mathbb{R}\backslash \{0\}$ under the continuous map $%
A\mapsto \det A$ of $L(\mathbb{R}^{n},\mathbb{R}^{n})$ to $\mathbb{R}$. The
group operation is composition (see, e.g. \cite%
{Arnold,Abraham,GaneshSprBig,GaneshADG}).
\begin{equation*}
(A,B)\in Gl(n,\mathbb{R})\times Gl(n,\mathbb{R})\mapsto A\circ B\in Gl(n,%
\mathbb{R})
\end{equation*}%
and the inverse map is
\begin{equation*}
A\in Gl(n,\mathbb{R})\mapsto A^{-1}\in Gl(n,\mathbb{R}).
\end{equation*}%
If we choose a basis in $\mathbb{R}^{n}$, we can represent each element $%
A\in Gl(n,\mathbb{R})$ by an invertible $(n\times n)-$-matrix. The group
operation is then matrix multiplication and the inversion is matrix
inversion. The identity is the identity matrix $I_{n}$. The group operations
are smooth since the formulas for the product and inverse of matrices are
smooth in the matrix components.

The Lie algebra of $Gl(n,\mathbb{R})$ is $\mathfrak{gl}(n)$, the vector
space $L(\mathbb{R}^{n},\mathbb{R}^{n})$\ of all linear transformations of $%
\mathbb{R}^{n}$, with the commutator bracket
\begin{equation*}
\lbrack A,B]=AB-BA.
\end{equation*}%
For every $A\in L(\mathbb{R}^{n},\mathbb{R}^{n})$,
\begin{equation*}
\gamma _{A}:t\in \mathbb{R\mapsto }\gamma _{A}(t)=\sum_{i=0}^{\infty }\frac{%
t^{i}}{i!}A^{i}\in Gl(n,\mathbb{R})
\end{equation*}%
is a one--parameter subgroup of $Gl(n,\mathbb{R})$, because
\begin{equation*}
\gamma _{A}(0)=I\qquad \text{and\qquad }\dot{\gamma}_{A}(t)=\sum_{i=0}^{%
\infty }\frac{t^{i-1}}{(i-1)!}A^{i}=\gamma _{A}(t)\,A
\end{equation*}
Hence $\gamma _{A}$ is an integral curve of the left--invariant
vector--field $X_{A}$. Therefore, the exponential map is given by
\begin{equation*}
\exp :A\in L(\mathbb{R}^{n},\mathbb{R}^{n})\mapsto \exp (A)\equiv \mathrm{e}%
^{A}=\gamma _{A}(1)=\sum_{i=0}^{\infty }\frac{A^{i}}{i!}\in Gl(n,\mathbb{R}).
\end{equation*}

For each $A\in Gl(n,\mathbb{R})$ the corresponding adjoint map
\begin{equation*}
Ad_{A}:L(\mathbb{R}^{n},\mathbb{R}^{n})\rightarrow L(\mathbb{R}^{n},\mathbb{R%
}^{n})
\end{equation*}%
is given by
\begin{equation*}
Ad_{A}B=A\cdot B\cdot A^{-1}.
\end{equation*}

\subsection{Rotational Lie Groups in Human/Humanoid Biomechanics}

Local kinematics at each rotational robot or (synovial) human joint, is
defined as a \emph{group action} of an $n$D constrained rotational Lie group
$SO(n)$ on the Euclidean space $\mathbb{R}^{n}$. In particular, there is an
action of $SO(2) -$-group in uniaxial human joints (cylindrical, or \emph{%
hinge joints}, like knee and elbow) and an action of $SO(3)-$-group in
three--axial human joints (spherical, or \emph{ball--and--socket joints},
like hip, shoulder, neck, wrist and ankle). In both cases, $SO(n)$ acts,
with its operators of rotation, on the vector $x=\{x^{\mu }\},\,(i=1,2,3)$
of external, Cartesian coordinates of the parent body--segment, depending,
at the same time, on the vector $q=\{q^{s}\},\,(s=1,\cdots ,n)$ on $n$
group--parameters, i.e., joint angles (see \cite%
{SIAM,GaneshSprSml,GaneshSprBig,GaneshADG}).

Each joint rotation $R\in SO(n)$ defines a map
\begin{equation*}
R:x^{\mu }\mapsto \dot{x}{}^{\mu }, \qquad R(x^{\mu },q^{s})=R_{q^{s}}x^{\mu
},
\end{equation*}
where $R_{q^{s}}\in SO(n)$ are joint group operators. The vector $%
v=\{v_{s}\},\,(s=1,\cdots ,n)$ of $n$ infinitesimal generators of these
rotations, i.e., joint angular velocities, given by
\begin{equation*}
v_{s}=-\left[\frac{\partial R(x^{\mu },q^{s})}{\partial q^{s}}\right]_{q=0}%
\frac{\partial }{\partial x^{\mu }}
\end{equation*}
constitute an $n$D Lie algebra $\mathfrak{so}(n)$ corresponding to the joint
rotation group $SO(n)$. Conversely, each joint group operator $R_{q^{s}}$,
representing a one--parameter subgroup of $SO(n)$, is defined as the
exponential map of the corresponding joint group generator $v_{s}$
\begin{equation*}
R_{q^{s}}=\exp (q^{s}v_{s})
\end{equation*}
This exponential map represents a solution of the joint operator
differential equation in the joint group--parameter space $\{q^{s}\}$
\begin{equation*}
\frac{dR_{q^{s}}}{dq^{s}}=v_{s}R_{q^{s}}.
\end{equation*}

\subsubsection{Uniaxial Group of Joint Rotations}

The uniaxial joint rotation in a single Cartesian plane around a
perpendicular axis, e.g., $xy-$plane about the $z$ axis, by an internal
joint angle $\theta ,$ leads to the following transformation of the joint
coordinates:
\begin{equation*}
x^{\prime }=x\cos \theta -y\sin \theta ,\qquad y^{\prime }=x\sin \theta
+y\cos \theta .
\end{equation*}
In this way, the joint $SO(2)-$group, given by
\begin{equation*}
SO(2)=\left\{ R_{\theta }=\left(
\begin{array}{cc}
\cos \theta & -\sin \theta \\
\sin \theta & \cos \theta%
\end{array}%
\right) |\theta \in \lbrack 0,2\pi ]\right\},
\end{equation*}%
acts in a canonical way on the Euclidean plane $\mathbb{R}^{2}$ by
\begin{equation*}
SO(2)=\left( \left(
\begin{array}{cc}
\cos \theta & -\sin \theta \\
\sin \theta & \cos \theta%
\end{array}%
\right) ,\left(
\begin{array}{c}
x \\
y%
\end{array}%
\right) \right) \longmapsto \left(
\begin{array}{cc}
x\cos \theta & -y\sin \theta \\
x\sin \theta & y\cos \theta%
\end{array}%
\right).
\end{equation*}
Its associated Lie algebra $\mathfrak{so}(2)$ is given by
\begin{equation*}
\mathfrak{so}(2)=\left\{ \left(
\begin{array}{cc}
0 & -t \\
t & 0%
\end{array}%
\right) |t\in \mathbb{R}\right\} ,
\end{equation*}%
since the curve $\gamma _{\theta }\in SO(2)$ given by
\begin{equation*}
\gamma _{\theta }:t\in \mathbb{R}\longmapsto \gamma _{\theta }(t)=\left(
\begin{array}{cc}
\cos t\theta & -\sin t\theta \\
\sin t\theta & \cos t\theta%
\end{array}%
\right) \in SO(2),
\end{equation*}%
passes through the identity $I_{2}=\left(
\begin{array}{cc}
1 & 0 \\
0 & 1%
\end{array}%
\right) $ and then
\begin{equation*}
\left. \frac{d}{dt}\right\vert _{t=0}\gamma _{\theta }(t)=\left(
\begin{array}{cc}
0 & -\theta \\
\theta & 0%
\end{array}%
\right) ,
\end{equation*}%
so that $I_{2}$ is a basis of $\mathfrak{so}(2)$, since $\dim \left(
SO(2)\right) =1$.

The \emph{exponential map} $\exp :\mathfrak{so}(2)\rightarrow SO(2)$ is
given by
\begin{equation*}
\exp \left(
\begin{array}{cc}
0 & -\theta \\
\theta & 0%
\end{array}
\right) =\gamma _{\theta }(1)=\left(
\begin{array}{cc}
\cos t\theta & -\sin t\theta \\
\sin t\theta & \cos t\theta%
\end{array}
\right).
\end{equation*}

The \emph{infinitesimal generator} of the action of $SO(2)$ on $\mathbb{R}%
^{2}$, i.e., joint angular velocity $v,$ is given by
\begin{equation*}
v=-y\frac{\partial }{\partial x}+x\frac{\partial }{\partial y},
\end{equation*}
since
\begin{equation*}
v_{\mathbb{R}^{2}}\left( x,y\right) =\left. \frac{d}{dt}\right\vert
_{t=0}\exp (tv)\left( x,y\right) =\left. \frac{d}{dt}\right\vert
_{t=0}\left(
\begin{array}{cc}
\cos tv & -\sin tv \\
\sin tv & \cos tv%
\end{array}%
\right) \left(
\begin{array}{c}
x \\
y%
\end{array}%
\right) .
\end{equation*}

The \emph{momentum map} $J:T^{\ast }\mathbb{R}^{2}\rightarrow \mathbb{R}$
associated to the lifted action of $SO(2)$ on $T^{\ast }\mathbb{R}^{2}\simeq
\mathbb{R}^{4}$ is given by
\begin{eqnarray*}
J\left( x,y,p_{1},p_{2}\right) &=&xp_{y}-yp_{x},\qquad \text{since} \\
J\left( x,y,p_{x},p_{y}\right) (\xi ) &=&(p_{x}dx+p_{y}dy)(v_{\mathbb{R}%
^{2}})=-vp_{x}y+-vp_{y}x.
\end{eqnarray*}

The Lie group $SO(2)$ acts on the symplectic manifold $(\mathbb{R}%
^{4},\omega =dp_{x}\wedge dx+dp_{y}\wedge dx)$ by
\begin{eqnarray*}
&&\qquad\phi \left( \left(
\begin{array}{cc}
\cos \theta & -\sin \theta \\
\sin \theta & \cos \theta%
\end{array}
\right) ,\left( x,y,p_{x},p_{y}\right) \right) \\
&=&\left( x\cos \theta -y\sin \theta ,\,x\sin \theta +y\cos \theta
,\,p_{x}\cos \theta -p_{y}\sin \theta ,\,p_{x}\sin \theta +p_{y}\cos \theta
\right) .
\end{eqnarray*}%
\bigbreak

\subsubsection{Three--Axial Group of Joint Rotations}

The three--axial $SO(3)-$group of human--like joint rotations depends on
three parameters, Euler joint angles $q^{i}=(\varphi ,\psi ,\theta ),$\
defining the rotations about the Cartesian coordinate triedar $(x,y,z)$
placed at the joint pivot point. Each of the Euler angles are defined in the
constrained range $(-\pi ,\pi )$, so the joint group space is a constrained
sphere of radius $\pi $ (see \cite{SIAM,GaneshSprSml,GaneshSprBig,GaneshADG}%
).

Let $G=SO(3)=\{A\in \mathcal{M}_{3\times 3}(\mathbb{R}):A^{t}A=I_{3},\det
(A)=1\}$ be the group of rotations in $\mathbb{R}^{3}$. It is a Lie group
and $\dim(G)=3$. Let us isolate its one--parameter joint subgroups, i.e.,
consider the three operators of the finite joint rotations $R_{\varphi
},R_{\psi },R_{\theta }\in SO(3),$ given by
\begin{equation*}
R_{\varphi } =\left[
\begin{array}{ccc}
1 & 0 & 0 \\
0 & \cos \varphi & -\sin \varphi \\
0 & \sin \varphi & \cos \varphi%
\end{array}
\right] , ~~ R_{\psi } =\left[
\begin{array}{ccc}
\cos \psi & 0 & \sin \psi \\
0 & 1 & 0 \\
-\sin \psi & 0 & \cos \psi%
\end{array}
\right] , ~~ R_{\theta } =\left[
\begin{array}{ccc}
\cos \theta & -\sin \theta & 0 \\
\sin \theta & \cos \theta & 0 \\
0 & 0 & 1%
\end{array}
\right]
\end{equation*}%
\noindent corresponding respectively to rotations about $x-$axis by an angle
$\varphi ,$ about $y-$axis by an angle $\psi ,$ and about $z-$axis by an
angle $\theta $.

The total three--axial joint rotation $A$ is defined as the product of above
one--parameter rotations $R_{\varphi },R_{\psi },R_{\theta },$ i.e., $%
A=R_{\varphi }\cdot R_{\psi }\cdot R_{\theta }$ is equal\footnote{%
Note that this product is noncommutative, so it really depends on the order
of multiplications.}
\begin{equation*}
A=\left[
\begin{array}{ccc}
\cos \psi \cos \varphi -\cos \theta \sin \varphi \sin \psi & \cos \psi \cos
\varphi +\cos \theta \cos \varphi \sin \psi & \sin \theta \sin \psi \\
-\sin \psi \cos \varphi -\cos \theta \sin \varphi \sin \psi & -\sin \psi
\sin \varphi +\cos \theta \cos \varphi \cos \psi & \sin \theta \cos \psi \\
\sin \theta \sin \varphi & -\sin \theta \cos \varphi & \cos \theta%
\end{array}
\right] .
\end{equation*}
However, the order of these matrix products matters: different order
products give different results, as the matrix product is \textit{%
noncommutative product}. This is the reason why Hamilton's \textit{%
quaternions}\footnote{%
Recall that the set of Hamilton's \textit{quaternions} $\mathbb{H}$
represents an extension of the set of complex numbers $\mathbb{C}$. We can
compute a rotation about the unit vector, $\mathbf{u}$ by an angle $\theta $%
. The quaternion $q$ that computes this rotation is
\begin{equation*}
q=\left( \cos \frac{\theta }{2}~,~u\sin \frac{\theta }{2}\right) .
\end{equation*}%
} are today commonly used to parameterize the $SO(3)-$group, especially in
the field of 3D computer graphics.

The one--parameter rotations $R_{\varphi },R_{\psi },R_{\theta }$ define
curves in $SO(3)$ starting from $I_{3}={\small \left(
\begin{array}{ccc}
1 & 0 & 0 \\
0 & 1 & 0 \\
0 & 0 & 1%
\end{array}
\right)} .$ Their derivatives in $\varphi =0,\psi =0$ and $\theta =0\,\ $%
belong to the associated \emph{tangent Lie algebra} $\mathfrak{so}(3)$. That
is the corresponding infinitesimal generators of joint rotations -- joint
angular velocities $v_{\varphi },v_{\psi },v_{\theta }\in \mathfrak{so}(3)$
-- are respectively given by
\begin{eqnarray*}
v_{\varphi } &=&{\small \left[
\begin{array}{ccc}
0 & 0 & 0 \\
0 & 0 & -1 \\
0 & 1 & 0%
\end{array}%
\right]} =-y\frac{\partial }{\partial z}+z\frac{\partial }{\partial y}%
,\qquad v_{\psi }={\small \left[
\begin{array}{ccc}
0 & 0 & 1 \\
0 & 0 & 0 \\
-1 & 0 & 0%
\end{array}%
\right]} =-z\frac{\partial }{\partial x}+x\frac{\partial }{\partial z}, \\
v_{\theta } &=&{\small \left[
\begin{array}{ccc}
0 & -1 & 0 \\
1 & 1 & 0 \\
0 & 0 & 0%
\end{array}%
\right]} =-x\frac{\partial }{\partial y}+y\frac{\partial }{\partial x}.
\end{eqnarray*}
Moreover, the elements are linearly independent and so
\begin{equation*}
\mathfrak{so}(3)=\left\{ \left[
\begin{array}{ccc}
0 & -a & b \\
a & 0 & -\gamma \\
-b & \gamma & 0%
\end{array}
\right] |a,b,\gamma \in \mathbb{R}\right\}.
\end{equation*}
The Lie algebra $\mathfrak{so}(3)$ is identified with $\mathbb{R}^{3}$ by
associating to each $v=(v_{\varphi },v_{\psi },v_{\theta })\in \mathbb{R}%
^{3} $ the matrix $v\in \mathfrak{so}(3)$ given by $v={\small \left[
\begin{array}{ccc}
0 & -a & b \\
a & 0 & -\gamma \\
-b & \gamma & 0%
\end{array}
\right]} . $ Then we have the following identities:

\begin{enumerate}
\item $\widehat{u\times v}=[\hat{u},v]$; and

\item $u\cdot v=-\frac{1}{2}\limfunc{Tr}(\hat{u}\cdot v)$.
\end{enumerate}

The exponential map $\exp :\mathfrak{so}(3)\rightarrow SO(3)$ is given by
\emph{Rodrigues relation}
\begin{equation*}
\exp (v)=I+\frac{\sin \left\Vert v\right\Vert }{\left\Vert v\right\Vert }v+%
\frac{1}{2}\left( \frac{\sin \frac{\left\Vert v\right\Vert }{2}}{\frac{%
\left\Vert v\right\Vert }{2}}\right) ^{2}v^{2}
\end{equation*}
where the norm $\left\Vert v\right\Vert $ is given by
\begin{equation*}
\left\Vert v\right\Vert =\sqrt{(v^{1})^{2}+(v^{2})^{2}+(v^{3})^{2}}.
\end{equation*}

The the dual, \emph{cotangent Lie algebra} $\mathfrak{so}(3)^{\ast },$
includes the three joint angular momenta $p_{\varphi },p_{\psi },p_{\theta
}\in \mathfrak{so}(3)^{\ast }$, derived from the joint velocities $v$ by
multiplying them with corresponding moments of inertia.

Note that the parameterization of $SO(3)-$rotations is the subject of
continuous research and development in many theoretical and applied fields
of mechanics, such as rigid body, structural, and multibody dynamics,
robotics, spacecraft attitude dynamics, navigation, image processing, etc.

\subsubsection{The Heavy Top}

Consider a rigid body moving with a fixed point but under the influence of
gravity. This problem still has a configuration space $SO(3)$, but the
symmetry group is only the circle group $S^{1}$, consisting of rotations
about the direction of gravity. One says that gravity has broken the
symmetry from $SO(3)$ to $S^{1}$. This time, eliminating the $S^{1}$
symmetry mysteriously leads one to the larger Euclidean group $SE(3)$ of
rigid motion of $\mathbb{R}^{3}$. Conversely, we can start with $SE(3)$ as
the configuration space for the rigid--body and `reduce out' translations to
arrive at $SO(3)$ as the configuration space. The equations of motion for a
rigid body with a fixed point in a gravitational field give an interesting
example of a system that is Hamiltonian. The underlying Lie algebra consists
of the algebra of infinitesimal Euclidean motions in $\mathbb{R}^{3}$ (see
\cite{Arnold,Abraham,GaneshSprBig,GaneshADG}).

The basic phase--space we start with is again $T^{\ast }SO(3)$,
parameterized by Euler angles and their conjugate momenta. In these
variables, the equations are in canonical Hamiltonian form. However, the
presence of gravity breaks the symmetry, and the system is no longer $SO(3)$
invariant, so it cannot be written entirely in terms of the body angular
momentum $p$. One also needs to keep track of $\Gamma $, the `direction of
gravity' as seen from the body. This is defined by $\Gamma =A^{-1}k$, where $%
k$ points upward and $A$ is the element of $SO(3)$ describing the current
configuration of the body. The equations of motion are
\begin{eqnarray*}
\dot{p}_{1} &=&\frac{I_{2}-I_{3}}{I_{2}I_{3}}p_{2}p_{3}+Mgl(\Gamma ^{2}\chi
^{3}-\Gamma ^{3}\chi ^{2}), \\
\dot{p}_{2} &=&\frac{I_{3}-I_{1}}{I_{3}I_{1}}p_{3}p_{1}+Mgl(\Gamma ^{3}\chi
^{1}-\Gamma ^{1}\chi ^{3}), \\
\dot{p}_{3} &=&\frac{I_{1}-I_{2}}{I_{1}I_{2}}p_{1}p_{2}+Mgl(\Gamma ^{1}\chi
^{2}-\Gamma ^{2}\chi ^{1}), \\
&&\text{and\qquad }\dot{\Gamma}\;=\;\Gamma \times \Omega ,
\end{eqnarray*}%
where $\Omega $ is the body's angular velocity vector, $I_{1},I_{2},I_{3}$
are the body's principal moments of inertia, $M$ is the body's mass, $g$ is
the acceleration of gravity, $\chi $ is the body fixed unit vector on the
line segment connecting the fixed point with the body's center of mass, and $%
l$ is the length of this segment.

\subsection{Euclidean Groups of Rigid Body Motion}

In this subsection we give description of two most important Lie groups in
classical mechanics in 2D and 3D, $SE(2)$ and $SE(3)$, respectively (see
\cite{Marsden,GaneshSprBig,GaneshADG}).

\subsubsection{Special Euclidean Group $SE(2)$ in the Plane}

\label{SE(2)}

The motion in uniaxial human joints is naturally modelled by the \emph{%
special Euclidean group in the plane}, $SE(2)$. It consists of all
transformations of $\mathbb{R}^{2}$ of the form $Az+a$, where $z,a\in
\mathbb{R}^{2}$, and%
\begin{equation*}
A\in SO(2)=\left\{ \text{matrices of the form }\left(
\begin{array}{cc}
\cos \theta & -\sin \theta \\
\sin \theta & \cos \theta%
\end{array}%
\right) \right\} .
\end{equation*}%
In other words, group $SE(2)$ consists of matrices of the form:\newline
$(R_{\theta },a) ={\small \left(
\begin{array}{cc}
R_{\theta } & a \\
0 & I%
\end{array}%
\right)} ,$ where $a\in \mathbb{R}^{2}$ and $R_{\theta }$ is the rotation
matrix:\newline
$R_{\theta }={\small \left(
\begin{array}{cc}
\cos \theta & -\sin \theta \\
\sin \theta & \cos \theta%
\end{array}%
\right)},$ while $I$\ is the $3\times 3$ identity matrix. The inverse $%
\left( R_{\theta },a\right) ^{-1}$\ is given by%
\begin{equation*}
\left( R_{\theta },a\right) ^{-1}=\left(
\begin{array}{cc}
R_{\theta } & a \\
0 & I%
\end{array}%
\right) ^{-1}=\left(
\begin{array}{cc}
R_{-\theta } & -R_{-\theta }a \\
0 & I%
\end{array}%
\right) .
\end{equation*}%
The Lie algebra $\mathfrak{se}(2)$ of $SE(2)$ consists of $3\times 3$ block
matrices of the form%
\begin{equation*}
\left(
\begin{array}{cc}
-\xi J & v \\
0 & 0%
\end{array}%
\right) ,\qquad \text{where}\qquad J=\left(
\begin{array}{cc}
0 & 1 \\
-1 & 0%
\end{array}%
\right), \qquad (J^{T}=J^{-1}=-J),
\end{equation*}%
with the usual commutator bracket. If we identify $\mathfrak{se}(2)$ with $%
\mathbb{R}^{3}$ by the isomorphism%
\begin{equation*}
\left(
\begin{array}{cc}
-\xi J & v \\
0 & 0%
\end{array}%
\right) \in \mathfrak{se}(2)\longmapsto (\xi ,v)\in \mathbb{R}^{3},
\end{equation*}%
then the expression for the Lie algebra bracket becomes

\begin{equation*}
\lbrack (\xi ,v_{1},v_{2}),(\zeta ,w_{1},w_{2})]=(0,\zeta v_{2}-\xi
w_{2},\xi w_{1}-\zeta v1)=(0,\xi J^{T}w-\zeta J^{T}v),
\end{equation*}%
where $v=(v_{1},v_{2})$ and $w=(w_{1},w_{2})$.

The \emph{adjoint group action} of%
\begin{equation*}
\left( R_{\theta },a\right) \left(
\begin{array}{cc}
R_{\theta } & a \\
0 & I%
\end{array}%
\right) \qquad \text{on}\qquad (\xi ,v)=\left(
\begin{array}{cc}
-\xi J & v \\
0 & 0%
\end{array}%
\right)
\end{equation*}%
is given by the \emph{group conjugation},%
\begin{equation*}
\left(
\begin{array}{cc}
R_{\theta } & a \\
0 & I%
\end{array}%
\right) \left(
\begin{array}{cc}
-\xi J & v \\
0 & 0%
\end{array}%
\right) \left(
\begin{array}{cc}
R_{-\theta } & -R_{-\theta }a \\
0 & I%
\end{array}%
\right) =\left(
\begin{array}{cc}
-\xi J & \xi Ja+R_{\theta }v \\
0 & 0%
\end{array}%
\right) ,
\end{equation*}%
or, in coordinates,
\begin{equation}
Ad_{\left( R_{\theta },a\right) }(\xi ,v)=(\xi ,\xi Ja+R_{\theta }v).
\label{adse2}
\end{equation}

In proving (\ref{adse2}) we used the identity $R_{\theta }J=JR_{\theta }$.
Identify the dual algebra, $\mathfrak{se}(2)^{\ast },$ with matrices of the
form $\left(
\begin{array}{cc}
\frac{\mu }{2}J & 0 \\
\alpha & 0%
\end{array}%
\right),$ via the nondegenerate pairing given by the trace of the product.
Thus, $\mathfrak{se}(2)^{\ast }$ is isomorphic to $\mathbb{R}^{3}$ via%
\begin{equation*}
\left(
\begin{array}{cc}
\frac{\mu }{2}J & 0 \\
\alpha & 0%
\end{array}%
\right) \in \mathfrak{se}(2)^{\ast }\longmapsto (\mu ,\alpha )\in \mathbb{R}%
^{3},
\end{equation*}%
so that in these coordinates, the pairing between $\mathfrak{se}(2)^{\ast }$
and $\mathfrak{se}(2)$ becomes%
\begin{equation*}
\left\langle (\mu ,\alpha ),(\xi ,v)\right\rangle =\mu \xi +\alpha \cdot v,
\end{equation*}%
that is, the usual dot product in $\mathbb{R}^{3}$. The \emph{coadjoint
group action} is thus given by%
\begin{equation}
Ad_{\left( R_{\theta },a\right) ^{-1}}^{\ast }(\mu ,\alpha )=(\mu -R_{\theta
}\alpha \cdot Ja+R_{\theta }\alpha ).  \label{ad*se2}
\end{equation}

Formula (\ref{ad*se2}) shows that the coadjoint orbits are the cylinders $%
T^{\ast }S_{\alpha }^{1}=\{(\mu ,\alpha )|\left\Vert \alpha \right\Vert=%
\text{const}\}$ if $\alpha \neq 0$ together with the points are on the $\mu
- $axis. The canonical cotangent bundle projection $\pi :T^{\ast }S_{\alpha
}^{1}\rightarrow S_{\alpha }^{1}$ is defined as $\pi (\mu ,\alpha )=\alpha .$

\subsubsection{Special Euclidean Group $SE(3)$ in the 3D Space}

\label{SE(3)}

The most common group structure in human biodynamics is the \emph{special
Euclidean group in 3D space}, $SE(3)$. Briefly, the Euclidean $SE(3)-$group
is defined as a semidirect (noncommutative) product of 3D rotations and 3D
translations, $SE(3):=SO(3)\rhd \mathbb{R}^{3}$ (see \cite%
{Marsden,GaneshSprBig,GaneshADG}). Its most important subgroups are the
following:\newline

{{\frame{$%
\begin{array}{cc}
\mathbf{Subgroup} & \mathbf{Definition} \\ \hline
\begin{array}{c}
SO(3),\text{ group of rotations} \\
\text{in 3D (a spherical joint)}%
\end{array}
&
\begin{array}{c}
\text{Set of all proper orthogonal } \\
3\times 3-\text{rotational matrices}%
\end{array}
\\ \hline
\begin{array}{c}
SE(2),\text{ special Euclidean group} \\
\text{in 2D (all planar motions)}%
\end{array}
&
\begin{array}{c}
\text{Set of all }3\times 3-\text{matrices:} \\
\left[
\begin{array}{ccc}
\cos \theta & \sin \theta & r_{x} \\
-\sin \theta & \cos \theta & r_{y} \\
0 & 0 & 1%
\end{array}%
\right]%
\end{array}
\\ \hline
\begin{array}{c}
SO(2),\text{ group of rotations in 2D} \\
\text{subgroup of }SE(2)\text{--group} \\
\text{(a revolute joint)}%
\end{array}
&
\begin{array}{c}
\text{Set of all proper orthogonal } \\
2\times 2-\text{rotational matrices} \\
\text{ included in }SE(2)-\text{group}%
\end{array}
\\ \hline
\begin{array}{c}
\mathbb{R}^{3},\text{ group of translations in 3D} \\
\text{(all spatial displacements)}%
\end{array}
& \text{Euclidean 3D vector space}%
\end{array}%
$}}} \bigskip

\paragraph{Lie Group $SE(3)$ and Its Lie Algebra}

An element of $SE(3)$ is a pair $(A,a)$ where $A\in SO(3)$ and $a\in \mathbb{%
R}^{3}.$ The action of $SE(3)$ on $\mathbb{R}^{3}$ is the rotation $A$
followed by translation by the vector $a$ and has the expression
\begin{equation*}
(A,a)\cdot x=Ax+a.
\end{equation*}
Using this formula, one sees that multiplication and inversion in $SE(3)$
are given by
\begin{equation*}
(A,a)(B,b)=(AB,Ab+a)\qquad \text{and\qquad }(A,a)^{-1}=(A^{-1},-A^{-1}a),
\end{equation*}
for $A,B\in SO(3)$ and $a,b\in \mathbb{R}^{3}.$ The identity element is $%
(l,0)$.

The Lie algebra of the Euclidean group $SE(3)$ is $\mathfrak{se}(3)=\mathbb{R%
}^{3}\times \mathbb{R}^{3}$ with the Lie bracket
\begin{equation}
\lbrack (\xi ,u),(\eta ,v)]=(\xi \times \eta ,\xi \times v-\eta \times u).
\label{lbse3}
\end{equation}

The Lie algebra of the Euclidean group has a structure that is a special
case of what is called a \emph{semidirect product}. Here it is the \emph{%
product of the group of rotations with the corresponding group of
translations}. It turns out that semidirect products occur under rather
general circumstances when the symmetry in $T^{\ast }G$ is broken (see \cite%
{Marsden,GaneshSprBig,GaneshADG}).

The dual Lie algebra of the Euclidean group $SE(3)$ is $\mathfrak{se}%
(3)^{\ast }=\mathbb{R}^{3}\times \mathbb{R}^{3}$ with the same Lie bracket (%
\ref{lbse3}).

\paragraph{Representation of $SE(3)$}

In other words, $SE(3):=SO(3)\rhd \mathbb{R}^{3}$ is the Lie group
consisting of isometries of $\mathbb{R}^{3}$.

Using homogeneous coordinates, we can represent $SE(3)$ as follows,
\begin{equation*}
SE(3)=\ \ \left\{ \left(
\begin{array}{cc}
R & p \\
0 & 1%
\end{array}%
\right) \in GL(4,\mathbb{R}):R\in SO(3),\,p\in \mathbb{R}^{3}\right\} ,
\end{equation*}%
with the action on $\mathbb{R}^{3}$ given by the usual matrix--vector
product when we identify $\mathbb{R}^{3}$ with the section $\mathbb{R}%
^{3}\times \{1\}\subset \mathbb{R}^{4}$. In particular, given%
\begin{equation*}
g=\left(
\begin{array}{cc}
R & p \\
0 & 1%
\end{array}%
\right) \in SE(3),
\end{equation*}%
and $q\in \mathbb{R}^{3}$, we have%
\begin{equation*}
g\cdot q=Rq+p,
\end{equation*}%
or as a matrix--vector product,%
\begin{equation*}
\left(
\begin{array}{cc}
R & p \\
0 & 1%
\end{array}%
\right) \left(
\begin{array}{c}
q \\
1%
\end{array}%
\right) =\left(
\begin{array}{c}
Rq+p \\
1%
\end{array}%
\right) .
\end{equation*}

\paragraph{Lie algebra of $SE(3)$}

The Lie algebra of $SE(3)$ is given by \
\begin{equation*}
\mathfrak{se}(3)=\ \ \left\{ \left(
\begin{array}{cc}
\mathbf{\omega} & v \\
0 & 0%
\end{array}%
\right) \in M_{4}(\mathbb{R}):\mathbf{\omega}\in \mathfrak{so}(3),\,v\in
\mathbb{R}^{3}\right\} ,
\end{equation*}%
where the attitude matrix $\mathbf{\omega}:\mathbb{R}^{3}\rightarrow
\mathfrak{so}(3)$ is given by%
\begin{equation*}
\mathbf{\omega}=\left(
\begin{array}{ccc}
0 & -\omega _{z} & \omega _{y} \\
\omega _{z} & 0 & -\omega _{x} \\
-\omega _{y} & \omega _{x} & 0%
\end{array}%
\right) .
\end{equation*}

\paragraph{The exponential map of $SE(3)$}

The exponential map, $\exp :\mathfrak{se}(3)\rightarrow SE(3)$, is given by%
\begin{equation*}
\exp \left(
\begin{array}{cc}
\mathbf{\omega} & v \\
0 & 0%
\end{array}%
\right) =\left(
\begin{array}{cc}
\exp (\mathbf{\omega}) & Av \\
0 & 1%
\end{array}%
\right) ,
\end{equation*}%
where

\begin{equation*}
A=I+\frac{1-\cos \left\Vert \omega \right\Vert }{\left\Vert \omega
\right\Vert ^{2}}\mathbf{\omega}+\frac{\left\Vert \omega \right\Vert -\sin
\left\Vert \omega \right\Vert }{\left\Vert \omega \right\Vert ^{3}}\mathbf{%
\omega}^{2},
\end{equation*}%
and $\exp (\mathbf{\omega})$ is given by the Rodriguez' formula,%
\begin{equation*}
\exp (\mathbf{\omega})=I+\frac{\sin \left\Vert \omega \right\Vert }{%
\left\Vert \omega \right\Vert }\mathbf{\omega}+\frac{1-\cos \left\Vert
\omega \right\Vert }{\left\Vert \omega \right\Vert ^{2}}\mathbf{\omega}^{2}.
\end{equation*}

In other words, the special Euclidean group $SE(3):=SO(3)\rhd \Bbb{R}^{3}$ 
is the Lie group consisting of isometries of the Euclidean 3D
space $\Bbb{R}^{3}$. 
An element of $SE(3)$ is a pair $(A,a)$ where $A\in SO(3)$ and $a\in \Bbb{R}%
^{3}.$ The action of $SE(3)$ on $\Bbb{R}^{3}$ is the rotation $A$ followed
by translation by the vector $a$ and has the expression
\[
(A,a)\cdot x=Ax+a.
\]

The Lie algebra of the Euclidean group $SE(3)$ is $\mathfrak{se}(3)=\Bbb{R}%
^{3}\times \Bbb{R}^{3}$ with the Lie bracket
\begin{equation}
\lbrack (\xi ,u),(\eta ,v)]=(\xi \times \eta ,\xi \times v-\eta \times u).
\label{lbse3}
\end{equation}

Using homogeneous coordinates, we can represent $SE(3)$ as follows,
\[
SE(3)=\ \ \left\{ \left(
\begin{array}{cc}
R & p \\
0 & 1
\end{array}
\right) \in GL(4,\Bbb{R}):R\in SO(3),\,p\in \Bbb{R}^{3}\right\} ,
\]
with the action on $\Bbb{R}^{3}$ given by the usual matrix--vector product
when we identify $\Bbb{R}^{3}$ with the section $\Bbb{R}^{3}\times
\{1\}\subset \Bbb{R}^{4}$. In particular, given
\[
g=\left(
\begin{array}{cc}
R & p \\
0 & 1
\end{array}
\right) \in SE(3),
\]
and $q\in \Bbb{R}^{3}$, we have
\[
g\cdot q=Rq+p,
\]
or as a matrix--vector product,
\[
\left(
\begin{array}{cc}
R & p \\
0 & 1
\end{array}
\right) \left(
\begin{array}{c}
q \\
1
\end{array}
\right) =\left(
\begin{array}{c}
Rq+p \\
1
\end{array}
\right) .
\]

The Lie algebra of $SE(3)$, denoted $\mathfrak{se}(3)$, is given by \
\[
\mathfrak{se}(3)=\ \ \left\{ \left(
\begin{array}{cc}
\omega & v \\
0 & 0
\end{array}
\right) \in M_{4}(\Bbb{R}):\omega\in \mathfrak{so}(3),\,v\in \Bbb{R}%
^{3}\right\} ,
\]
where the attitude (or, angular velocity) matrix $\omega:\Bbb{R}%
^{3}\rightarrow \mathfrak{so}(3)$ is given by
\[
\omega=\left(
\begin{array}{ccc}
0 & -\omega _{z} & \omega _{y} \\
\omega _{z} & 0 & -\omega _{x} \\
-\omega _{y} & \omega _{x} & 0
\end{array}
\right) .
\]

The \emph{exponential map}, $\exp :\mathfrak{se}(3)\rightarrow
SE(3)$, is given by
\[
\exp \left(
\begin{array}{cc}
\omega & v \\
0 & 0
\end{array}
\right) =\left(
\begin{array}{cc}
\exp (\omega) & Av \\
0 & 1
\end{array}
\right) ,
\]
where

\[
A=I+\frac{1-\cos \left\Vert \omega \right\Vert }{\left\Vert \omega
\right\Vert ^{2}}\omega+\frac{\left\Vert \omega \right\Vert -\sin
\left\Vert \omega \right\Vert }{\left\Vert \omega \right\Vert
^{3}} \omega^{2},
\]
and $\exp (\omega)$ is given by the \emph{Rodriguez' formula},
\[
\exp (\omega)=I+\frac{\sin \left\Vert \omega \right\Vert }{%
\left\Vert \omega \right\Vert }\omega+\frac{1-\cos \left\Vert
\omega \right\Vert }{\left\Vert \omega \right\Vert
^{2}}\omega^{2}.
\]

\subsection{Basic Mechanical Examples}

\subsubsection{$SE(2)-$Hovercraft}

Configuration manifold is $(\theta,x,y)\in SE(2)$, given by matrix
\begin{equation*}
P=\left[
\begin{array}{ccc}
\cos \theta & \sin \theta & x \\
-\sin \theta & \cos \theta & y \\
0 & 0 & 1%
\end{array}
\right] .
\end{equation*}

Kinematic equations of motion in Lie algebra $\mathfrak{se}(2)$:
\begin{equation*}
\dot{P}=P\left[
\begin{array}{ccc}
0 & \omega & v_{x} \\
-\omega & 0 & v_{y} \\
0 & 0 & 0%
\end{array}
\right] ,\qquad (\omega =\dot{\theta},\,v_{x}=\dot{x},v_{y}=\dot{y}).
\end{equation*}

$\bigskip $Kinetic energy:
\begin{equation*}
E_{k}=\frac{1}{2}m(v_{x}^{2}+v_{y}^{2})+\frac{1}{2}I\omega ^{2},
\end{equation*}
where $m,I$ are mass and inertia moment of the hovercraft.

Dynamical equations of motion:
\begin{eqnarray*}
m\dot{v}_{x} &=&m\omega v_{y}+u_{1}, \\
m\dot{v}_{y} &=&-m\omega v_{x}+u_{2}, \\
I\dot{\omega} &=&\tau u_{2},
\end{eqnarray*}
where $\tau =-h$ is the torque applied at distance $h$ from the
center--of--mass, while $u_{1},u_{2}$ are control inputs.

\subsubsection{$SO(3)-$Satellite}

Configuration manifold is rotation matrix $R\in SO(3)$, with associated
angular--velocity (attitude) matrix $\mathbf{\omega}=(\omega_1,\omega_2,%
\omega_3)\in \mathfrak{so}(3)\approx \mathbb{R}^{3}$ given by
\begin{equation*}
\mathbf{\omega}\in \mathfrak{so}(3)\longmapsto \left[
\begin{array}{ccc}
0 & -\omega _{3} & \omega _{2} \\
\omega _{3} & 0 & -\omega _{1} \\
-\omega _{2} & \omega _{1} & 0%
\end{array}
\right] .
\end{equation*}

Kinematic equation of motion in $\mathfrak{so}(3)$:
\begin{equation*}
\dot{R}=R\mathbf{\omega},
\end{equation*}

Kinetic energy:
\begin{equation*}
E_{k}=\frac{1}{2}\mathbf{\omega }^{T}\mathbf{I\omega },
\end{equation*}
where inertia tensor $\mathbf{I}$ is given by diagonal matrix,
\begin{equation*}
\mathbf{I}=diag\{I_{1},I_{2},I_{3}\}.
\end{equation*}

Dynamical Euler equations of motion:
\begin{equation*}
\mathbf{I\dot{\omega}}=\mathbf{I\omega \times \omega } + \tau_iu^i,
\end{equation*}
where $\times$ is the cross--product in 3D, $\tau_i$ are three external
torques and $u^i=u^i(t)$ are control inputs.

\subsubsection{$SE(3)-$Submarine}

The motion of a rigid body in incompressible, irrotational and inviscid
fluid is defined by the configuration manifold $SE(3)$, given by a pair of
rotation matrix and translation vector, $(R,p)\in SE(3)$, such that angular
velocity (attitude) matrix and linear velocity vector, $(\mathbf{\omega },%
\mathbf{v})\in\mathfrak{se}(3)\approx \mathbb{R}^{6}$.\bigskip

Kinematic equations of motion in $\mathfrak{se}(3)$:
\begin{equation*}
\dot{p}=R\mathbf{v},\qquad \dot{R}=R\mathbf{\omega}.
\end{equation*}

Kinetic energy (symmetrical):
\begin{equation*}
E_{k}=\frac{1}{2}\mathbf{v}^{T}\mathbf{Mv}+\frac{1}{2}\mathbf{\omega }^{T}%
\mathbf{I\omega },
\end{equation*}
where mass and inertia matrices are diagonal (for a neutrally buoyant
ellipsoidal body with uniformly distributed mass),
\begin{eqnarray*}
\mathbf{M} &=&diag\{m_{1},m_{2},m_{3}\}, \\
\mathbf{I} &=&diag\{I_{1},I_{2},I_{3}\}.
\end{eqnarray*}

Dynamical \textit{Kirchhoff equations of motion} read:
\begin{equation*}
\mathbf{M\dot{v}=Mv\times \omega },\qquad \mathbf{I\dot{\omega}}=\mathbf{%
I\omega \times \omega }+\mathbf{Mv\times v}.
\end{equation*}
By including the body--fixed external forces and torques, $f_{i},\tau _{i}$,
with input controls $u^{i}=u^{(}t)$, the dynamical equations become:
\begin{eqnarray*}
\mathbf{M\dot{v}} &=&\mathbf{Mv\times \omega }+f_{i}u^{i}, \\
\mathbf{I\dot{\omega}} &=&\mathbf{I\omega \times \omega }+\mathbf{Mv\times v}%
+\tau _{i}u^{i}.
\end{eqnarray*}

\subsection{Newton--Euler $SE(3)-$Dynamics}

\subsubsection{$SO(3):$ Euler Equations of Rigid Rotations}

Unforced Euler equations read in vector form
\begin{equation*}
\mathbf{\dot{\,\pi }}\equiv \mathbf{I\dot{\omega}}=\mathbf{\pi \times \omega
},\qquad \text{with}\quad\mathbf{I}=diag\{I_{1},I_{2},I_{3}\}
\end{equation*}
and in scalar form
\begin{equation*}
\begin{matrix}
I_{1}\dot{\omega}_{1}=(I_{2}-I_{3})\omega _{2}\omega _{3} \\
I_{2}\dot{\omega}_{2}=(I_{3}-I_{1})\omega _{3}\omega _{1} \\
I_{3}\dot{\omega}_{3}=(I_{1}-I_{2})\omega _{1}\omega _{2}%
\end{matrix}%
.
\end{equation*}
Using rotational kinetic--energy Lagrangian
\begin{equation*}
L(\mathbf{\omega })=E_{k}^{rot}={\frac{1}{2}}\mathbf{\omega }^{t}\mathbf{%
I\omega }=\frac{1}{2}(I_{1}\omega _{1}^{2}+I_{2}\omega _{2}^{2}+I_{3}\omega
_{3}^{2})\qquad (^t=\text{`transpose'})
\end{equation*}
Regarding the angular momentum ~$\mathbf{\pi }=\partial _{\mathbf{\omega }}L=%
\mathbf{I\omega }$ = $(I_{1}\omega _{1},I_{2}\omega _{2},I_{3}\omega _{3})$
as a vector, we can derive unforced Euler equations: $\mathbf{\dot{\,\pi }}=%
\mathbf{\pi \times \omega }$ as a system of \textit{%
Euler--Lagrange--Kirchhoff equations}
\begin{equation*}
\frac{d}{dt}\partial _{\mathbf{\omega }}L=\partial _{\mathbf{\omega }%
}L\times \mathbf{\omega }.
\end{equation*}

Forced Euler equations read in vector form
\begin{equation*}
\mathbf{\dot{\,\pi }}+\mathbf{\omega }\times \mathbf{\pi }=\mathbf{T}
\end{equation*}
and in scalar form
\begin{equation*}
\begin{matrix}
I_{1}\dot{\omega}_{1}+(I_{3}-I_{2})\omega _{2}\omega _{3} & = & T_{1}\, \\
I_{2}\dot{\omega}_{2}+(I_{1}-I_{3})\omega _{3}\omega _{1} & = & T_{2} \\
I_{3}\dot{\omega}_{3}+(I_{2}-I_{1})\omega _{1}\omega _{2} & = & T_{3}%
\end{matrix}%
\end{equation*}

\subsubsection{$SE(3):$ Coupled Newton--Euler Equations}

Forced coupled Newton--Euler equations read in vector form
\begin{eqnarray*}
&&\mathbf{\dot{p}}~ \mathbf{\equiv M\dot{v}=F+p\times \omega },\qquad \text{%
with}~~\mathbf{M}=diag\{m_{1},m_{2},m_{3}\} \\
&&\mathbf{\dot{\pi}}~ \mathbf{\equiv I\dot{\omega}=T+\pi \times \omega
+p\times v},\qquad \mathbf{I}=diag\{I_{1},I_{2},I_{3}\},
\end{eqnarray*}
with principal inertia moments given in Cartesian coordinates ($x,y,z$) by
density $\rho-$dependent volume integrals
\begin{equation*}
I_{1}=\iiint \rho (z^{2}+y^{2})dxdydz, ~~ I_{2}=\iiint \rho
(x^{2}+y^{2})dxdydz, ~~ I_{3}=\iiint \rho (x^{2}+y^{2})dxdydz,
\end{equation*}

In tensor form, the forced--coupled Newton--Euler equations read
\begin{eqnarray*}
\dot{p}_{i} &\equiv &M_{ij}\dot{v}^{j}=F_{i}+\varepsilon _{ik}^{j}p_{j}{%
\omega }^{k}, \\
\dot{\pi}_{i} &\equiv &I_{ij}\dot{\omega}^{j}=T_{i}+\varepsilon _{ik}^{j}\pi
_{j}\omega ^{k}+\varepsilon _{ik}^{j}p_{j}v^{k},
\end{eqnarray*}
where the permutation symbol $\varepsilon _{ik}^{j}$ is defined as
\begin{equation*}
\varepsilon _{ik}^{j}=
\begin{cases}
+1 & \text{if }(i,j,k)\text{ is }(1,2,3),(3,1,2)\text{ or }(2,3,1), \\
-1 & \text{if }(i,j,k)\text{ is }(3,2,1),(1,3,2)\text{ or }(2,1,3), \\
0 & \text{otherwise: }i=j\text{ or }j=k\text{ or }k=i.%
\end{cases}%
\end{equation*}

In scalar form these equations read
\begin{equation*}
\begin{array}{c}
\dot{p}_{_1}={F_1}-{m_3}{v_3}{\omega _2}+{m_2}{v_2}{\omega _3} \\
\dot{p}_{_2}={F_2}+{m_3}{v_3}{\omega _1}-{m_1}{v_1}{\omega _3} \\
\dot{p}_{_3}={F_3}-{m_2}{v_2}{\omega _1}+{m_1}{v_1}{\omega _2} \\
\\
\dot{\pi}_{_1}={T_1}+({m_2}-{m_3}){v_2}{v_3}+({I_2}-{I_3}){\omega _2}{\omega
_3} \\
\dot{\pi}_{_2}={T_2}+({m_3}-{m_1}){v_1}{v_3}+({I_3}-{I_1}){\omega _1}{\omega
_3} \\
\dot{\pi}_{_3}={T_3}+({m_1}-{m_2}){v_1}{v_2}+({I_1}-{I_2}){\omega _1}{\omega
_2}%
\end{array}%
.
\end{equation*}
These coupled rigid--body equations can be derived from the \textit{%
Newton--Euler kinetic energy}
\begin{equation*}
E_{k}={\frac{1}{2}}\mathbf{v}^{t}\mathbf{Mv}+{\frac{1}{2}}\mathbf{\omega }%
^{t}\mathbf{I\omega }
\end{equation*}
or, in tensor form
\begin{equation*}
E={\frac{1}{2}}M_{ij}\dot{v}^{i}\dot{v}^{j}+{\frac{1}{2}}I_{ij}\dot{\omega}%
^{i}\dot{\omega}^{j}.
\end{equation*}
Using the \emph{Kirchhoff--Lagrangian equations}
\begin{eqnarray*}
\frac{d}{{dt}}\partial _{\mathbf{v}}E_{k} &=&\partial _{\mathbf{v}%
}E_{k}\times \mathbf{\omega }+\mathbf{F} \\
{\frac{d}{{dt}}}\partial _{\mathbf{\omega }}E_{k} &=&\partial _{\mathbf{%
\omega }}E_{k}\times \mathbf{\omega }+\partial _{\mathbf{v}}E_{k}\times
\mathbf{v}+\mathbf{T},
\end{eqnarray*}
or, in tensor form
\begin{eqnarray*}
\frac{d}{dt}\partial _{v^{i}}E &=&\varepsilon _{ik}^{j}\left( \partial
_{v^{j}}E\right) \omega ^{k}+F_{i}, \\
\frac{d}{dt}\partial _{{\omega }^{i}}E &=&\varepsilon _{ik}^{j}\left(
\partial _{{\omega }^{j}}E\right) {\omega }^{k}+\varepsilon _{ik}^{j}\left(
\partial _{v^{j}}E\right) v^{k}+T_{i}
\end{eqnarray*}
we can derive linear and angular momentum covectors
\begin{equation*}
\mathbf{p}=\partial _{\mathbf{v}}E_{k}{,\qquad \mathbf{\pi }=\partial _{%
\mathbf{\omega }}E_{k}}
\end{equation*}
or, in tensor form
\begin{equation*}
p_{i}=\partial _{v^{i}}E{,\qquad }\pi _{i}=\partial _{{\omega }^{i}}E,
\end{equation*}
and in scalar form
\begin{eqnarray*}
\mathbf{p} &=&[p_{1},p_{2},p_{3}]=[m_{1}v_{1},m_{2}v_{2},m_{2}v_{3}] \\
\mathbf{\pi } &=&[\pi _{1},\pi _{2},\pi _{3}]=[I_{1}\omega _{1},I_{2}\omega
_{2},I_{3}\omega _{3}],
\end{eqnarray*}
with their respective time derivatives, in vector form
\begin{equation*}
~\mathbf{\dot{p}}=\frac{d}{dt}\mathbf{p=}\frac{d}{{dt}}\partial _{\mathbf{v}%
}E_{k}{,\qquad \mathbf{\dot{\pi}}=}\frac{d}{dt}\mathbf{\pi =}{\frac{d}{{dt}}}%
\partial _{\mathbf{\omega }}E_{k}
\end{equation*}
or, in tensor form
\begin{equation*}
~\dot{p}_{i}=\frac{d}{dt}p_{i}=\frac{d}{dt}\partial _{v^{i}}E{,\qquad \dot{%
\pi}_{i}=}\frac{d}{dt}\pi _{i}=\frac{d}{dt}\partial _{{\omega }^{i}}E,
\end{equation*}
and in scalar form
\begin{eqnarray*}
\mathbf{\dot{p}} &=&[\dot{p}_{1},\dot{p}_{2},\dot{p}_{3}]=[m_{1}\dot{v}%
_{1},m_{2}\dot{v}_{2},m_{3}\dot{v}_{3}] \\
{\mathbf{\dot{\pi}}} &=&[\dot{\pi}_{1},\dot{\pi}_{2},\dot{\pi}_{3}]=[I_{1}%
\dot{\omega}_{1},I_{2}\dot{\omega}_{2},I_{3}\dot{\omega}_{3}].
\end{eqnarray*}

In addition, for the purpose of biomechanical injury prediction/prevention,
we have linear and angular jolts, respectively given in vector form by
\begin{eqnarray*}
~&&\mathbf{\dot{F}=\ddot{p}-\dot{p}\times \omega -p\times \dot{\omega}}\qquad
\\
~&&\mathbf{\dot{T}=\ddot{\pi}}~\mathbf{-\dot{\pi}\times \omega -\pi \times
\dot{\omega}-\dot{p}\times v-p\times \dot{v}},
\end{eqnarray*}
or, in tensor form\footnote{%
In this paragraph the overdots actually denote the absolute Bianchi
(covariant) derivatives, so that the jolts retain the proper covector
character, which would be lost if ordinary time derivatives are used.
However, for simplicity, we stick to the same notation.}
\begin{eqnarray*}
~\dot{F}_{i} &=&\ddot{p}_{i}-\varepsilon _{ik}^{j}\dot{p}_{j}{\omega }%
^{k}-\varepsilon _{ik}^{j}p_{j}{\dot{\omega}}^{k}, \\
~\dot{T}_{{i}} &=&\ddot{\pi}_{i}~-\varepsilon _{ik}^{j}\dot{\pi}_{j}\omega
^{k}-\varepsilon _{ik}^{j}\pi _{j}{\dot{\omega}}^{k}-\varepsilon _{ik}^{j}%
\dot{p}_{j}v^{k}-\varepsilon _{ik}^{j}p_{j}\dot{v}^{k},
\end{eqnarray*}
where the linear and angular jolt covectors are
\begin{eqnarray*}
\mathbf{\dot{F}} &\equiv &\dot{F}_{i}=\mathbf{M\ddot{v}}\,\equiv \mathbf{\,}%
M_{ij}\ddot{v}^{j}=[\dot{F}_{1},\dot{F}_{2},\dot{F}_{3}], \\
\mathbf{\dot{T}} &\equiv &\dot{T}_{{i}}=\mathbf{I\ddot{\omega}\equiv \,}%
I_{ij}\ddot{\omega}^{j}=[\dot{T}_{{1}},\dot{T}_{{2}},\dot{T}_{{3}}],
\end{eqnarray*}
where
\begin{equation*}
\mathbf{\ddot{v}}=\ddot{v}^{{i}}=[\ddot{v}^{{1}},\ddot{v}^{{2}},\ddot{v}^{{3}%
}]^{t},\qquad \mathbf{\ddot{\omega}}=\ddot{\omega}^{{i}}=[\ddot{\omega}^{{1}%
},\ddot{\omega}^{{2}},\ddot{\omega}^{{3}}]^{t},
\end{equation*}
are linear and angular jerk vectors.

In scalar form, the $SE(3)-$jolt expands as
\begin{eqnarray*}
&&\left\{
\begin{array}{l}
\dot{F}_{{1}}=\ddot{p}_{1}-m_{{2}}\omega _{{3}}\dot{v}_{{2}}+m_{{3}}\left( {%
\omega }_{{2}}\dot{v}_{{3}}+v_{{3}}\dot{\omega}_{{2}}\right) -m_{{2}}v_{{2}}{%
\dot{\omega}}_{{3}}, \\
\dot{F}_{{2}}=\ddot{p}_{2}+m_{{1}}\omega _{{3}}\dot{v}_{{1}}-m_{{3}}\omega _{%
{1}}\dot{v}_{{3}}-m_{{3}}v_{{3}}\dot{\omega}_{{1}}+m_{{1}}v_{{1}}\dot{\omega}%
_{{3}}, \\
\dot{F}_{{3}}=\ddot{p}_{3}-m_{{1}}\omega _{{2}}\dot{v}_{{1}}+m_{{2}}\omega _{%
{1}}\dot{v}_{{2}}-v_{{2}}\dot{\omega}_{{1}}-m_{{1}}v_{{1}}\dot{\omega}_{{2}},%
\end{array}
\right. \\
&& \\
&&\left\{
\begin{array}{l}
\dot{T}_{{1}}=\ddot{\pi}_{1}-(m_{{2}}-m_{{3}})\left( v_{{3}}\dot{v}_{{2}}+v_{%
{2}}\dot{v}_{{3}}\right) -(I_{{2}}-I_{{3}})\left( \omega _{{3}}\dot{\omega}_{%
{2}}+{\omega }_{{2}}{\dot{\omega}}_{{3}}\right) , \\
\dot{T}_{{2}}=\ddot{\pi}_{2}+(m_{{1}}-m_{{3}})\left( v_{{3}}\dot{v}_{{1}}+v_{%
{1}}\dot{v}_{{3}}\right) +(I_{{1}}-I_{{3}})\left( {\omega }_{{3}}{\dot{\omega%
}}_{{1}}+{\omega }_{{1}}{\dot{\omega}}_{{3}}\right) , \\
\dot{T}_{{3}}=\ddot{\pi}_{3}-(m_{{1}}-m_{{2}})\left( v_{{2}}\dot{v}_{{1}}+v_{%
{1}}\dot{v}_{{2}}\right) -(I_{{1}}-I_{{2}})\left( {\omega }_{{2}}{\dot{\omega%
}}_{{1}}+{\omega }_{{1}}{\dot{\omega}}_{{2}}\right).%
\end{array}
\right.
\end{eqnarray*}

\subsection{Symplectic Group in Hamiltonian Mechanics}

\label{Sp}

Here we give a brief description of symplectic group (see \cite%
{Marsden,GaneshSprBig,GaneshADG}).

Let $J=\left(
\begin{array}{cc}
0 & I \\
-I & 0%
\end{array}
\right) , $ with $I$ the $n\times n$ identity matrix. Now, $A\in L(\mathbb{R}%
^{2n}, \mathbb{R}^{2n})$ is called a \emph{symplectic matrix} if $A^{T}J%
\mathbf{\,}A=J$. Let $Sp(2n,\mathbb{R})$ be the set of $2n\times 2n$
symplectic matrices. Taking determinants of the condition $A^{T}J\mathbf{\,}%
A=J$ gives $\det A=\pm 1$, and so $A\in GL(2n,\mathbb{R})$. Furthermore, if $%
A,B\in Sp(2n,R)$, then $(AB)^{T}J(AB)=B^{T}A^{T}JAB=J$. Hence, $AB\in Sp(2n,%
\mathbb{R})$, and if $A^{T}J\mathbf{\,}A=J$, then $%
JA=(A^{T})^{-1}J=(A^{-1})^{T}J$, so $J=\left( A-1\right) ^{T}JA^{-1}$, or $%
A^{-1}\in Sp(2n,\mathbb{R})$. Thus, $Sp(2n,\mathbb{R}) $ is a group.

The \emph{symplectic Lie group}
\begin{equation*}
Sp(2n,\mathbb{R})=\left\{ A\in GL(2n,\mathbb{R}):A^{T}J\mathbf{\,}A=J\right\}
\end{equation*}
is a noncompact, connected Lie group of dimension $2n^{2}+n$. Its Lie
algebra
\begin{equation*}
\mathfrak{sp}(2n,\mathbb{R})=\left\{ A\in L(\mathbb{R}^{2n},\mathbb{R}%
^{2n}):A^{T}J\mathbf{\,}A=J=0\right\} ,
\end{equation*}
called the \emph{symplectic Lie algebra}, consists of the $2n\times 2n$
matrices $A$ satisfying $A^{T}J\mathbf{\,}A=0$.

Consider a particle of mass $m$ moving in a potential $V(q)$, where $%
q^{i}=(q^{1},q^{2},q^{3})\in \mathbb{R}^{3}$. Newtonian second law states
that the particle moves along a curve $q(t)$ in $\mathbb{R}^{3}$ in such a
way that $m\ddot{q}^{i}=-\limfunc{grad}V(q^{i})$. Introduce the 3D--momentum
$p_{i}=m\dot{q}^{i}$, and the energy (Hamiltonian)
\begin{equation*}
H(q,p)=\frac{1}{2m}\sum_{i=1}^{3}p_{i}^{2}+V(q).
\end{equation*}
Then
\begin{eqnarray*}
\frac{\partial H}{\partial q^{i}} &=&\frac{\partial V}{\partial q^{i}}=-m%
\ddot{q}^{i}=-\dot{p}_{i},\text{ \ \ and} \\
\frac{\partial H}{\partial p_{i}} &=&\frac{1}{m}p_{i}=\dot{q}^{i},\qquad
(i=1,2,3),
\end{eqnarray*}

and hence Newtonian law \fbox{$F=m\ddot{q}^{i}$} is equivalent to \textit{%
Hamiltonian equations}
\begin{equation*}
\dot{q}^{i}=\frac{\partial H}{\partial p_{i}},\qquad \dot{p}_{i}=-\frac{%
\partial H}{\partial q^{i}}.
\end{equation*}

Now, writing $z=(q^{i},p_{i})$,
\begin{equation*}
J\limfunc{grad}H(z)=\left(
\begin{array}{cc}
0 & I \\
-I & 0%
\end{array}
\right) \left(
\begin{array}{c}
\frac{\partial H}{\partial q^{i}} \\
\frac{\partial H}{\partial p_{i}}%
\end{array}
\right) =\left( \dot{q}^{i},\dot{p}_{i}\right) =\dot{z},
\end{equation*}
so the \textit{complex Hamiltonian equations} read
\begin{equation*}
\dot{z}=J\limfunc{grad}H(z).
\end{equation*}
Now let $f:\mathbb{R}^{3}\times \mathbb{R}^{3}\rightarrow \mathbb{R}%
^{3}\times \mathbb{R}^{3}$ and write $w=f(z)$. If $z(t)$ satisfies the
complex Hamiltonian equations then $w(t)=f(z(t))$ satisfies $\dot{w}=A^{T}%
\dot{z}$, where $A^{T}=[\partial w^{i}/\partial z^{j}]$ is the Jacobian
matrix of $f$. By the chain rule,

\begin{equation*}
\dot{w}=A^{T}J\limfunc{grad}_{z}H(z)=A^{T}J\mathbf{\,}A\limfunc{grad}%
_{w}H(z(w)).
\end{equation*}%
Thus, the equations for $w(t)$ have the form of Hamiltonian equations with
energy $K(w)=H(z(w))$ iff $A^{T}J\mathbf{\,}A=J$, that is, iff $A$ is
symplectic. A nonlinear transformation $f$ is canonical iff its Jacobian
matrix is symplectic. $Sp(2n,\mathbb{R})$ is the linear invariance group of
classical mechanics.

\section{Medical Applications: Prediction of Injuries}

\subsection{General Theory of Musculo--Skeletal Injury Mechanics}

The prediction and prevention of traumatic brain injury, spinal
injury and musculo-skeletal injury is a very important aspect of preventive
medical science. Recently, in a series of papers \cite{ivbrain,ivspine,ivgen}, we have proposed a new coupled loading-rate hypothesis as a unique cause of
all above injuries. This new hypothesis states that the unique cause of brain, spinal
and musculo-skeletal injuries is a Euclidean Jolt, which is an impulsive loading that
strikes any part of the human body (head, spine or any bone/joint) -- in
several coupled degrees-of-freedom simultaneously. It never goes in a single
direction only. Also, it is never a static force. It is always an impulsive
translational and/or rotational force coupled to some mass eccentricity.
This is, in a nutshell, our universal Jolt theory of all mechanical injuries.%
\begin{figure}[h]
\centerline{\includegraphics[width=13cm]{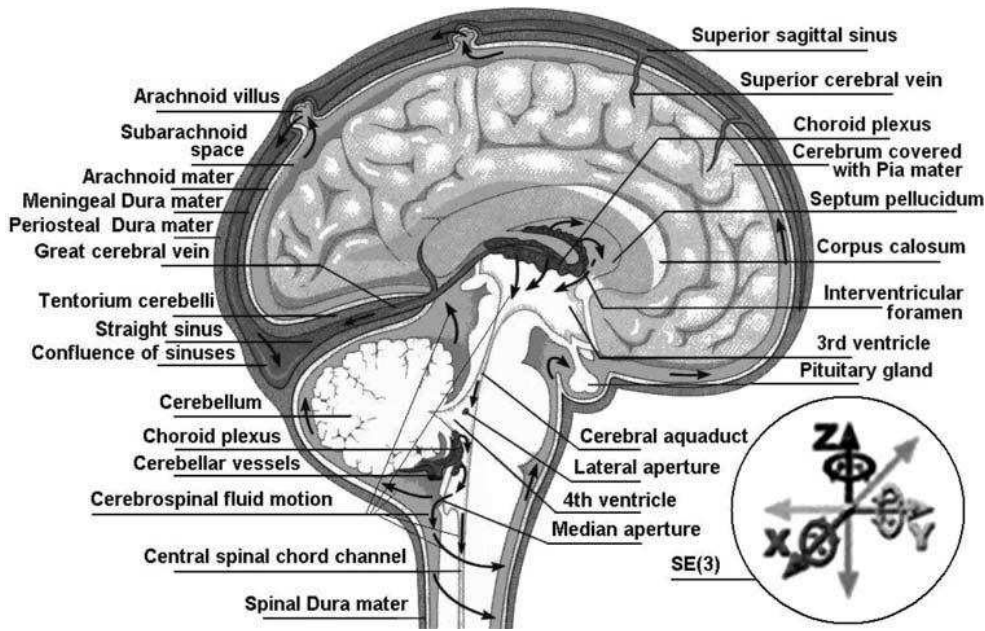}}
\caption{Human brain and its
$SE(3)-$group of microscopic
three-dimensional motions
within the cerebrospinal fluid
inside the cranial cavity.}
\label{BrainInj}
\end{figure}

To show this, based on the previously defined covariant force law, we have
firstly formulated the fully coupled Newton--Euler dynamics of:

1.\quad Brain's micro-motions within the cerebrospinal fluid inside the
cranial cavity;

2.\quad Any local inter-vertebral motions along the spine; and

3.\quad Any local joint motions in the human musculo-skeletal system.\newline

Then, from it, we have defined the essential concept of \textbf{Euclidean Jolt%
}, which is the main cause of all mechanical injuries. The Euclidean Jolt
has two main components:

1.\quad Sudden motion, caused either by an accidental impact or slightly
distorted human movement; and

2.\quad Unnatural mass distribution of the human body (possibly with some
added masses), which causes some mass eccentricity from the natural
physiological body state.
\begin{figure}[h]
\centerline{\includegraphics[width=13cm]{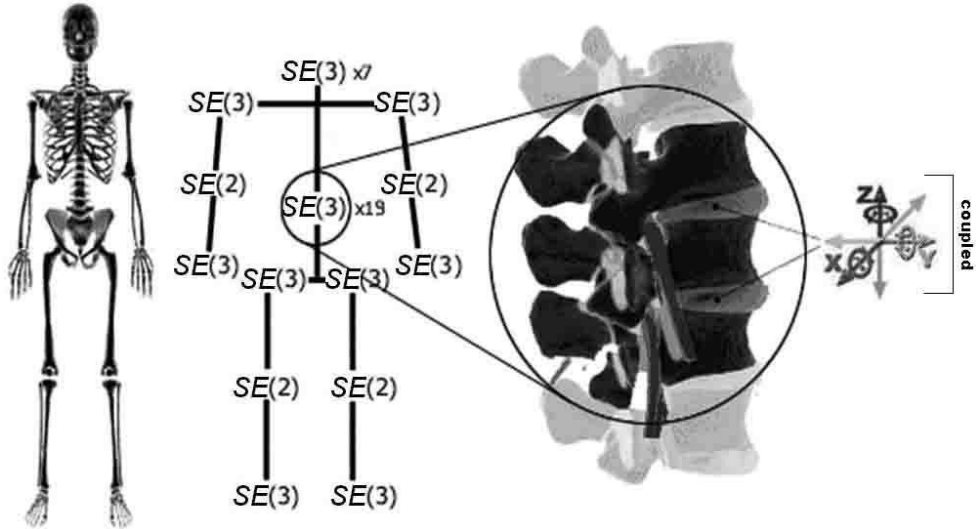}}
\caption{Human body representation in terms of SE(3)/SE(2)-groups of rigid-body motion, with
the vertebral column represented as a chain of 26 flexibly-coupled SE(3)-groups.}
\label{SpineSE3}
\end{figure}

What does this all mean? We will try to explain it in \textquotedblleft plain
English\textquotedblright. As we live in a 3D space, one could think that
motion of any part of the human body, either caused by an accidental impact
or by voluntary human movement, \textquotedblleft just obeys classical
mechanics in 6 degrees-of-freedom: three translations and three
rotations\textquotedblright . However, these 6 degrees-of-freedom are not
independent motions as it is suggested by the standard term
\textquotedblleft degrees-of-freedom\textquotedblright . In reality, these
six motions of any body in space are coupled. Firstly, three rotations are
coupled in the so-called rotation group (or matrix, or quaternion).
Secondly, three translations are coupled with the rotation group to give the
full Euclidean group of rigid body motions in space. A simple way to see
this is to observe someone throwing an object in the air or hitting a tennis
ball: how far and where it will fly depends not only on the standard
\textquotedblleft projectile\textquotedblright\ mechanics, but also on its
local \textquotedblleft spin\textquotedblright\ around all three axes
simultaneously. Every golf and tennis player knows this simple fact. Once
the spin is properly defined we have a \textquotedblleft fully coupled
Newton--Euler dynamics\textquotedblright\ -- to start with.\newline

The covariant force law for any biodynamical system (which we introduced
earlier in our biodynamics books and papers, see our references in the cited
papers above) goes one step beyond the Newton--Euler dynamics. It states:
\begin{eqnarray*}
\mathbf{Euclidean\ Force\ covector\ field\ } \qquad\mathbf{=\newline
} \\
\mathbf{Body\ mass\ distribution\ \times \ Euclidean\ Acceleration\
vector\ field}
\end{eqnarray*}

This is a nontrivial biomechanical generalization of the fundamental
Newton's definition of the force acting on a single particle. Unlike
classical engineering mechanics of multi-body systems, this fundamental law
of biomechanics proposes that forces acting on a multi-body system and
causing its motions are fundamentally different physical quantities from the
resulting accelerations. In simple words, forces are massive quantities
while accelerations are massless quantities. More precisely, the
acceleration vector field includes all linear and angular accelerations of
individual body segments. When we couple them all with the total body's
mass-distribution matrix of all body segments (including all masses and
inertia moments), we get the force co-vector field, comprising all the
forces and torques acting on the individual body segments. In this way, we
have defined the 6-dimensional Euclidean force for an arbitrary
biomechanical system.
\begin{figure}[h]
\centerline{\includegraphics[width=7cm]{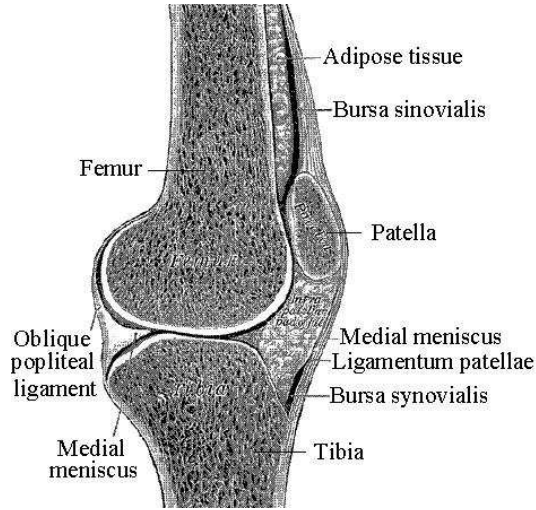}}
\caption{Schematic latero-frontal view of the left knee joint. Although designed to perform mainly
flexion/extension (strictly in the sagittal plane) with some restricted medial/lateral rotation in
the semi-flexed position, it is clear that the knee joint really has at least six-degrees-of-freedom,
including three micro-translations. The injury actually occurs when some of these microscopic
translations become macroscopic, which normally happens only after an external jolt.}
\label{KneeInj}
\end{figure}

Now, for prediction of injuries, we need to take the rate-of-change (or
derivative, with respect to time) of the Euclidean biomechanical force
defined above. In this way, we get the Euclidean Jolt, which is the sudden
change (in time) of the 6-dimensional Euclidean force:
\begin{eqnarray*}
\mathbf{Euclidean\ Jolt\ covector\ field\ } \qquad\mathbf{=\ } \\
\mathbf{Body\ mass\ distribution\ \times \ Euclidean\ Jerk\ vector\ field}
\end{eqnarray*}

And again, it consists of two components: (i) massless linear and angular
jerks (of all included body segments), and (ii) their mass distribution. For
the sake of simplicity, we can say that the mass distribution matrix
includes all involved segmental masses and inertia moments, as well as
\textquotedblleft eccentricities\textquotedblright\ or \textquotedblleft
pathological leverages\textquotedblright\ from the normal physiological
state.\newline

Therefore, the unique cause of all brain, spine and musculo-skeletal
injuries has two components:

1.\quad Coupled linear and angular jerks; ~and

2.\quad Mass distribution with \textquotedblleft
eccentricities\textquotedblright .\newline

\noindent In other words, ~\textbf{there are no injuries in static
conditions without any mass eccentricities; all injuries are caused by
mutually coupled linear and angular jerks, which are also coupled with the
involved human mass distribution.\newline
}

The Euclidean Jolt causes two forms of discontinuous brain, spine or
musculo-skeletal injury:

1.\quad Mild rotational disclinations; ~and

2.\quad Severe translational dislocations (or, fractures).\newline

In the cited papers above, we have developed the soft-body dynamics of
biomechanical disclinations and dislocations, caused by the Euclidean Jolt,
using the Cosserat multipolar viscoelastic continuum model.\\

Implications of the new universal theory are various, as follows.\newline

\noindent\textbf{A.} ~~The research in traumatic brain injury (TBI, see Figure \ref{BrainInj}) has so
far identified the rotation of the brain-stem as the main cause of the TBI
due to various crashes/impacts. The contribution of our universal Jolt theory
to the TBI research is the following:\newline

1.\quad Rigorously defined this brain rotation as a mechanical disclination
of the brain-stem tissue modelled by the Cosserat multipolar soft-body model;

2.\quad Showing that brain rotation is never uni-axial but always
three-axial;

3.\quad Showing that brain rotation is always coupled with translational
dislocations. This is a straightforward consequence of our universal Jolt
theory.\newline

These apparently `obvious' facts are actually \textsl{radically new:} we
cannot separately analyze rapid brain's rotations from translations, because
they are in reality always coupled.\newline

One practical application of the brain Jolt theory is in design of helmets.
Briefly, a `hard' helmet saves the skull but not the brain; alternatively, a
`soft' helmet protects the brain from the collision jolt but does not
protect the skull. A good helmet is both `hard' and `soft'. A proper helmet
would need to have both a hard external shell (to protect the skull) and a
soft internal part (that will dissipate the energy from the collision jolt
by its own destruction, in the same way as a car saves its passengers from
the collision jolt by its own destruction).\newline

Similarly, in designing safer car air-bags, the two critical points will be
(i) their placement within the car, and (ii) their \textquotedblleft
soft-hard characteristics\textquotedblright , similar to the helmet
characteristics described above.\newline

\noindent\textbf{B.} ~~In case of spinal injury (see Figure \ref{SpineSE3}), the contribution of our
universal Jolt theory is the following:\newline

1.\quad The spinal injury is always localized at the certain vertebral or
inter-vertebral point;

2.\quad In case of severe translational injuries (vertebral fractures or
discus herniae) they can be identified using X-ray or other medical imaging
scans; in case of microscopic rotational injuries (causing the back-pain
syndrome) they cannot be identified using current medical imaging scans;

3.\quad There is no spinal injury without one of the following two causes:

\qquad a.~~~Impulsive rotational + translational loading caused by either
fast human movements or various
crashes/impacts; and/or

\qquad b.~~~Static eccentricity from the normal physiological spinal form,
caused by external loading;

\qquad c.~~~Any spinal injury is caused by a combination of the two points
above: impulsive rotational + translational loading and static eccentricity.\newline

This is a straightforward consequence of our universal Jolt theory. We cannot
separately analyze translational and rotational spinal injuries. Also, there
are no \textquotedblleft static injuries\textquotedblright\ without
eccentricity. Indian women have for centuries carried bulky loads on their
heads without any spinal injuries; they just prevented any load
eccentricities and any jerks in their motion.\newline

The currently used \textquotedblleft Principal loading
hypothesis\textquotedblright\ that describes spinal injuries in terms of
spinal tension, compression, bending, and shear, covers only a small subset
of all spinal injuries covered by our universal Jolt theory. To prevent
spinal injuries we need to develop spinal jolt awareness: ability to control
all possible impulsive spinal loadings as well as static eccentricities.%
\newline

\noindent\textbf{C.} ~~In case of general musculo-skeletal injury (see Figure \ref{KneeInj} for the particular case of knee injury), the
contribution of our universal Jolt theory is the following:\newline

1.\quad The injury is always localized at the certain joint or bone and
caused by an impulsive loading, which hits this particular joint/bone in
several coupled degrees-of-freedom simultaneously;

2.\quad Injury happens when most of the body mass is hanging on that joint;
for example, in case of a knee injury, when most of the body mass is on one
leg with a semi-flexed knee --- and then, caused by some external shock, the
knee suddenly \textquotedblleft jerks\textquotedblright\ (this can happen in
running, skiing, and ball games, as well as various crashes/impacts); or, in
case of shoulder injury, when most of the body mass is hanging on one arm
and then it suddenly jerks.\newline

To prevent these injuries we need to develop musculo-skeletal jolt
awareness. For example, never overload a flexed knee and avoid any kind of
uncontrolled motions (like slipping) or collisions with external objects.

\subsection{Analytical Mechanics of Traumatic Brain Injury (TBI)}

\subsubsection{The $SE(3)-$jolt: the cause of TBI}

In this subsection we give a brief on TBI mechanics. For more details and references, see \cite{ivbrain}.

In the language of modern dynamics, the
microscopic motion of human brain within the skull is governed by
the Euclidean SE(3)--group of 3D motions (see next subsection).
Within brain's SE(3)--group we have both SE(3)--kinematics
(consisting of SE(3)--velocity and its two time derivatives:
SE(3)--acceleration and SE(3)--jerk) and SE(3)--dynamics
(consisting of SE(3)--momentum and its two time derivatives:
SE(3)--force and SE(3)--jolt), which is brain's kinematics
$\times $ brain's mass--inertia distribution.

Informally, the external SE(3)--jolt\footnote{The mechanical
SE(3)--jolt concept is based on the mathematical concept of
higher--order tangency (rigorously defined in terms of jet bundles
of the head's configuration manifold), as follows:
When something hits the human head, or the head hits some external
body, we have a collision. This is naturally described by the
SE(3)--momentum, which is a nonlinear coupling of 3 linear
Newtonian momenta with 3 angular Eulerian momenta. The tangent to
the SE(3)--momentum, defined by the (absolute) time derivative, is
the SE(3)--force. The second-order tangency is given by the
SE(3)--jolt, which is the tangent to the SE(3)--force, also
defined by the time derivative.} is a sharp and sudden change in
the SE(3)--force acting on brain's mass--inertia distribution
(given by brain's mass and inertia matrices). That is, a
`delta'--change in a 3D force--vector coupled to a 3D
torque--vector, striking the head--shell with the brain immersed
into the cerebrospinal fluid. In other words, the SE(3)--jolt is a
sudden, sharp and discontinues shock in all 6 coupled dimensions
of brain's continuous micro--motion within the cerebrospinal
fluid (Figure \ref{BrainInj}), namely within the three Cartesian
($x,y,z$)--translations and the three corresponding Euler angles
around the Cartesian axes: roll, pitch and yaw. If the SE(3)--jolt
produces a mild shock to the brain (e.g., strong head shake), it
causes mild TBI, with temporary disabled associated sensory-motor
and/or cognitive functions and affecting respiration and movement.
If the SE(3)--jolt produces a hard shock (hitting the head with
external mass), it causes severe TBI, with the total loss of
gesture, speech and movement.

The SE(3)--jolt is the absolute
time--derivative of the covariant force 1--form (or, co-vector
field). The fundamental law of biomechanics is the \emph{covariant
force law}:
$$\text{Force co-vector field}=\text{Mass distribution}\times \text{%
Acceleration vector--field},$$ which is formally written (using the Einstein
summation convention, with indices labelling the three Cartesian
translations and the three corresponding Euler angles):
\begin{equation*}
F_{{\mu}}=m_{{\mu}{\nu}}a^{{\nu}},\qquad ({\mu,\nu}%
=1,...,6)
\end{equation*}
where $F_{{\mu}}$ denotes the 6 covariant components of the external
``pushing''\ SE(3)--force co-vector field, $m_{{\mu}{\nu}}$
represents the 6$\times $6 covariant components of brain's
inertia--metric tensor, while $a^{{\nu}}$ corresponds to the 6
contravariant components of brain's internal SE(3)--acceleration
vector-field.

Now, the covariant (absolute, Bianchi) time-derivative
$\frac{{D}}{dt}(\cdot )$ of the covariant SE(3)--force $F_{{\mu
}}$ defines the corresponding external ``striking" SE(3)--jolt
co-vector field:
\begin{equation}
\frac{{D}}{dt}(F_{{\mu }})=m_{{\mu }{\nu }}\frac{{D}}{dt}(a^{{\nu }})=m_{{%
\mu }{\nu }}\left( \dot{a}^{{\nu }}+\Gamma _{\mu \lambda }^{{\nu }}a^{{\mu }%
}a^{{\lambda }}\right) ,  \label{Bianchi}
\end{equation}%
where ${\frac{{D}}{dt}}{(}a^{{\nu }})$ denotes the 6 contravariant
components of brain's internal SE(3)--jerk vector-field and
overdot ($\dot{~}$) denotes the time derivative. $\Gamma _{\mu
\lambda }^{{\nu }}$ are the Christoffel's symbols of the
Levi--Civita connection for the SE(3)--group, which are zero in
case of pure Cartesian translations and nonzero in case of
rotations as well as in the full--coupling of translations and
rotations.

In the following, we elaborate on the SE(3)--jolt concept (using
vector and tensor methods) and its biophysical TBI consequences in
the form of brain's dislocations and disclinations.

\subsubsection{$SE(3)-$group of brain's micro--motions within the CSF}

The brain and the CSF together exhibit periodic microscopic
translational and rotational motion in a pulsatile fashion to and
from the cranial cavity, in the frequency range of normal heart
rate (with associated periodic squeezing of brain's
ventricles). This micro--motion is mathematically
defined by the Euclidean (gauge) $SE(3)-$group. 

In other words, the gauge $SE(3)-$group of Euclidean micro-motions
of the brain immersed in the cerebrospinal fluid within the
cranial cavity, contains matrices of the form {\small $ \left(
\begin{array}{cc}
{\bf R} & {\bf b} \\
0 & 1
\end{array}
\right), $} where ${\bf b}$ is brain's 3D micro-translation
vector and ${\bf R}$ is brain's 3D rotation matrix, given by
the product ${\bf R}=R_{\varphi }\cdot R_{\psi }\cdot R_{\theta }$
of brain's three Eulerian micro-rotations,
$\text{roll}=R_{\varphi },~\text{pitch}=R_{\psi
},~\text{yaw}=R_{\theta }$,
performed respectively about the $x-$axis by an angle $%
\varphi ,$ about the $y-$axis by an angle $\psi ,$ and about the
$z-$axis by an angle $\theta $,
{\small
\begin{equation*}
R_{\varphi } =\left[
\begin{array}{ccc}
1 & 0 & 0 \\
0 & \cos \varphi & -\sin \varphi \\
0 & \sin \varphi & \cos \varphi%
\end{array}
\right] , ~~ R_{\psi } =\left[
\begin{array}{ccc}
\cos \psi & 0 & \sin \psi \\
0 & 1 & 0 \\
-\sin \psi & 0 & \cos \psi%
\end{array}
\right] , ~~ R_{\theta } =\left[
\begin{array}{ccc}
\cos \theta & -\sin \theta & 0 \\
\sin \theta & \cos \theta & 0 \\
0 & 0 & 1%
\end{array}
\right].
\end{equation*}}

Therefore, brain's natural $SE(3)-$dynamics within the
cerebrospinal fluid is given by the coupling of Newtonian
(translational) and Eulerian (rotational) equations of
micro-motion.

\subsubsection{Brain's natural $SE(3)-$dynamics}

To support our coupled loading--rate hypothesis, we formulate the coupled
Newton--Euler dynamics of brain's micro-motions within the scull's $%
SE(3)-$group of motions. The forced Newton--Euler equations read in vector
(boldface) form
\begin{eqnarray}
\text{Newton} &:&~\mathbf{\dot{p}}~\mathbf{\equiv M\dot{v}=F+p\times \omega }%
,  \label{vecForm} \\
\text{Euler} &:&~\mathbf{\dot{\pi}}~\mathbf{\equiv I\dot{\omega}=T+\pi
\times \omega +p\times v},  \notag
\end{eqnarray}
where $\times $ denotes
the vector cross product,\footnote{%
Recall that the cross product $\mathbf{u\times v}$ of two vectors $\mathbf{u}
$ and $\mathbf{v}$ equals $\mathbf{u\times v}=uv\func{sin}\theta \mathbf{n}$%
, where $\theta $ is the angle between $\mathbf{u}$ and $\mathbf{v}$, while $%
\mathbf{n}$ is a unit vector perpendicular to the plane of $\mathbf{u}$ and $%
\mathbf{v}$ such that $\mathbf{u}$ and $\mathbf{v}$ form a right-handed
system.}
\begin{equation*}
\mathbf{M}\equiv M_{ij}=diag\{m_{1},m_{2},m_{3}\}\qquad \text{and}\qquad
\mathbf{I}\equiv I_{ij}=diag\{I_{1},I_{2},I_{3}\},\qquad(i,j=1,2,3)
\end{equation*}
are brain's (diagonal) mass and inertia matrices,\footnote{%
In reality, mass and inertia matrices ($\mathbf{M,I}$) are not diagonal but
rather full $3\times 3$ positive--definite symmetric matrices with coupled
mass-- and inertia--products. Even more realistic, fully--coupled
mass--inertial properties of a brain immersed in (incompressible,
irrotational and inviscid) cerebrospinal fluid are defined by the single
non-diagonal $6\times 6$ positive--definite symmetric mass--inertia matrix $%
\mathcal{M}_{SE(3)}$, the so-called material metric tensor of the $SE(3)-$%
group, which has all nonzero mass--inertia coupling products. In other
words, the $6\times 6$ matrix $\mathcal{M}_{SE(3)}$ contains: (i) brain's own mass plus the added mass matrix associated with the fluid, (ii)
brain's own inertia plus the added inertia matrix associated with the
potential flow of the fluid, and (iii) all the coupling terms between linear
and angular momenta. However, for simplicity, in this paper we shall
consider only the simple case of two separate diagonal $3\times 3$ matrices (%
$\mathbf{M,I}$).} defining brain's mass--inertia distribution, with
principal inertia moments given in Cartesian coordinates ($x,y,z$) by volume
integrals
\begin{equation*}
I_{1}=\iiint \rho (z^{2}+y^{2})dxdydz,~~I_{2}=\iiint \rho
(x^{2}+z^{2})dxdydz,~~I_{3}=\iiint \rho (x^{2}+y^{2})dxdydz,
\end{equation*}
dependent on brain's density $\rho =\rho (x,y,z)$,
\begin{equation*}
\mathbf{v}\equiv v^{i}=[v_{1},v_{2},v_{3}]^{t}\qquad \text{and\qquad }%
\mathbf{\omega }\equiv {\omega }^{i}=[\omega _{1},\omega _{2},\omega
_{3}]^{t}
\end{equation*}
(where $[~]^{t}$ denotes the vector transpose) are brain's linear and
angular velocity vectors\footnote{%
In reality, $\mathbf{\omega }$ is a $3\times 3$ \emph{attitude matrix}. However, for simplicity, we will stick to the (mostly)
symmetrical translation--rotation vector form.} (that is, column vectors),
\begin{equation*}
\mathbf{F}\equiv F_{i}=[F_{1},F_{2},F_{3}]\qquad \text{and}\qquad \mathbf{T}%
\equiv T_{i}=[T_{1},T_{2},T_{3}]
\end{equation*}
are gravitational and other external force and torque co-vectors (that is,
row vectors) acting on the brain within the scull,
\begin{eqnarray*}
\mathbf{p} &\equiv &p_{i}\equiv \mathbf{Mv}%
=[p_{1},p_{2},p_{3}]=[m_{1}v_{1},m_{2}v_{2},m_{2}v_{2}]\qquad \text{and} \\
\mathbf{\pi } &\equiv &\pi _{i}\equiv \mathbf{I\omega }=[\pi _{1},\pi
_{2},\pi _{3}]=[I_{1}\omega _{1},I_{2}\omega _{2},I_{3}\omega _{3}]
\end{eqnarray*}
are brain's linear and angular momentum co-vectors.

In tensor form, the forced Newton--Euler equations (\ref{vecForm}) read
\begin{eqnarray*}
\dot{p}_{i} &\equiv &M_{ij}\dot{v}^{j}=F_{i}+\varepsilon _{ik}^{j}p_{j}{%
\omega }^{k},\qquad(i,j,k=1,2,3) \\
\dot{\pi}_{i} &\equiv &I_{ij}\dot{\omega}^{j}=T_{i}+\varepsilon _{ik}^{j}\pi
_{j}\omega ^{k}+\varepsilon _{ik}^{j}p_{j}v^{k},
\end{eqnarray*}
where the permutation symbol $\varepsilon _{ik}^{j}$ is\ defined as
\begin{equation*}
\varepsilon _{ik}^{j}=
\begin{cases}
+1 & \text{if }(i,j,k)\text{ is }(1,2,3),(3,1,2)\text{ or }(2,3,1), \\
-1 & \text{if }(i,j,k)\text{ is }(3,2,1),(1,3,2)\text{ or }(2,1,3), \\
0 & \text{otherwise: }i=j\text{ or }j=k\text{ or }k=i.
\end{cases}
\end{equation*}

In scalar form, the forced Newton--Euler equations (\ref{vecForm}) expand as
\begin{eqnarray}
\text{Newton} &:&\left\{
\begin{array}{c}
\dot{p}_{_{1}}={F_{1}}-{m_{3}}{v_{3}}{\omega _{2}}+{m_{2}}{v_{2}}{\omega _{3}%
} \\
\dot{p}_{_{2}}={F_{2}}+{m_{3}}{v_{3}}{\omega _{1}}-{m_{1}}{v_{1}}{\omega _{3}%
} \\
\dot{p}_{_{3}}={F_{3}}-{m_{2}}{v_{2}}{\omega _{1}}+{m_{1}}{v_{1}}{\omega _{2}%
}
\end{array}
\right. ,  \label{scalarForm} \\
\text{Euler} &:&\left\{
\begin{array}{c}
\dot{\pi}_{_{1}}={T_{1}}+({m_{2}}-{m_{3}}){v_{2}}{v_{3}}+({I_{2}}-{I_{3}}){%
\omega _{2}}{\omega _{3}} \\
\dot{\pi}_{_{2}}={T_{2}}+({m_{3}}-{m_{1}}){v_{1}}{v_{3}}+({I_{3}}-{I_{1}}){%
\omega _{1}}{\omega _{3}} \\
\dot{\pi}_{_{3}}={T_{3}}+({m_{1}}-{m_{2}}){v_{1}}{v_{2}}+({I_{1}}-{I_{2}}){%
\omega _{1}}{\omega _{2}}
\end{array}
\right. ,  \notag
\end{eqnarray}
showing brain's individual mass and inertia couplings.

Equations (\ref{vecForm})--(\ref{scalarForm}) can be derived from the
translational + rotational kinetic energy of the brain\footnote{%
In a fully--coupled Newton--Euler brain dynamics, instead of equation (\ref
{Ek}) we would have brain's kinetic energy defined by the inner product:
\begin{equation*}
E_{k}=\frac{1}{2}\left[ \QOVERD( ) {\mathbf{p}}{\mathbf{\pi }}\left|
\mathcal{M}_{SE(3)}\right. \QOVERD( ) {\mathbf{p}}{\mathbf{\pi }}\right] .
\end{equation*}
}
\begin{equation}
E_{k}={\frac{1}{2}}\mathbf{v}^{t}\mathbf{Mv}+{\frac{1}{2}}\mathbf{\omega }%
^{t}\mathbf{I\omega },  \label{Ek}
\end{equation}
or, in tensor form
\begin{equation*}
E={\frac{1}{2}}M_{ij}{v}^{i}{v}^{j}+{\frac{1}{2}}I_{ij}{\omega}%
^{i}{\omega}^{j}.
\end{equation*}

For this we use the \emph{Kirchhoff--Lagrangian equations}:
\begin{eqnarray}
\frac{d}{{dt}}\partial _{\mathbf{v}}E_{k} &=&\partial _{\mathbf{v}%
}E_{k}\times \mathbf{\omega }+\mathbf{F},  \label{Kirch} \\
{\frac{d}{{dt}}}\partial _{\mathbf{\omega }}E_{k} &=&\partial _{\mathbf{%
\omega }}E_{k}\times \mathbf{\omega }+\partial _{\mathbf{v}}E_{k}\times
\mathbf{v}+\mathbf{T},  \notag
\end{eqnarray}
where $\partial _{\mathbf{v}}E_{k}=\frac{\partial E_{k}}{\partial \mathbf{v}}%
,~\partial _{\mathbf{\omega }}E_{k}=\frac{\partial E_{k}}{\partial \mathbf{%
\omega }}$; in tensor form these equations read
\begin{eqnarray*}
\frac{d}{dt}\partial _{v^{i}}E &=&\varepsilon _{ik}^{j}\left( \partial
_{v^{j}}E\right) \omega ^{k}+F_{i}, \\
\frac{d}{dt}\partial _{{\omega }^{i}}E &=&\varepsilon _{ik}^{j}\left(
\partial _{{\omega }^{j}}E\right) {\omega }^{k}+\varepsilon _{ik}^{j}\left(
\partial _{v^{j}}E\right) v^{k}+T_{i}.
\end{eqnarray*}

Using (\ref{Ek})--(\ref{Kirch}), brain's linear and angular momentum
co-vectors are defined as
\begin{equation*}
\mathbf{p}=\partial _{\mathbf{v}}E_{k}{,\qquad \mathbf{\pi }=\partial _{%
\mathbf{\omega }}E_{k},}
\end{equation*}
or, in tensor form
\begin{equation*}
p_{i}=\partial _{v^{i}}E{,\qquad }\pi _{i}=\partial _{{\omega }^{i}}E,
\end{equation*}
with their corresponding time derivatives, in vector form
\begin{equation*}
~\mathbf{\dot{p}}=\frac{d}{dt}\mathbf{p=}\frac{d}{dt}\partial _{\mathbf{v}}E{%
,\qquad \mathbf{\dot{\pi}}=}\frac{d}{dt}\mathbf{\pi =}\frac{d}{dt}\partial _{%
\mathbf{\omega }}E,
\end{equation*}
or, in tensor form
\begin{equation*}
~\dot{p}_{i}=\frac{d}{dt}p_{i}=\frac{d}{dt}\partial _{v^{i}}E{,\qquad \dot{%
\pi}_{i}=}\frac{d}{dt}\pi _{i}=\frac{d}{dt}\partial _{{\omega }^{i}}E,
\end{equation*}
or, in scalar form
\begin{equation*}
\mathbf{\dot{p}}=[\dot{p}_{1},\dot{p}_{2},\dot{p}_{3}]=[m_{1}\dot{v}%
_{1},m_{2}\dot{v}_{2},m_{3}\dot{v}_{3}],\qquad {\mathbf{\dot{\pi}}}=[\dot{\pi%
}_{1},\dot{\pi}_{2},\dot{\pi}_{3}]=[I_{1}\dot{\omega}_{1},I_{2}\dot{\omega}%
_{2},I_{3}\dot{\omega}_{3}].
\end{equation*}

While brain's healthy $SE(3)-$dynamics within the cerebrospinal fluid is
given by the coupled Newton--Euler micro--dynamics, the TBI is actually
caused by the sharp and discontinuous change in this natural $SE(3)$
micro-dynamics, in the form of the $SE(3)-$jolt, causing brain's
discontinuous deformations.

\subsubsection{Brain's traumatic dynamics: the $SE(3)-$jolt}

The $SE(3)-$jolt, the actual cause of the TBI (in the form of the brain's
plastic deformations), is defined as a coupled Newton+Euler jolt; in
(co)vector form the $SE(3)-$jolt reads\footnote{%
Note that the derivative of the cross--product of two vectors follows the
standard calculus product--rule: $\frac{d}{dt}(\mathbf{u\times v})=\mathbf{%
\dot{u}\times v+u\times \dot{v}.}$}
\begin{equation*}
SE(3)-\text{jolt}:\left\{
\begin{array}{l}
\text{Newton~jolt}:\mathbf{\dot{F}=\ddot{p}-\dot{p}\times \omega -p\times
\dot{\omega}}~,\qquad \\
\text{Euler~jolt}:\mathbf{\dot{T}=\ddot{\pi}}~\mathbf{-\dot{\pi}\times
\omega -\pi \times \dot{\omega}-\dot{p}\times v-p\times \dot{v}},
\end{array}
\right.
\end{equation*}
where the linear and angular jolt co-vectors are
\begin{equation*}
\mathbf{\dot{F}\equiv M\ddot{v}}=[\dot{F}_{{1}},\dot{F}_{{2}},\dot{F}_{{3}%
}],\qquad \mathbf{\dot{T}\equiv I\ddot{\omega}}=[\dot{T}_{{1}},\dot{T}_{{2}},%
\dot{T}_{{3}}],
\end{equation*}
where
\begin{equation*}
\mathbf{\ddot{v}}=[\ddot{v}_{{1}},\ddot{v}_{{2}},\ddot{v}_{{3}}]^{t},\qquad
\mathbf{\ddot{\omega}}=[\ddot{\omega}_{{1}},\ddot{\omega}_{{2}},\ddot{\omega}%
_{{3}}]^{t},
\end{equation*}
are linear and angular jerk vectors.

In tensor form, the $SE(3)-$jolt reads\footnote{%
In this paragraph the overdots actually denote the absolute
Bianchi (covariant) time-derivative (\ref{Bianchi}), so that the
jolts retain the proper covector character, which would be lost if
ordinary time derivatives are used. However, for the sake of
simplicity and wider readability, we stick to the same overdot
notation.}
\begin{eqnarray*}
~\dot{F}_{i} &=&\ddot{p}_{i}-\varepsilon _{ik}^{j}\dot{p}_{j}{\omega }%
^{k}-\varepsilon _{ik}^{j}p_{j}{\dot{\omega}}^{k}, \qquad(i,j,k=1,2,3) \\
~\dot{T}_{{i}} &=&\ddot{\pi}_{i}~-\varepsilon _{ik}^{j}\dot{\pi}_{j}\omega
^{k}-\varepsilon _{ik}^{j}\pi _{j}{\dot{\omega}}^{k}-\varepsilon _{ik}^{j}%
\dot{p}_{j}v^{k}-\varepsilon _{ik}^{j}p_{j}\dot{v}^{k},
\end{eqnarray*}
in which the linear and angular jolt covectors are defined as
\begin{eqnarray*}
\mathbf{\dot{F}} &\equiv &\dot{F}_{i}=\mathbf{M\ddot{v}}\,\equiv \mathbf{\,}%
M_{ij}\ddot{v}^{j}=[\dot{F}_{1},\dot{F}_{2},\dot{F}_{3}], \\
\mathbf{\dot{T}} &\equiv &\dot{T}_{{i}}=\mathbf{I\ddot{\omega}\equiv \,}%
I_{ij}\ddot{\omega}^{j}=[\dot{T}_{{1}},\dot{T}_{{2}},\dot{T}_{{3}}],
\end{eqnarray*}
where \ $\mathbf{\ddot{v}}=\ddot{v}^{{i}},$ and $\mathbf{\ddot{\omega}}=%
\ddot{\omega}^{{i}}$ are linear and angular jerk vectors.

In scalar form, the $SE(3)-$jolt expands as
\begin{eqnarray*}
\text{Newton~jolt} &:&\left\{
\begin{array}{l}
\dot{F}_{{1}}=\ddot{p}_{1}-m_{{2}}\omega _{{3}}\dot{v}_{{2}}+m_{{3}}\left( {%
\omega }_{{2}}\dot{v}_{{3}}+v_{{3}}\dot{\omega}_{{2}}\right) -m_{{2}}v_{{2}}{%
\dot{\omega}}_{{3}}, \\
\dot{F}_{{2}}=\ddot{p}_{2}+m_{{1}}\omega _{{3}}\dot{v}_{{1}}-m_{{3}}\omega _{%
{1}}\dot{v}_{{3}}-m_{{3}}v_{{3}}\dot{\omega}_{{1}}+m_{{1}}v_{{1}}\dot{\omega}%
_{{3}}, \\
\dot{F}_{{3}}=\ddot{p}_{3}-m_{{1}}\omega _{{2}}\dot{v}_{{1}}+m_{{2}}\omega _{%
{1}}\dot{v}_{{2}}-v_{{2}}\dot{\omega}_{{1}}-m_{{1}}v_{{1}}\dot{\omega}_{{2}},
\end{array}
\right. \\
&& \\
\text{Euler~jolt} &:&\left\{
\begin{array}{l}
\dot{T}_{{1}}=\ddot{\pi}_{1}-(m_{{2}}-m_{{3}})\left( v_{{3}}\dot{v}_{{2}}+v_{%
{2}}\dot{v}_{{3}}\right) -(I_{{2}}-I_{{3}})\left( \omega _{{3}}\dot{\omega}_{%
{2}}+{\omega }_{{2}}{\dot{\omega}}_{{3}}\right) , \\
\dot{T}_{{2}}=\ddot{\pi}_{2}+(m_{{1}}-m_{{3}})\left( v_{{3}}\dot{v}_{{1}}+v_{%
{1}}\dot{v}_{{3}}\right) +(I_{{1}}-I_{{3}})\left( {\omega }_{{3}}{\dot{\omega%
}}_{{1}}+{\omega }_{{1}}{\dot{\omega}}_{{3}}\right) , \\
\dot{T}_{{3}}=\ddot{\pi}_{3}-(m_{{1}}-m_{{2}})\left( v_{{2}}\dot{v}_{{1}}+v_{%
{1}}\dot{v}_{{2}}\right) -(I_{{1}}-I_{{2}})\left( {\omega }_{{2}}{\dot{\omega%
}}_{{1}}+{\omega }_{{1}}{\dot{\omega}}_{{2}}\right).
\end{array}
\right.
\end{eqnarray*}

We remark here that the linear and angular momenta ($\mathbf{p,\pi
}$), forces ($\mathbf{F,T}$) and jolts ($\mathbf{\dot{F},\dot{T}}$) are
co-vectors (row vectors), while the linear and angular velocities ($\mathbf{%
v,\omega }$), accelerations ($\mathbf{\dot{v},\dot{\omega}}$) and jerks ($%
\mathbf{\ddot{v},\ddot{\omega}}$) are vectors (column vectors).
This bio-physically means that the `jerk' vector should not be
confused with the `jolt' co-vector. For example, the `jerk'\ means
shaking the head's own mass--inertia matrices (mainly in the
atlanto--occipital and atlanto--axial joints), while the
`jolt'means actually hitting the head with some external
mass--inertia matrices included in the `hitting'\ SE(3)--jolt, or
hitting some external static/massive body with the head (e.g., the
ground -- gravitational effect, or the wall -- inertial effect).
Consequently, the mass-less `jerk' vector\ represents a (translational+rotational) \textit{%
non-collision effect} that can cause only weaker brain injuries, while the
inertial `jolt'\ co-vector represents a (translational+rotational) \textit{%
collision effect} that can cause hard brain injuries.

\subsubsection{Brain's dislocations and disclinations caused by the $SE(3)-$%
jolt}

Recall from introduction that for mild TBI, the best injury predictor
is considered to be the product of brain's strain and strain rate, which is the standard
isotropic viscoelastic continuum concept. To improve this standard concept,
in this subsection, we consider human brain as a 3D anisotropic
multipolar \emph{Cosserat viscoelastic continuum}, exhibiting
coupled--stress--strain elastic properties. This non-standard
continuum model is suitable for analyzing plastic (irreversible)
deformations and fracture mechanics in multi-layered
materials with microstructure (in which slips and bending of
layers introduces additional degrees of freedom, non-existent in
the standard continuum models.

The $SE(3)-$jolt $(\mathbf{\dot{F},\dot{T}})$ causes two types of
brain's rapid discontinuous deformations:

\begin{enumerate}
\item  The Newton jolt $\mathbf{\dot{F}}$ can cause
micro-translational \emph{dislocations}, or discontinuities in the
Cosserat translations;

\item  The Euler jolt $\mathbf{\dot{T}}$ can cause micro-rotational \emph{%
disclinations}, or discontinuities in the Cosserat rotations.
\end{enumerate}

To precisely define brain's dislocations and disclinations,
caused by the $SE(3)-$jolt $(\mathbf{\dot{F},\dot{T}})$, we first
define the coordinate co-frame, i.e., the set of basis 1--forms
$\{dx^{i}\}$, given in local coordinates
$x^{i}=(x^{1},x^{2},x^{3})=(x,y,z)$, attached to brain's
center-of-mass. Then, in the coordinate co-frame $\{dx^{i}\}$ we
introduce the following set of brain's plastic--deformation--related $%
SE(3)-$based differential $p-$forms\footnote{%
Differential $p-$forms are totally skew-symmetric covariant
tensors, defined using the exterior wedge--product and exterior
derivative. The proper definition of exterior derivative $d$ for a
$p-$form $\beta $ on a smooth
manifold $M$, includes the \textit{Poincar\'{e} lemma}: {$%
d(d\beta )=0$}, and validates the \textit{general Stokes formula}
\[
\int_{\partial M}\beta =\int_{M}d\beta ,
\]
where $M$ is a $p-$dimensional \emph{manifold with a boundary} and
$\partial M$ is its $(p-1)-$dimensional \emph{boundary}, while the
integrals have appropriate dimensions.
\par
A $p-$form $\beta $ is called \emph{closed} if its exterior
derivative is equal to zero,
\[
d\beta =0.
\]
From this condition one can see that the closed form (the
\emph{kernel} of the exterior derivative operator $d$) is
conserved quantity. Therefore, closed $p-$forms possess certain
invariant properties, physically corresponding to the conservation
laws.
\par
A $p-$form $\beta $ that is an exterior derivative of some $(p-1)-$form $%
\alpha $,
\[
\beta =d\alpha ,
\]
is called \emph{exact} (the \emph{image} of the exterior
derivative operator $d$). By \textit{Poincar\'{e} lemma}, exact
forms prove to be closed automatically,
\[
d\beta =d(d\alpha )=0.
\]
This lemma is the foundation of the de Rham cohomology theory}:\newline\\
$~~~~$the \emph{dislocation current }1--form, $\mathbf{J}=J_{i}\,dx^{i};$%
\newline
$~~~~$the \emph{dislocation density }2--form, $\mathbf{\alpha }=\frac{1}{2}%
\alpha _{ij}\,dx^{i}\wedge dx^{j};$\newline
$~~~~$the \emph{disclination current }2--form, $\mathbf{S}=\frac{1}{2}%
S_{ij}\,dx^{i}\wedge dx^{j};$ ~and\newline
$~~~~$the \emph{disclination density }3--form, $\mathbf{Q}=\frac{1}{3!}%
Q_{ijk}\,dx^{i}\wedge dx^{j}\wedge dx^{k}$,

where $\wedge $ denotes the exterior wedge--product. These four $SE(3)-$based differential forms
satisfy the following set of continuity equations:
\begin{eqnarray}
&&\mathbf{\dot{\alpha}}=\mathbf{-dJ-S,}  \label{dis1} \\
&&\mathbf{\dot{Q}}=\mathbf{-dS,}  \label{dis2} \\
&&\mathbf{d\alpha }=\mathbf{Q,}  \label{dis3} \\
&&\mathbf{dQ}=\mathbf{0,}\qquad   \label{dis4}
\end{eqnarray}
where $\mathbf{d}$ denotes the exterior derivative.

In components, the simplest, fourth equation (\ref{dis4}),
representing the \emph{Bianchi identity}, can be rewritten as
\[
\mathbf{dQ}=\partial _{l}Q_{[ijk]}\,dx^{l}\wedge dx^{i}\wedge
dx^{j}\wedge dx^{k}=0,
\]
where $\partial _{i}\equiv\partial /\partial x^{i}$, while $\theta
_{\lbrack ij...]}$ denotes the skew-symmetric part of $\theta
_{ij...}$.

Similarly, the third equation (\ref{dis3}) in components reads
\begin{eqnarray*}
\frac{1}{3!}Q_{ijk}\,dx^{i}\wedge dx^{j}\wedge dx^{k} &=&\partial
_{k}\alpha
_{\lbrack ij]}\,dx^{k}\wedge dx^{i}\wedge dx^{j},\text{\qquad or} \\
Q_{ijk} &=&-6\partial _{k}\alpha _{\lbrack ij]}.
\end{eqnarray*}

The second equation (\ref{dis2}) in components reads
\begin{eqnarray*}
\frac{1}{3!}\dot{Q}_{ijk}\,dx^{i}\wedge dx^{j}\wedge dx^{k}
&=&-\partial
_{k}S_{[ij]}\,dx^{k}\wedge dx^{i}\wedge dx^{j},\text{\qquad or} \\
\dot{Q}_{ijk} &=&6\partial _{k}S_{[ij]}.
\end{eqnarray*}

Finally, the first equation (\ref{dis1}) in components reads
\begin{eqnarray*}
\frac{1}{2}\dot{\alpha}_{ij}\,dx^{i}\wedge dx^{j} &=&(\partial _{j}J_{i}-%
\frac{1}{2}S_{ij})\,dx^{i}\wedge dx^{j},\text{\qquad or} \\
\dot{\alpha}_{ij}\, &=&2\partial _{j}J_{i}-S_{ij}\,.
\end{eqnarray*}

In words, we have:

\begin{itemize}
\item  The 2--form equation (\ref{dis1}) defines the time derivative $%
\mathbf{\dot{\alpha}=}\frac{1}{2}\dot{\alpha}_{ij}\,dx^{i}\wedge
dx^{j}$ of the dislocation density $\mathbf{\alpha }$ as the
(negative) sum of the
disclination current $\mathbf{S}$ and the curl of the dislocation current $%
\mathbf{J}$.

\item  The 3--form equation (\ref{dis2}) states that the time derivative $%
\mathbf{\dot{Q}=}\frac{1}{3!}\dot{Q}_{ijk}\,dx^{i}\wedge
dx^{j}\wedge dx^{k}$ of the disclination density $\mathbf{Q}$ is
the (negative) divergence of the disclination current
$\mathbf{S}$.

\item  The 3--form equation (\ref{dis3}) defines the disclination density $%
\mathbf{Q}$ as the divergence of the dislocation density
$\mathbf{\alpha }$, that is, $\mathbf{Q}$ is the \emph{exact}
3--form.

\item  The Bianchi identity (\ref{dis4}) follows from equation
(\ref{dis3})
by \textit{Poincar\'{e} lemma} and states that the disclination density $%
\mathbf{Q}$ is conserved quantity, that is, $\mathbf{Q}$ is the
\emph{closed} 3--form. Also, every 4--form in 3D space is zero.
\end{itemize}

From these equations, we can derive two important conclusions:

\begin{enumerate}
\item  Being the derivatives of the
dislocations, brain's disclinations are higher--order tensors,
and thus more complex quantities, which means that they present
a higher risk for the severe TBI than dislocations --- a fact
which \emph{is} supported by the literature (see review of
existing TBI--models given in Introduction of \cite{ivbrain}).

\item  Brain's dislocations and disclinations are mutually
coupled by the underlaying $SE(3)-$group, which means that we
cannot separately analyze translational and rotational TBIs --- a
fact which \emph{is not} supported by the literature.
\end{enumerate}

For more medical details and references, see \cite{ivbrain}.


\begin{thebibliography}{99}
\bibitem{Arnold} Arnold, V.I., Mathematical Methods of Classical Mechanics.
Springer, New York, (1978)

\bibitem{AbrahamMeh} Abraham, R., Marsden, J., Foundations of Mechanics.
Benjamin, Reading, MA, (1978)

\bibitem{Abraham} Abraham, R., Marsden, J., Ratiu, T., Manifolds, Tensor
Analysis and Applications. Springer, New York, (1988)

\bibitem{Marsden} Marsden, J.E., Ratiu, T.S., Introduction to Mechanics and
Symmetry, A Basic Exposition of Classical Mechanical Systems. (2nd ed),
Springer, New York, (1999)

\bibitem{SIAM} Ivancevic, V. Symplectic Rotational Geometry in Human
Biomechanics. SIAM Rev. 46(3), 455--474, (2004)

\bibitem{GaneshSprSml} Ivancevic, V., Ivancevic, T., Human-Like
Biomechanics, A Unified Mathematical Approach to Human Biomechanics and
Humanoid Robotics. Springer, Berlin, (2005)

\bibitem{GaneshWSc} Ivancevic, V., Ivancevic, T., Natural Biodynamics. World
Scientific, Singapore (2006)

\bibitem{GaneshSprBig} Ivancevic, V., Ivancevic, T., Geometrical Dynamics of
Complex Systems, A Unified Modelling Approach to Physics, Control,
Biomechanics, Neurodynamics and Psycho-Socio-Economical Dynamics. Springer,
Dordrecht, (2006)

\bibitem{GaneshADG} Ivancevic, V., Ivancevic, T., Applied Differential
Geometry, A Modern Introduction. World Scientific, Singapore, (2007)

\bibitem{Switzer} Switzer, R.K., Algebraic Topology -- Homology and
Homotopy. (in Classics in Mathematics), Springer, New York, (1975)

\bibitem{de Rham} de Rham, G., Differentiable Manifolds. Springer, Berlin,
(1984)

\bibitem{LangDg} Lang, S., Fundamentals of Differential Geometry. Graduate
Texts in Mathematics, Springer, New York, (1999)

\bibitem{LangMan} Lang, S., Introduction to Differentiable Manifolds (2nd
ed.). Graduate Texts in Mathematics, Springer, New York, (2002)

\bibitem{Dieudonne} Dieudonne, J.A., Foundations of Modern Analysis (in four
volumes). Academic Press, New York, (1969)

\bibitem{Dieudonne2} Dieudonne, J.A., A History of Algebraic and
Differential Topology 1900-1960. Birkh\'{'}auser, Basel, (1988)

\bibitem{Spivak} Spivak, M., Calculus on Manifolds, A Modern Approach to
Classical Theorems of Advanced Calculus. HarperCollins Publishers, (1965)

\bibitem{SpivakDG} Spivak, M., A comprehensive introduction to differential
geometry, Vol.I-V, Publish or Perish Inc., Berkeley, (1970-75)

\bibitem{Bruhat} Choquet-Bruhat, Y., DeWitt-Morete, C., Analysis, Manifolds
and Physics (2nd ed). North-Holland, Amsterdam, (1982)

\bibitem{Choquet} Choquet-Bruhat, Y., DeWitt-Morete, C., Analysis, Manifolds
and Physics, Part II, 92 Applications (rev. ed). North-Holland, Amsterdam,
(2000)

\bibitem{Bott} Bott, R., Tu, L.W., Differential Forms in Algebraic Topology.
Graduate Texts in Mathematics, Springer, New York, (1982)

\bibitem{Chevalley} Chevalley, C., Theory of Lie groups, Princeton Univ.
Press, Princeton, (1946)

\bibitem{Helgason} Helgason, S., Differential Geometry, Lie Groups and
Symmetric Spaces. (2nd ed.) American Mathematical Society, Providence, RI,
(2001)

\bibitem{Gilmore} Gilmore, R., Lie Groups, Lie Algebras and Some of their
Applications (2nd ed.), Dover, (2002)

\bibitem{Fulton} Fulton, W., Harris, J., Representation theory. A first
course, Graduate Texts in Mathematics, Springer, New York, (1991)

\bibitem{Bourbaki} Bourbaki, N., Elements of Mathematics, Lie Groups and Lie
Algebras, Springer, (2002)

\bibitem{Conway} Conway, J.H., Curtis, R.T., Norton, S.P., Parker, R.A.,
Wilson, R.A., Atlas of Finite Groups: Maximal Subgroups and
Ordinary Characters for Simple Groups. Clarendon Press, Oxford, (1985)

\bibitem{Schafer} Schafer, R.D., An Introduction to Nonassociative Algebras.
Dover, New York, (1996)

\bibitem{ivbrain} V.G. Ivancevic, New mechanics of traumatic brain injury,
Cogn. Neurodyn. \textbf{3}:281-293, (2009)\\
\underline{http://www.springerlink.com/content/p27023577564202h/?p=4351a9d0d76a4fd4b45d6720dad056f3\&pi=8}

\bibitem{ivspine} V.G. Ivancevic, New mechanics of spinal injury, IJAM,
\textbf{1}(2): 387--401, (2009)\\
\underline{http://www.worldscinet.com/ijam/01/0102/S1758825109000174.html}

\bibitem{ivgen} V.G. Ivancevic, New mechanics of generic musculo-skeletal
injury, BRL, \textbf{4}(3):273--287, (2009)\\
\underline{http://www.worldscinet.com/brl/04/0403/S1793048009001022.html}
\end{thebibliography}
\end{document}